\documentclass{amsart}
\usepackage{bbm}
\usepackage{stmaryrd }
\usepackage{amsfonts}
\usepackage{amscd}
\usepackage{amssymb}
\usepackage{latexsym}
 \input{xy}
\xyoption{all}

\usepackage[bookmarks,colorlinks,pagebackref]{hyperref}

\theoremstyle{plain}

\newtheorem*{conjectuur*}{Conjecture}

\swapnumbers
\newtheorem{theorem}[subsection]{Theorem}
\newtheorem{corollary}[subsection]{Corollary}
\newtheorem{lemma}[subsection]{Lemma}
\newtheorem{proposition}[subsection]{Proposition}

\theoremstyle{definition}
\newtheorem{definition}[subsection]{Definition}

\newtheorem{example}[subsection]{Example}

\theoremstyle{remark}

\newtheorem{remark}[subsection]{Remark}

\swapnumbers

\newcommand \ul[1]{\seq{#1}\natural }
\newcommand\cp[1]{{#1}\mathstrut_\sharp}
\newcommand\pp[1]{{#1}\mathstrut_\flat}

\newcommand \after{\circ}

\newcommand \ch{characteristic}

\newcommand \homo{homomorphism}

\renewcommand\iff{if and only if}
\newcommand\implication[2]{\eqref{#1}~$\Rightarrow$~\eqref{#2}}
\newcommand\into{\hookrightarrow}
\newcommand \inv[1]{{#1^{-1}}}
\newcommand \inverse[2]{{#1^{-1}(#2)}}
\newcommand \iso{\cong}

\newcommand \los{\L o\'s' Theorem}

\newcommand \nat{\mathbb N}
\newcommand \norm[1]{\left|#1\right|}

\newcommand \op\operatorname

\newcommand \range [2]{#1,\dots,#2}
\newcommand \restrict [2]{\left.#1\right|_{{#2}}}
\newcommand \rij[2]{(#1_1,\dots,#1_{#2})}
\let\sub\subseteq

\newcommand \zet{\mathbb Z}

\newcommand\dan{\to}

\newcommand\en{\wedge}

\newcommand\of{\vee}

\renewcommand\ul[1]{#1_\natural}
\newcommand \gr{Gro\-then\-dieck ring}
\newcommand\tamered[1]{\mathcal {#1}^{\text{tame}}}
\newcommand \viso{virtually isomorphic}
\newcommand\visoism{virtual isomorphism}
\newcommand \tame{tame}
\newcommand \diag{\ul\omega}
\newcommand \anmodel{\mathcal R^{\text{an}}}
\newcommand \anlang{L^{\text{an}}}
\newcommand \Taylor{Taylor}

\newcommand \ofin{o-finitistic}
\newcommand \powarch[1]{\text{Arch}^{\text{pow}}(\mathcal {#1})}
\newcommand \enh[1]{\Delta#1}
\newcommand \real{\mathbb R}

\newcommand \trace[1]{\op{tr}(#1)}
\newcommand \rankgr[1]{{\mathfrak Z(#1)}}
\newcommand \disc[1]{{\mathfrak D(\mathcal #1)}}
\newcommand \discomin[1]{{\mathfrak D^{\text{viso}}(\mathcal #1)}}

\newcommand \ind[1]{\mathcal {#1}_{\text{ind}}}
\newcommand \defext[2]{{#2}^{{\mathcal #1}}}
\newcommand \discrk[2]{\chi_{#1}(#2)}
\newcommand \trunc[2]{ #1_{\downharpoonleft_{#2}}}
\newcommand \rk[2]{\op{rk}_{#1}(#2)}

\newcommand \euler[2]{\chi_{\mathcal #1}(#2)}
\newcommand \eulerm[2]{\mu_{\mathcal #1}(#2)}
\newcommand \isoclass[1]{{\langle #1\rangle}}
\newcommand \overspill[1]{\ul{\real}{\langle #1\rangle}}
\newcommand \class[1]{{[ #1]}}
\newcommand \node[1]{\op{Node}( #1)}
\newcommand \vertcomp[1]{\op{Vert}( #1)}
\newcommand \fr[1]{\op{fr}( #1)}
\newcommand \lef{\mathbbm L}
\newcommand\fiber[3]{{#1}_{{#3}}[#2]}
\newcommand \grotan{{\mathbf {Gr}^{\text{an}}}}
\newcommand \grotth[2]{{\mathbf {Gr}_{#2}(\mathcal {#1})}}

\newcommand \grotomin[2]{{\mathbf {Gr}^{\text{virt}}_{#2}(\mathcal {#1})}}
\newcommand \grot[1]{{\mathbf {Gr}(\mathcal {#1})}}
\newcommand \grots[1]{{\mathbf {Gr}^{\text{s}}(\mathcal {#1})}}
\newcommand\bary[1]{\beta({#1})}
\newcommand\locdef[1]{\mathcal {#1}^{\text{loc}}}
\newcommand\hardy[1]{\mathbf H(\mathcal {#1})}
\newcommand\hardyfin[1]{\mathbf H^{\text{fin}}(\mathcal {#1})}
\newcommand\omin{o-minimalistic}
\newcommand\Omin{O-minimalistic}
\newcommand\Ded{{\text{DCTC}}}
\newcommand\ominth{T^{\text{omin}}}

\newcommand\oo[2]{{\, ]#1,#2[\, }}
\newcommand\co[2]{{\, [#1,#2[\, }}
\newcommand\oc[2]{{\, ]#1,#2]\, }}
\newcommand\cc[2]{{\, [#1,#2]\, }}

\newcommand \tuple[1]{\mathbf {#1}}

\title {O-minimalism}
\author{Hans Schoutens}
\date{\today} 
\thanks{Partially supported by  PSC-CUNY grant \#62247-00 40}

\begin{document}
\begin{abstract}  
An ordered  structure is called \omin\ if it has all the first-order features of an o-minimal structure. We propose a theory, $\Ded$ (Definable Completeness/Type Completeness), that   describes many properties of  \omin\ structures (dimension theory, monotonicity, Hardy structures, quasi-cell decomposition). Failure of cell decomposition leads to the related notion of a \tame\ structure, and we give a criterium for an \omin\ structure to be \tame. To any \omin\ structure, we can associate its \gr, which in the non-o-minimal case is a non-trivial invariant. To study this invariant, we identify a third \omin\ property, the Discrete Pigeonhole Principle, which in turn allows us to define  discretely valued Euler \ch{s}.
\end{abstract}

\maketitle
\setcounter{tocdepth}{1}
\tableofcontents

\section{Introduction}
O-minimality has been studied extensively (see \cite{vdDomin} for some of the literature). It also has been   generalized in many ways (weak o-minimality, quasi-o-minimality, d-minimality, local o-minimality, o-minimal open cores, etc.) These generalizations attempt to bring into the fold certain ordered structures that fail some of the good finiteness properties of o-minimality, but still behave tamely.\footnote{The concept of `tameness' is quite vague and often depends on the particular author's taste, the present author not excluded; see \S\ref{s:tame}.} We offer a different perspective in this paper, where our point of departure is the observation that, in contrast to ultra{power}s, an ultra{product} of o-minimal structures need no longer be o-minimal; let us call it   \emph{ultra-o-minimal}. This leads to two natural questions: (i)  under which conditions on the o-minimal components is an ultra-o-minimal structure again o-minimal? And (ii), what properties do ultra-o-minimal structures have? We will give one answer to (i) by proving a criterion in terms of Euler \ch{s} in Theorem~\ref{T:ominul} below, but we will be mainly concerned with (ii).

Given a language $L$ with an order relation, let  $\ominth:=\ominth(L)$ be  the intersection of  the theories of all o-minimal $L$-structures.  Models of  $\ominth$ will be called \emph{\omin}; they are precisely the elementary substructures of ultra-o-minimal ones.  O-minimalism is in essence a non-standard feature, as any \omin\ expansion  of the reals is already o-minimal (Corollary~\ref{C:realsomin}). In the first half of this paper, we will focus  on two elementary properties, \emph{definable completeness} (=every definable subset has an infimum) and \emph{type completeness}\footnote{This is a slightly stronger version of what is called in the literature \emph{local o-minimality}, but which agrees with it in the case of an expansion of an ordered field.} (=every one-sided type of a  point, including the ones at infinity, is complete). We denote by $\Ded$ these axiom schemes on one-variable definable sets (where the dependence on parameters has to be quantified out to get   sentences in the language $L$). I do not know whether $\Ded$ is equal to $\ominth$, but in \S\ref{s:gr}, I will formulate a third \omin\ (first-order) property, the Discrete Pigeonhole Principle (DPP=any definable, injective map from a discrete set to itself is bijective), which as of yet, I do not know how to derive from $\Ded$. In fact, it is not clear if we can axiomatize o-minimalism by first-order conditions on one-variable formulae only (note that DPP is a priori not of this form).

But even without knowing fully what the theory $\ominth$ of o-minimalism is, we can deduce at least its first-order properties.  To give just an example, suppose $\mathcal M$   expands an ordered field. If it is o-minimal, then every definable map is continuous, and in fact, differentiable, at all but finitely many points. This is not an elementary statement, so we have to decide with which first-order concept we should replace finiteness   (continuity and differentiability are elementary expressible). As we shall argue below,  the right candidate  is  discreteness (boundedness and closedness then follow automatically). So,   for each formula $\varphi(x,y,\tuple z)$ in $n+2$ variables, we can check, for a given $n$-tuple $\tuple c$, whether the set defined by $\varphi(x,y,\tuple c)$ is the graph of a map $f$, and then we can  express that, if $f$ is not   differentiable at $x$, it must be at every other point in some neighborhood of $x$ where $f$ is defined. It follows that in an \omin\ expansion of an ordered field, every definable map is differentiable outside a bounded, closed, discrete subset (whence outside some finite subset on any compact interval).  In fact, this can be proven already within $\Ded$ by mimicking the proof in \cite[Chapt.~7]{vdDomin}.

So this paper investigates properties of \omin\ structures that are modified versions of the corresponding properties of o-minimal structures. Therefore, whereas most papers on generalizing o-minimality are searching for weakenings that would include certain tamely behaving structures, our hands are tied and we have to obey by the properties of o-minimalism. Thus, to the chagrin of some of my esteemed colleagues, this means that we have to discard the structure $(\mathbb Q,<,+,\mathbb Z)$ as it is not \omin, although it is definably complete and locally o-minimal. However, it fails to have the type completeness property at infinity, which forces every discrete set to be bounded.  
We derive most results already from $\Ded$---albeit only in detail for two variables, leaving   higher arities to the reader, commenting on it occasionally---such as the Monotonicity Theorem (Theorem~\ref{T:disctu}), Fiber Dimension Theorem (Corollary~\ref{C:fibdim}), Quasi-Cell Decomposition (Theorem~\ref{T:plane}), Hardy structures on germs at infinity (Theorem~\ref{T:Hardy}), definable simpleness of \omin\ expansions of groups (Theorem~\ref{T:Dedlimgr}), etc. These are the analogues of the o-minimal concepts except that `finite' has to be  replaced  by `discrete' (which, as already remarked, always implies bounded and closed). This program, however, does not always pan out. For instance, while decomposing into cells, we seem to run into infinite disjunctions, leading to the notion of a \emph{quasi-cell}, which is only locally a cell. However, there is a large class of definable subsets, called \emph\tame\ subsets, that have a `definable' cell decomposition, that is to say, loosely speaking, they admit a cell decomposition in `discretely' many cells (see \S\ref{s:tame} for the precise definition). A \emph\tame\ structure is then one in which every definable subset is \tame,   and we show that it is always at least a model of $\Ded$ (at present I have no examples of a non-\omin\ \tame\ structure, nor of a non-\tame\ \omin\ structure, but presumably these are different notions). 
Any \omin\ expansion of a real closed field  by one-variable functions   is \tame, or more generally, by functions having only a discrete set of discontinuities. 

In \S\ref{s:gr}, we study the \gr\ of an \omin\ structure: it is equal to the ring of integers \iff\ the structure is o-minimal (in which case it corresponds to the Euler \ch). In the last section, we give some criteria for an expansion of an \omin\ structure by a set to be again \omin. For discrete subsets, this leads to the notion of an \emph\ofin\ set, that is to say, a set enjoying all first-order properties of an arbitrary finite set in an o-minimal structure. This notion is particularly interesting when it comes to classifying definable subsets up to `virtual' isomorphism, that is to say, definable in some \omin\ expansion; the corresponding \gr\ is then called the \emph{virtual \gr}. However, a priori, the treatment depends on a choice of `context', that is to say, of an ultra-o-minimal elementary extension. Using this technology, we can now associate to each definable, discrete subset of $M$ a (discretely valued) Euler \ch\ defined on its virtual \gr. This allows us to  calculate  explicitly this  virtual \gr\ in the special case of a \tame, \omin\ expansion of an ordered field admitting a power dominant discrete subset (Corollary~\ref{C:gromintame}).

The last section is an application to the study of analytic sets. In the o-minimal context, (sub)analytic sets are normally understood to be given by analytic functions supported on the unit box (often simply called \emph{restricted analytic functions}), as the corresponding structure $\real_{\text{an}}$ is o-minimal, and admits quantifier elimination in an appropriate language by the seminal work of \cite{DvdD}. There is a good reason to restrict to compact support, as the global sine function defines $\zet$, and hence can never be part of an o-minimal expansion. Our approach here is to look at subsets of $\real^k$ that can be uniformly approximated on compact sets by   $\real_{\text{an}}$-definable subsets. More precisely, we call a subset $X\sub \real^k$ a \emph{Taylor} set, if the ultraproduct over all $n$  of the truncations $\trunc Xn:=\{\tuple x\in X| \norm {\tuple x}\leq n\}$ is definable in $\ul\anmodel$, where the  latter structure is obtained as the ultraproduct of the scalings of $\real_{\text{an}}$ by a factor $n$ (that is to say, for each $n$, the expansion  of $\real$ by power series converging on  $ \norm {\tuple x}\leq n$). Any subset definable by a quantifier free formula using convergent power series, whence in particular, any globally analytic variety, is \Taylor. A discrete subset is \Taylor\ \iff\ it is closed, and any such set satisfies the Discrete Pigeonhole Principle with respect to \Taylor\ maps. However, we can now also define sets by analytic parameterization, like the spiral with polar coordinates $R=\exp\theta$, for $\theta>0$ (in contrast, the spiral obtained by allowing   $\theta$   to be negative  as well  is not \Taylor!). We use our \omin\ results---which, as noted above, are a priori non-standard in nature---to give a geometric treatment of the class of \Taylor\ sets: to a \Taylor\ set $X$, we associate an $\ul\anmodel$-definable subset $\pp X$, called its \emph{protopower},   given as the ultraproduct of its truncations. We obtain a good dimension theory, a monotonicity theorem, a (partly conjectural, locally finite)  cell-decomposition, and a corresponding \gr, all indicative of the tameness of the class of \Taylor\ sets, albeit not   first-order.

\subsection*{Notations and conventions}
Definable means definable with parameters, unless stated explicitly otherwise. 
Throughout this paper, $L$ denotes some language containing  a distinguished binary relation symbol $<$ and any $L$-structure $\mathcal M$ will be a  dense linear order without endpoints.    
 We introduce two new symbols  $-\infty$ and $\infty$, and, given an $L$-structure $\mathcal M$, we let $M_\infty:=M\cup\{\pm\infty\}$. When needed, $\tt U$ denotes some  predicate (often unary), and we will write $(\mathcal M,X)$ for the $L(\tt U)$-structure in which $X$ is the interpretation of $\tt U$. 
 
 We will use the following ISO convention for intervals:   \emph{open}  $\oo ab$ (which we always assume to be non-empty, that is to say, $a<b$), \emph{closed}  $\cc ab$ (including the singleton $\{a\}=\cc aa$), \emph{half-open} $\oc ab$ or $\co ab$, and their infinite variants like $\oo{-\infty}a$, $\oc{-\infty}a$, $\oo a{\infty}$, and $\co a\infty$, with $a,b\in M$.   Note that the usage of $\infty$ here is only informal since these are definable subsets in the language without the extra constants $\pm\infty$ by formulae of the form $x<a$, etc.: any interval is  definable (with parameters). The union and the intersection of two non-disjoint intervals are again   intervals.   
 Note that  in $\mathbb Q$ the set of all rational numbers $q$ with $3<q<\pi$ is not an interval, as it is only an infinite conjunction of definable subsets.   Given a subset $Y\sub M$ and a point $b\in M$, we will sometimes use   notations like $Y_{<b}:=Y\cap \oo{-\infty}b$.

When taking ultraproducts, we rarely ever mention the underlying index set or (non-principal) ultrafilter. We use the notation introduced in \cite{SchUlBook}, denoting ultraproducts with a subscript $\natural$. Thus, we write $\ul\nat$, $\ul\zet$, and  $\ul\real$ for the ultrapower of the set of natural numbers $\nat$, integers $\zet$, and reals $\real$ respectively. On occasion we need the (countable) ultraproduct of the diagonal sequence $(n)_n$ in $\ul\nat$, which we denote suggestively by $\diag$.

\section{O-minimality}

 %
%

 
 \begin{definition}[O-minimality]\label{D:omin}
An $L$-structure $\mathcal M$ is called \emph{o-minimal}, if every definable subset $
Y\sub M$ is a finite union of open intervals and points. 
\end{definition} 
 
 In other words, any $L$-formula (with parameters) in one free variable is equivalent with a (quantifier free) formula using only $<$ (and parameters). The main feature of o-minimality is that this dearth (up to equivalence) of one-variable formulae puts severe restrictions on its many-variable formulae.  Real closed fields are o-minimal; for more examples, we refer to the (vast) literature (\cite{vdDomin} is a good place to start).

 \begin{lemma}\label{L:ded}
Given a definable subset $Y\sub M$ in an o-minimal structure $\mathcal M$, the infimum and the supremum of $Y$ exist in $M_\infty$. 
\end{lemma} 
\begin{proof}
Follows from the fact that the endpoints of an interval are its  infimum and  supremum. Recall the convention that the infimum and the supremum of the empty set are respectively $\infty$ and $-\infty$. 
\end{proof} 

\begin{example}\label{E:nonomin}
Because of quantifier elimination, the structure $(\mathbb Q,<)$ is o-minimal. Expanding this to the ordered field $(\mathbb Q, +,-,\cdot;<;0,1)$, however, destroys this. Indeed,  the set of all $x>0$ such that $x^2>2$ has no infimum in $\mathbb Q$ (its infimum in $\real$, of course, is $\sqrt 2$), contradicting Lemma~\ref{L:ded}.
\end{example} 

 We will always view an ordered structure in its \emph{order topology}  in which the basic open subsets are the open (possibly unbounded) intervals. Note that this is always a Hausdorff topology by density. The Cartesian powers $M^n$ are equipped with the product topology, for which the open boxes form a basis, where an \emph{open box} is any product of open intervals. 
Recall that the \emph{interior} $Y^\circ$ of a subset $Y$ in a topological space $M$ consists of all points $x\in Y$ for which there exists an open $U\sub Y$ with $x\in U$; the \emph{exterior} is the interior of the complement $M \setminus Y$; the \emph{closure} $\bar Y$ is the complement of the exterior; the \emph{frontier} $\fr Y$ is the difference $\bar Y\setminus Y$; and the \emph{boundary}  is the difference $\partial Y:=\bar Y\setminus Y^\circ$. In the ordered case, if $Y\sub M$ is definable, say, by a formula $\phi(x)$, then its interior is also definable, namely by the formula $\phi^\circ(x)$  given as the conjunction of $\phi(x)$ and the formula $\exists z_1<x<z_2,  \forall y: z_1<y<z_2 \dan \phi(y)$; and likewise for its closure, exterior, frontier, and boundary. Hence if $\mathcal M$ is o-minimal, then  $Y^\circ$ is a finite union of open intervals, 
 and $\partial Y$ is  finite. 

A definable subset $X\sub M^n$   is called   \emph{definably connected} if it can not be written as a disjoint union of two   open definable subsets of $X$. The image of a definably connected subset under a definable and continuous function is definably connected. Any interval is definably connected. However, this notion does not always agree with the notion of topological connectedness:  the interval $\oo34$ in $\mathbb Q$ is definably connected, but the decomposition in the two disjoint (non-definable) opens $3<x<\pi$ and $\pi<x<4$ shows that it is not topologically connected (in fact, $\mathbb Q$ is totally disconnected as a topological space). Over the reals, by Dedekind completeness, both notions agree. It is an easy exercise that in an o-minimal structure, a subset of $M$ is definably connected \iff\ it is an interval.


\section{The theory $\Ded$}
Given an $L$-formula $\varphi(\tuple z, x)$ in $n+1$ variables, we will associate several formulae to it with the intention to  express certain properties of o-minimality (here the   $\tuple z$ will be thought of as parameters). Each formula comes with a dual formula upon reversing the order (or, equivalently, by taking the negation), but since we will not need most of these dual ones, we will omit them here. First, we express the Dedekind closure property from Lemma~\ref{L:ded}, that is to say, we construct a formula $\varphi^{\text{inf}}(\tuple z)$ 
 expressing for each tuple of parameters $\tuple z$ that $Y:=\varphi(\tuple z,\mathcal M)$ has an infimum in the $L$-structure $\mathcal M$. 
We allow here infinity as a value, so we will have to treat this case separately. Let $\psi(\tuple z)$ be the formula saying that for each $y$ there exists $x<y$ such that $\varphi(\tuple z, x)$, expressing that the infimum of $Y$ is $-\infty$. Then $\varphi^{\text{inf}}$ is the disjunction of $\psi$ and the formula saying that there exists $x$ such that $x\leq y$ for all $y\in Y$ and such that if $t$ is another element satisfying that $t\leq y$ for all $y\in Y$, then $t\leq x$ (in other words, $x=\inf(Y)$). 
By Lemma~\ref{L:ded}, any o-minimal structure satisfies 
\begin{equation}\label{eq:infsup}
(\forall\tuple z)\varphi^{\text{inf}}(\tuple z).\tag{DC}
\end{equation}

Next we express the following left limit property of intervals: for any point there is an open interval to its left  which is either entirely contained or entirely disjoint from the given interval. In terms of characteristic functions of definable subsets of $M$, this says that any such function has a left limit in each point. More precisely, let $\varphi^-(\tuple z, x)$ be the formula expressing that there exists $y<x$ such that either $\oo yx$ is contained in or disjoint from $Y$. 
Again, we must treat infinity separately:  let  $\varphi^{-\infty}(\tuple z)$ express the fact that there exists $y$ such that $\oo{-\infty}y$ is either contained in or disjoint from $Y$. Since these properties hold for any type of interval, and are preserved under finite unions, we showed that in any o-minimal structure, we have 
\begin{equation}\label{eq:event}
(\forall x,\tuple z)\varphi^-(\tuple z, x)\en\varphi^{-\infty}(\tuple z).\tag{TC}
\end{equation}
%
%
%
%

\begin{definition}[O-minimalism]\label{D:ominic}
For a fixed language with a binary order symbol $<$, we define the theory   $\Ded$ as the extension of  the theory of dense linear orders without endpoints  by the two axiom schemes given by \eqref{eq:infsup} (\emph{definable completeness}) and \eqref{eq:event} (\emph{type completeness}),  with    $\varphi$ running over all $L_{n+1}$-formulae. Put differently, in a model $\mathcal M$ of $\Ded$, any definable subset $Y\sub M$ has an infimum and its \ch\ function has a left limit at each point.

More generally, we call any model of  the first-order theory $\ominth$ of the class of o-minimal $L$-structures an \emph\omin\ structure. In other words, an $L$-structure is \omin, if it satisfies every sentence in $L$ that holds in every o-minimal $L$-structure. Since this includes  \eqref{eq:infsup} and \eqref{eq:event}, an \omin\ structure is a model of $\Ded$. The converse is not clear, although we will see below an axiom, \eqref{eq:DPP}, which is \omin, but is currently not known to follow from $\Ded$. 
\end{definition}

\begin{lemma}\label{L:ominred}
A reduct of an \omin\ structure is again \omin.
\end{lemma}
\begin{proof}
Let $L\sub L'$ be languages, let $\mathcal M'$ be an \omin\ $L'$-structure, and let $\mathcal M:=\restrict{\mathcal M'}L$ be its $L$-reduct. To show that $\mathcal M$ is \omin, take a sentence in $\ominth_L$ and let $\mathcal N'$ be any o-minimal $L'$-structure. Since its reduct $\restrict{\mathcal N'}L$ is also o-minimal, $\sigma$ holds in the latter, whence also in $\mathcal N'$ itself. As this holds for all o-minimal $L'$-structures, $\sigma$ also holds in $\mathcal M'$. Since $ \sigma$ only mentions $ L$-symbols, it must therefore already hold in the reduct $\mathcal M$, as we needed to show.
\end{proof}

We will call an ultraproduct of o-minimal $L$-structures an \emph{ultra-o-minimal} structure. Using a well-known elementarity criterion  via ultraproducts, we have:

\begin{corollary}\label{C:ulomin}
An $L$-structure is \omin\ \iff\ it is elementary equivalent with (equivalently, an elementary substructure of) an ultra-o-minimal structure.\qed
\end{corollary}

This produces many examples of \omin\ structures  which fail to be o-minimal.

\begin{example}\label{E:ulomin}
Let $L$ be the language of ordered fields together with a unary predicate $\tt U$ For each $n$, let $\mathcal R_n:=(\real,\{0,1,\dots,n\}$ be the expansion of field $\real$. Since $\{0,1,2,\dots,n\}$  is finite whence definable without the predicate $\tt U$, each $\mathcal R_n$ is   o-minimal, and therefore their ultraproduct $\ul{\mathcal R}$ is \omin\   by Corollary~\ref{C:ulomin}. The set $D:=\tt U(\ul{\mathcal R})$ is discrete but not finite, so $\ul{\mathcal R}$ cannot be o-minimal. Note that $D$ contains $\nat$ and that $\diag$ is its maximum. In fact, $D=(\ul\nat)_{\leq\diag}$.
\end{example}

By \cite[Corollary 1.5]{MillIVT}, the definable completeness axiom~\eqref{eq:infsup} is equivalent with $M$ itself being definably connected, and also with the validity of the \emph{Intermediate Value Theorem} (IVT is the statement that  if $f\colon \cc ab\to M$ is a definable and continuous function, then $f$ assumes each value between $f(a)$ and $f(b)$). Since an ordered field (with no additional structure) is real closed \iff\ it is o-minimal, \iff\ it satisfies IVT for polynomials, we get:

\begin{corollary}\label{C:ominfield}
An \omin\ ordered (pure) field is o-minimal.\qed
\end{corollary}

This, however, tells us nothing about proper expansions of ordered fields. 
We can reformulate axiom~\eqref{eq:event} as follows, justifying its name. Given $Y\sub M$ and $a\in M$, we say that \emph{$a^-$ belongs to $Y$}, if there exists an open interval $\oo ba\sub Y$ (similarly, \emph{$a^+$ belongs to $Y$}, if $\oo ab\sub Y$ for some $b>a$). Thinking of $a^-$   as a partial type (that is to say, consisting of all formulae $b<x<a$  in the variable $x$, where $b$ runs over all elements less than  $a$), if $Y$ is defined by a formula $\varphi$, then   $a^-$ belongs to $Y$ \iff\  any realization of the type of $a^-$  in any elementary extension of $\mathcal M$ satisfies   $\varphi$. Therefore,   \eqref{eq:event} says that if $Y$ is definable, then   $a^-$ belongs either to $Y$ or to $M\setminus Y$, for any $a\in M$ (and the same assertion for $a^+$), or equivalently, that $a^-$ is a complete type. We extend this terminology to the two types $(-\infty)^+$ and $\infty^-$ in the obvious way. Recall that an $L$-structure $\mathcal M$ is called \emph{locally o-minimal}, if for each definable subset $Y\sub M$ and each $a\in Y$, there exists an open interval $I$ containing $a$ such that $Y\cap I$ is a finite union of intervals. If this is the case, then we may shrink $I$ so that $I\cap Y$ is an interval, from which it is now clear that this is equivalent with \eqref{eq:event} at $a$. In other words, we showed the first part of the following result, whereas the second is a well-known consequence whose proof we just reproduce here (for more on local o-minimality, see\cite{ToffVoz,VozThesis}):

\begin{proposition}\label{P:locomin}
An \omin\ structure $\mathcal M$ is locally o-minimal. Let $K\sub M$ be   compact and $Y\sub M$ definable.  If either $K$ is   open or $Y$ is contained in $ K$, then $K\cap Y$  is a finite union of intervals. 
\end{proposition}
\begin{proof}
Given a definable subset $Y\sub M$, by assumption, we can find in the open case, for each $x\in K$, an open interval $I_x\sub K$ such that $Y\cap I_x$ is an interval. Since $K$ is compact and the $I_x$ cover $K$, there exist finitely many points $x_1,\dots,x_n\in K$ such that $K= I_{x_1}\cup\dots\cup I_{x_n}$ and hence $K\cap Y$ is a finite union of intervals. If $K$ is arbitrary, then we cannot arrange for all $I_x$ to be contained in $K$, and so we only get $K\sub  I_{x_1}\cup\dots\cup I_{x_n}$. But since $Y\sub K$, the same conclusion can be drawn.
\end{proof}

The next corollary improves \cite[2.13(3)]{DMS} as it does not assume any underlying field structure.

\begin{corollary}\label{C:realsomin}
If $\mathcal M$ is an \omin\ structure with underlying order that of the reals, then it is o-minimal.
\end{corollary}
\begin{proof}
Identify $M$ with $\real$, and let $Y\sub \real$ be definable. Depending whether $(-\infty)^+$ or $\infty^-$ belong to $Y$ or not, we may assume after possibly removing one or two unbounded intervals that $Y$ is bounded, whence contained in some closed, bounded interval $K:=\cc ab$. Hence $Y=Y\cap K$ is a finite union of intervals by Proposition~\ref{P:locomin}.  
\end{proof}

\begin{remark}\label{R:realsomin}
From the proof it is clear that we have the following more general result: if an \omin\ structure has the Heine-Borel property, meaning that any closed bounded set is compact, then it is o-minimal.
\end{remark}

Note that \eqref{eq:event} is stronger than local o-minimality, since we also have this condition at $\pm \infty$, which seems to the author an omission in the original definition of \cite{ToffVoz}.  Notwithstanding, since local o-minimality is mostly studied in expansions of ordered fields, the condition at $\infty$ then also holds since we can map the latter to zero by taking reciprocals and apply local o-minimality there.  We start with gathering some easy facts on  one-variable definable subsets in an \omin\ structure.

%

\begin{proposition}\label{P:Dedlim}
For     a definable subset   $Y\sub M$ in an \omin\ structure $\mathcal M$ (or, more generally, in any model of $\Ded$), we have: 
\begin{enumerate}
\item\label{i:infbound} The infimum of $Y$ is either infinite or a topological boundary point.
\item\label{i:intint} If $a,b\in \partial Y$ and $\oo ab\cap \partial Y=\emptyset$, then $\oo ab$ is either disjoint from or entirely contained in $Y$. 
\item\label{i:defconn} If $Y$ is definably connected, then it  is an interval.
\item\label{i:dmin} $Y$ either has a non-empty interior or is discrete.
\item \label{i:disc} If $Y$ is discrete, then it has a minimum and a maximum, and it is closed and bounded.
\item\label{i:secbound} The topological boundary $\partial Y$ is discrete, closed, and bounded.
\end{enumerate}
 \end{proposition}
 \begin{proof}
 To prove \eqref{i:infbound}, let $l\in M$ be the infimum of $Y$. By \eqref{eq:event} and the previous remark, $l^-$ either belongs to $Y$ or to $M\setminus Y$. The former case is excluded since $l$ is the infimum of $Y$. In particular, $l$ is not an interior point of $Y$. If $l^+$ does not belong to $Y$, then $l$ is an isolated point of $Y$, and hence belongs to the (topological) boundary. In the remaining case, $l$ lies in the closure of $Y$, since some open interval $\oo lx$ lies inside $Y$. To prove \eqref{i:intint}, suppose there exists $x\in \oo ab\cap Y$. By \eqref{eq:event}, either $x^-$ belongs to $Y$ or to $M\setminus Y$. In the latter case, there exists  $z<x$ such that $\oo zx$ is disjoint from $Y$. However, $x$ is then not an interior point of $Y$, whence must belong to its topological boundary, contradiction. So $x^-$ belongs to $Y$, which means that  the set of all $z\in\oo ax$ such that $\oo zx\sub Y$ is non-empty. The infimum of this set must   be a topological boundary point of $Y$ by  \eqref{i:infbound}, and hence must be equal to $a$, showing that  $\oc ax\sub Y$. Arguing the same with $x^+$, then shows that also $\co xb\sub Y$, as we needed to prove. 

To prove \eqref{i:defconn}, let $Y\sub M$ be definably connected. Let $l$ and $h$ be its respective infimum and supremum (including the case that these are infinite).
The case  $l=h$ is trivial, so assume $l<h$. If there were some $x\in\oo lh$ not belonging to $Y$, then $Y$ would be   the union of the two definable, non-empty, disjoint open subsets  $Y_{<x}$ and $Y_{>x}$, contradiction, Hence, $Y$ is  an interval with endpoints $l$ and $h$.
%
%
To prove \eqref{i:dmin}, assume $Y$ is not discrete. Hence there exists $a\in Y$ which is not isolated, that is to say, such that any open interval containing $a$ has some other point in common with $Y$. If both $a^-$ and $a^+$ belong  to $M\setminus Y$, then there are $x<a<y$ such that $\oo xa,\oo ay$ are disjoint from $Y$, contradicting that $a$ is not isolated. Hence, say, $a^-$ belongs to $Y$ and $Y$ has non-empty interior.  

Assume next that $Y$ is discrete and let $l$ be its infimum (including possibly the case $l=-\infty$). If $l^+$ belongs to $Y$, then $\oo lz\sub Y$ for some $z>l$, contradicting discreteness. So $l^+$ does not belong to $Y$, which forces $l\in Y$. In particular, $l$ is finite, proving  the first part of \eqref{i:disc}, and in particular, that $Y$ is bounded. To show that $Y$ is closed,  suppose it is not. Let $x\notin Y$ be a point in its closure. Since $Y\cup\{ x\}$ is definable but not discrete, it must have interior by \eqref{i:dmin}, so that some open interval $I$ is contained in  $Y\cup\{ x\}$. But then  $I\cap Y=I\setminus\{x\}$ is not discrete, contradiction. To see \eqref{i:secbound}, it suffices by \eqref{i:disc} to show that $\partial Y$ is discrete. Let $b\in \partial Y$. We have to show that $b$ is an isolated point of $\partial Y$, and this will clearly hold if $b$ is an isolated point of either $Y$ or   $M\setminus Y$. In the remaining case, exactly one from among $b^-$ and $b^+$  belongs to $Y$, say,  $b^-$. Hence, there exist $x<b<y$ so that $\oo xb\sub Y$ and $\oo by\cap Y=\emptyset$. Since any point in $\oo xb$ is interior to $Y$ and any point in $\oo by$ exterior to $Y$,  we get $\partial Y\cap\oo xy=\{b\}$, as we needed to prove. 
\end{proof}

\begin{remark}\label{R:discrete}
Although any discrete, definable subset $Y$ of  an \omin\ structure has a maximum, non-definable subsets of $Y$  need not have a maximum. For instance, the set $D\cap \real=\mathbb N$ from Example~\ref{E:ulomin} has no maximum, whence cannot be definable. Nonetheless, we can define in general  a successor function $\sigma_Y$ on $Y$ by letting $\sigma_Y(b)$ be the minimum of (the definable subset) $Y_{>b}$, for any non-maximal $b$ in $Y$. 
\end{remark}

\begin{corollary}\label{C:ddc}
The theory $\Ded$ is equivalent to type completeness \eqref{eq:event} in conjunction with discrete definable completeness, that is to say, the weaker version of \eqref{eq:infsup} which only requires definable  discrete  sets to have an infimum.
\end{corollary}
\begin{proof}
Let $Y\sub M$ be a definable set in a model $\mathcal M$ of the weaker system from the assertion. Inspecting the argument in the proof of \eqref{i:secbound}, type completeness already implies that $\partial Y$ is discrete. Hence $\partial Y$ has an infimum $b$, and it is now not hard to show that $b$ is also the infimum of $Y$.
\end{proof}

Recall the  definition of  \emph{d-minimality} from \cite{Milldmin}: any definable subset $Y\sub M$ is the union of an open subset and finitely many discrete subsets (where the later finiteness is uniform in the parameters). Hence \eqref{i:dmin} and Corollary~\ref{C:ulomin} immediate yield: 

\begin{corollary}\label{C:dmin}
Any   \omin\ structure is d-minimal. \qed 
\end{corollary}

The reals with a unary predicate defining the integer powers of $2$ is d-minimal, but cannot be \omin, since the latter set is discrete but not closed. The following result fails in general in d-minimal structures, but holds immediately in \omin\ structures since the union of finitely many closed discrete subsets is again discrete.

\begin{corollary}\label{C:uniondisc}
In an \omin\ structure $\mathcal M$, any finite union of definable discrete subsets of $M$ is discrete.\qed
\end{corollary}

Using \eqref{i:secbound}, \eqref{i:intint} and Remark~\ref{R:discrete}, we get immediately the following  structure theorem for one-variable definable subsets:

\begin{theorem}\label{T:onevar}
In an \omin\ structure, every definable subset is  the disjoint union of a  closed, bounded, discrete set and (possibly infinitely many) disjoint open intervals. 

Conversely, if every definable subset $Y\sub M$ of an $L$-structure $\mathcal M$ is a disjoint union of open intervals and a single closed, bounded, discrete set, then $\mathcal M$ is  a model of $\Ded$.
\end{theorem}
\begin{proof}
We only need to show the second assertion. Let us show that $M$ is definably connected. Indeed, if $U_1$ and $U_2$ are disjoint definable open sets covering $M$, then this would yield a covering of $M$ by disjoint open intervals. However, considering what the endpoints would be, this can only be the trivial covering, showing that one of the $U_i$ must be empty. By  \cite[Corollary 1.5]{MillIVT}, definable connectedness implies definable completeness. To prove type completeness, let $Y\sub M$ be definable and $b\in \partial Y$ (boundary points are the only points in which it can fail). There is nothing to prove if  $b$ is an isolated point of $Y$ or of $M\setminus Y$, so assume it is not. Decompose $Y=U\cup D$ into definable subsets with $U$ a disjoint union of open intervals and $D$ closed, bounded, and discrete. Let $\oo pq$ and $\oo uv$ be the open interval in $U$ immediately to the left and to the right of $b$ respectively. Since $b$ is not isolated but in the boundary,  it must be equal to exactly one of $q$ or $u$, that is to say,    either   $q=b<u$ or $q<u=b$. Say the latter holds, so that $b^+$ belongs to $Y$. Since $D\cup\{b\}$ is discrete, we can find an open interval $I$ containing $b$ such that $I\cap (D\cup\{b\})=\{b\}$. Shrinking  $I$ if necessary, we can make it disjoint from $\oo pq$, and hence $I_{<b}\cap Y=\emptyset$, showing that $b^-$ belongs to $M\setminus Y$. We need to verify this also at $b=\pm\infty$, where the same argument works in view of the boundedness of $D$.
\end{proof}

%
%

We conclude this section with an entirely topological characterization of o-minimality and o-minimalism, that is to say, without reference to the order relation. 

\begin{theorem}\label{T:top}
Let  $\mathcal M$ be a (densely ordered) $L$-structure, and let $Y\sub M$ range in the following over all non-empty definable subsets.
\begin{enumerate}
\item\label{i:omintop} $\mathcal M$ is  o-minimal \iff\ all $\partial Y$ are non-empty and finite.
\item\label{i:ominismtop} $\mathcal M$ is  \omin\ \iff\ all $\partial Y$ are non-empty, discrete, and contained in some definable subset  which is definably connected but whose complement is not.
\end{enumerate}
\end{theorem}
\begin{proof}
The direct implication in \eqref{i:ominismtop} is clear since $\partial Y$ is discrete and bounded by \eqref{i:disc}, whence contained in some closed interval $U:=\cc ab$, which is clearly definably connected whereas its complement is not. For the converse, to prove \eqref{eq:infsup}, we may assume $Y$ is bounded from below, and then replace $Y$ with its upward closure, that is to say, all $x$ such that $x\geq y$ for some $y\in Y$. By assumption,  $\partial Y$ contains at least one element, which necessarily must be the infimum of $Y$. To prove \eqref{eq:event} at a point $a\in M$, we may replace $Y$ by $Y_{\leq a}$, and assume $Y_{>a}$ is empty. Let $b$ be the supremum of $Y$. If  $b<a$, then $a^-$ must  belong to $M\setminus Y$ and we are done. So assume $a=b$, whence $a\in\partial Y$. By discreteness, there exists an open interval $I$ around $a$ such that $I\cap \partial Y=\{a\}$. Let $u\in I_{<a}$ and assume $u\notin Y$. Hence the infimum $v$ of $Y_{>u}$ must belong to $\partial Y$, and hence by the choice of $I$ be equal to $a$. On the other hand,  since $a$ is the supremum of $Y$, there must be $u<x<a$ such that $x\in Y$, contradicting the previous observation. Hence $\oo ua\sub Y$, as we needed to show. We also have to verify \eqref{eq:event} at $\pm\infty$, and this will follow once we show that $\partial Y$ is bounded. Choose $U$ containing $\partial Y$  such that $U$ is definably connected and $M\setminus U$ is not definably connected. It suffices to   show that $U$ is bounded. However, under definable completeness \eqref{eq:infsup}, any definably connected subset is an interval (see the proof of \eqref{i:defconn} above), and so if the complement of $U$ is not an interval, then $U$ must be a bounded interval.  

The direct implication in \eqref{i:omintop} is also immediate, and for the converse, we know already from what we just proved that $\mathcal M$ is \omin. By Theorem~\ref{T:onevar}, we can write $Y$ as a disjoint union of open intervals and a discrete set $D$. If two intervals in this decomposition would have a common endpoint which also belongs to $D$, then we replace these three objects by their union, which is then again an open interval. After doing this for all possible endpoints, $\partial Y$ now contains all the endpoints of the remaining intervals as well as points in $D$, and therefore both are finite collections, showing that  $Y$ is a finite union of open intervals and points.
\end{proof}

\begin{remark}\label{R:Pod}
We could apply the criterion given by \eqref{i:omintop} in other situations where there is an underlying definable topology (or, at least, with a basis of definable opens). For instance, a field viewed with its Zariski topology satisfies the criterion of \eqref{i:omintop} \iff\ it is minimal. Since any field is compact with respect to its Zariski topology, discrete sets are finite, and hence a field satisfying the condition in \eqref{i:ominismtop} is likewise minimal. One might be tempted to deduce from this Podewski's conjecture: an elementary extension of a minimal field is again minimal. But this is erroneous, because, unlike the ordered case, neither being discrete nor being definably connected is a first-order condition. This stems from the fact that the collection of basic Zariski open subsets does not belong to a single definable family, unlike the ordered case, so that one cannot quantify over basic open sets. 
\end{remark}

\begin{corollary}\label{C:opencore}
An \omin\ structure with o-minimal open core is itself o-minimal.
\end{corollary}
\begin{proof}
Recall that the \emph{open core} (see, for instance, \cite{DMS,MiSpCore}) of an ordered structure $\mathcal M$ is the reduct in which the definable sets are the definable open subsets of $\mathcal M$. Suppose that $\mathcal M$ is \omin\ and its open core is o-minimal. Let $Y\sub M$ be a definable subset. By \eqref{i:ominismtop}, its boundary $\partial Y$ is non-empty and discrete. Since $\partial Y$ is closed, it is definable in the open core, and therefore   finite by o-minimality. The result now follows from \eqref{i:omintop}. 
\end{proof}

This also follows from the description of the one-variable definable subsets in a structure with o-minimal core as the disjoint union of a finite set and of finitely many dense and co-dense subsets of intervals.

\section{Definable maps}
Next we study definable maps, where we call a map  $f\colon Y\sub M^n\to M^k$  \emph{definable} if its   graph $\Gamma(f)\sub M^{n+k}$  is a definable subset. Note that since    its domain $Y\sub M^n$ is the projection of its graph, it too is definable. Similarly, the set $\Gamma^*(f)$ consisting of all $(f(x),x)\in M^{k+n}$ is definable and is called the \emph{reverse graph} of $f$. If $k=n=1$, we speak of \emph{one-variable} maps.

\begin{lemma}\label{L:imdisc}
Let $f\colon Y\to M$ be a definable map in an \omin\ structure  $\mathcal M$. 
\begin{enumerate}
\item\label{i:imdisc} If $Y$ is discrete, then so is its image $f(Y)$.
\item\label{i:fibdisc} If $f(Y)$ and each   fiber of $f$ is   discrete, then so is $Y$.
\end{enumerate}
\end{lemma}
\begin{proof}
Suppose \eqref{i:imdisc} does not hold, so that $Y\sub M$ is discrete but not $f(Y)$. Let $H$ be the (non-empty, definable, discrete) subset of all $x\in Y$ such that $f(Y_{\geq x})$ is non-discrete, and let $h$ be its maximum. 
Since $h$ cannot be the maximal element of $Y\sub M^n$ lest $f(Y_{\geq h})$ be a singleton, we can find its successor $\sigma(h)\in  Y$ by  \eqref{i:succdisc}. But $f(Y_{\geq h})=\{f(h)\}\sqcup f(Y_{\geq \sigma(h)})$, so that neither $f(Y_{\geq \sigma(h)})$  can be discrete by Corollary~\ref{C:uniondisc}, contradicting  maximality.

Assume next that \eqref{i:fibdisc} is false, so that $Y$ is non-discrete, but  $Z:=f(Y)$ and all $\inverse fu$ are discrete. This time, let $H$ be the subset of all $x\in Z$ such that $\inverse f{Z_{\geq x}}$ is non-discrete, and let $h$ be its maximum. Again $h$ must be non-maximal in $Z$, and so admits an immediate successor $\sigma(h)\in Z$. Since both subsets on the right hand side of 
$$
\inverse f{Z_{\geq h}}=\inverse fh\sqcup \inverse f{Z_{\geq \sigma(h)}}
$$
are discrete, the first by assumption and the second by maximality,   so must their union be by Corollary~\ref{C:uniondisc}, contradiction. 
\end{proof}

We will say that  a one-variable function $f\colon Y\to M$ has a \emph{ jump discontinuity} at a point $c$ if the left and right limit of $f$ at $c$ exist, but are different.

\begin{theorem}[Monotonicity]\label{T:disctu}
 The set of discontinuities of a one-variable definable map $f\colon Y\to M$  in an \omin\ structure  $\mathcal M$   is discrete, closed, and bounded, and consists entirely of jump discontinuities.  Moreover, there is a  definable discrete, closed,  bounded subset    $D\sub Y$ so that in between any two consecutive points of $D\cup\{\pm\infty\}$, the map is \emph{monotone}, that is to say, either strictly increasing, strictly decreasing, or constant.
\end{theorem}
\begin{proof}
We start with proving that all discontinuities are jump discontinuities, or equivalently, that $f$ has a left limit in each point $a\in Y$. For each $y<a$, let $w(y)$ be the supremum of $f(\co ya)$ and let $b$ be the infimum of $w(Y_{< a})$. I claim that $b$ is the left limit of $f$ at $a$. To this end, choose $p<b<q$, and we need to show that there is some $x<a$ with $p<f(x)<q$. If $b^+$ does not belong   $w(Y_{< a})$, and therefore belongs to its complement, then $b$ is an isolated point of  $w(Y_{< a})$,  implying that $f$ takes constant value $b$ on some interval $\oo ya$,  so that $b$ is indeed  the left limit at $a$. In the remaining case, we can find $u>b$ such that $\oc bu\sub w(Y_{< a})$. We may choose $u<q$. In particular, $u=w(y)$ for some $y<a$. Since $b$ is strictly less than the supremum $u=w(y)$, we can find $x\in \co ya$ such that $b< f(x)\leq u$, whence $p<f(x)<q$, as required. 

Let $C\sub Y$ be the definable subset given as the union of the interior of all fibers, that is to say, $x\in C$ \iff\  $x$ is an interior point of $\inverse f{f(x)}$. Being an open set, $C$ is a disjoint union of open intervals, and $f$ is constant, whence continuous on each of these open intervals. Every fiber of the restriction of $f$ to $Y\setminus C$  must have empty interior, whence is discrete by \eqref{i:dmin}. So upon replacing $f$ by this restriction, we may reduce to the case that $f$ has discrete fibers. There is nothing to show if $Y$ is discrete, and so without loss of generality, we may assume $Y$ is an open interval. For fixed $a\in Y$, let $L_a$ (respectively, $H_a$) be the set of all $x\in Y$ such that $f(x)<f(a)$ (respectively, $f(a)<f(x)$). Since $Y$ is the disjoint union of $L_a$, $H_a$, and $\inverse f{f(a)}$ with the latter being discrete, $a^-$ must belong  to one of the first two sets by Corollary~\ref{C:uniondisc}, and depending on which is the case, we will denote this symbolically by writing respectively $f(a^-)< f(a)$ or $f(a^-)> f(a)$ (with a similar convention for $a^+$). Let $L_-$ (respectively, $H_-$, $L_+$, and $H_+$) be the set of all $a\in Y$ such that $f(a^-)<f(a)$ (respectively, $f(a^-)>f(a)$, $f(a^+)<f(a)$, and $f(a^+)>f(a)$), so that $Y$ is the disjoint union of these four definable subsets. Let $D$ be the union of the topological boundaries of these four sets, a discrete set by Corollary~\ref{C:uniondisc}. If $b<c$ are consecutive elements in $D$, then $\oo bc$ must belong to exactly one of these four sets by \eqref{i:intint}, say, to $L_-$. It is now easy to see that in that case $f$ is strictly increasing on $\oo bc$. This then settles the last assertion.

Let $S$ be the (definable) subset of all discontinuities of $f$. To prove that $S$ is discrete, we need to show by \eqref{i:dmin} that it has   empty interior, and this will follow if we can show that any open interval $I\sub Y$ contains a point at which $f$ is continuous. By what we just proved, by shrinking $I$ if necessary, we may assume $f$ is monotone on $I$, say, strictly increasing. Note that $f$ is then in particular injective. By Lemma~\ref{L:imdisc}, the image $f(I)$ contains an open interval $J$. Since $f$ is strictly increasing, $\inverse fJ$ is also an open interval, and $f$ restricts to a bijection between $\inverse fJ$ and $J$. We leave it to the reader to verify that any strictly increasing bijection between intervals is continuous. 
\end{proof}

\begin{remark}\label{R:lim}
Given a definable map $f$ and a point $a$, we denote its left and right limit simply by $f^-(a)$ or $f^+(a)$ respectively, even if these values are infinite (to be distinguished from the symbol $f(a^-)$ which occurred above in formulae of the form $f(a^-)<f(a)$). Note that we even have this property at $\pm\infty$, so that we can define $f^+(-\infty)$ and $f^-(\infty)$, which we then simply abbreviate as $f(-\infty)$ and $f(\infty)$ respectively. 
\end{remark}

\begin{corollary}\label{C:ctuconn}
In an \omin\ structure $\mathcal M$, a definable map $f\colon I\to M$ with domain an open interval $I$ is continuous \iff\ its graph is definably connected.
\end{corollary}
\begin{proof}
Let $C$ be the graph of $f$.  
If $f$ is not continuous, then it has a jump discontinuity  at some point $a\in I$ by Theorem~\ref{T:disctu}. Without loss of generality, we may assume $f^-(a)<f^+(a)$. Let $c$ be some element between these two limits and different from $f(a)$. By definition of one-sided limit, there exist $p<a<q$ such that $f(x)<c$ whenever $p<x<a$, and $f(x)>c$ whenever $a<x<q$. Consider the  two open subsets
\begin{align*}
U_-&:=(I_{<a}\times M)\cup (I_{<q}\times M_{<c}) \\
U_+&:=(I_{>a}\times M)\cup (I_{>p}\times M_{>c}).
\end{align*}
It is not hard to check that $C$ is contained in their union but disjoint from their intersection, showing that it is not definably connected. 

Conversely, assume $f$ is continuous but $C$ is not definably connected, so that there exist definable open subsets $U$ and $U'$ whose union contains $C$ but whose intersection is disjoint from $C$. Since the projections $\pi(C\cap U)$ and $\pi(C\cap U')$  onto the first coordinate are definable subsets    partitioning $I$, they must  have a common boundary point $b\in I$ by Proposition~\ref{P:Dedlim}. Since $(b,f(b))$ belongs to either $U$ or $U'$,  say, to $U$, there exists a box $J\times J'\sub U$ containing $(b,f(b))$. By continuity, we may assume $f(J)\sub J'$. This implies that $(x,f(x))\in U$, for all $x\in J$, and hence that $J\sub\pi(C\cap U)$, contradicting that $b$ is a boundary point of the latter.
%
\end{proof}

\begin{remark}\label{R:genctu}
Without proof, we claim that the above results extend to arbitrary dimensions: given a definable map $f\colon X\sub M^n\to M^k$, the set of discontinuities of $f$ is nowhere dense in $X$. For instance (with terminology to be defined below), if $X\sub M^2$ has dimension two, for each $a,b\in M$, let $D_a$ and $E_b$ be some discrete sets, as given by Theorem~\ref{T:disctu}, such that between any two consecutive points the respective maps $y\mapsto f(a,y)$ and $x\mapsto f(x,b)$ are continuous and monotone. Let $D$  and $E$ be the respective union of all $\{a\}\times D_a$ and all $E_b\times\{b\}$. By Corollary~\ref{C:fibdim} and Proposition~\ref{P:nondisc2} below, both $D$ and $E$ are one-dimensional, closed subsets, and hence $X':=X\setminus(D\cup E)$ is open and dense in $X$. It is now not hard to show that $f$ is continuous on $X'$ (see \cite[Chapt.~3, Lemma 2.16]{vdDomin}). 
\end{remark}

\section{Discrete sets}
We start our analysis of multi-variable definable subsets, with a special emphasis on definable subsets of the \emph{plane} $M^2$ (and address the general case through some sporadic remarks). Since projections play an important role, we introduce some notation. Fix $n$ and let $I\sub \{1,\dots,n\}$ of size $\norm I:=e$. We let $\pi_I\colon M^n\to M^e$ be the projection $\rij an\mapsto (a_{i_1},\dots,a_{i_e})$, where $I=\{i_1<i_2<\dots<i_e\}$. When $I$ is a singleton $\{i\}$, we just write $\pi_i$ for the projection onto the $i$-th coordinate.  Given a tuple $\tuple a=\rij ae\in M_e$, the ($I$-)\emph{fiber} of $X$ above $\tuple a$ is the set 
$$
\fiber X{\tuple a}I:=\pi_{I^c}\left(\inverse{\pi_I}{\tuple a}\cap X\right),
$$
 where $I^c$ is the complement of $I$.  In other words, $\fiber X{\tuple a}I$ is the set of all $\tuple b\in M^{n-e}$ such that $\tuple c\in X$, where $\tuple c$ is obtained from $\tuple b$ by inserting $a_{i_k}$ at the $k$-th spot. In case $I$ is of the form $\{1,\dots,e\}$, for some $e$, we omit $I$ from the notation, since the length of the tuple $\tuple a$ then determines the projection, and we refer to it as a \emph{principal projection}, with a similarly nomenclature for fibers. 
 Thus, for example, the \emph{principal fiber} $\fiber Xa{}$   is the set $\fiber Xa1$ of all $n-1$-tuples $\tuple b$ such that $(a,\tuple b)\in X$.  
Recall that by \eqref{i:disc} any definable discrete subset of $M$ is closed and bounded. The same is true in higher dimensions, for which we first prove:

\begin{theorem}\label{T:projdisc}
 A   definable subset  $X\sub M^n$ in an \omin\ structure $\mathcal M$   is  $X$ is discrete \iff\ all projections $\pi_1(X), \dots, \pi_n(X)$ are discrete. 
 \end{theorem}
\begin{proof}
Suppose all projections are discrete and let $\rij an\in X$. Hence we can find open intervals $I_k$, for $k=\range 1n$,  such that $I_k\cap \pi_k(X)=\{a_k\}$. The open box $I_1\times \cdots\times I_n$ then intersects $X$ only  in the point $\rij an$, proving that $X$ is discrete. 
To prove the converse, we will induct on $n$, proving simultaneously   the following three properties for $X\sub M^n$ discrete:
\begin{enumerate}
\item\label{i:projdisc} $\pi_1(X), \dots, \pi_n(X)$ are discrete;
\item\label{i:totord} $X$ with the induced lexicographical ordering has a minimal element;
\item\label{i:succdisc} for this ordering, there exists a definable map $\sigma_X$ of $X$, sending every non-maximal element in $X$ to its immediate successor.
\end{enumerate}
All three properties have been established by Proposition~\ref{P:Dedlim} when $n=1$, so assume they hold for $n-1$. Assume towards a contradiction that $\pi_1(X)$ is not discrete. For each $a\in\pi_1(X)$,
the fiber $\fiber Xa{}$ (that is to say, the set of all $\tuple b\in M^{n-1}$ such that $(a,\tuple b)\in X$), is discrete since $a\times \fiber Xa{}\sub X$. 
By the induction hypothesis for \eqref{i:totord}, in its lexicographical order, $\fiber Xa{}$ has a minimum, denoted $f(a)$, yielding a definable map $f\colon \pi_1(X)\to M^{n-1}$ whose graph lies in $X$.
  By Theorem~\ref{T:disctu}, each composition $\pi_i\after f\colon \pi_1(X)\to M$, for $i=\range1{n-1}$,  is continuous  outside a discrete set. The union of these discrete sets is again discrete by Corollary~\ref{C:uniondisc}, and hence, since $\pi_1(X)$ is assumed non-discrete, there is a common point $a$ at which all $f_i$ are continuous, whence also $f$. By the discreteness of $X$, we can find an open interval $I$ and an open box $U\sub M^{n-1}$   containing respectively $a$ and $f(a)$ such that $(I\times U)\cap X=\{(a,f(a))\}$. By continuity, we can find an open interval $J\sub I$ containing $a$ such that $f(J)\sub U$. However, this means that for any  $u\in J$ different from $a$, we have $f(u)\in U$, whence $(u,f(u))\in (I\times U)\cap X=\{(a,f(a))\}$, forcing $u=a$, contradiction.

To prove \eqref{i:totord}, we now have established that $\pi_1(X)$ is discrete, whence has a minimum $l$. The minimum of $X$ in the lexicographical ordering is then easily seen to be $(l,\op{min}(\fiber Xl{}))$. To define $\sigma_X$, let $\tuple a=(a,\tuple b)\in X$. For $\tuple a\neq \max(X)$,   either  $\tuple b$ is not the maximum of $\fiber Xa{}$ and we set $\sigma(\tuple a):=(a,\tuple b')$ where $\tuple b':=\sigma_{\fiber Xa{}}(\tuple b)$; or otherwise, $a$ is not the maximum of $\pi_1(X)$ and we set $\sigma(\tuple a):=(a',\op{min}(\fiber X{a'}{}))$, where $a':=\sigma_{\pi_1(X)}(a)$. Note that the existence of $a'$ and $\tuple b'$ follow from the induction hypothesis on \eqref{i:succdisc}. We leave it to the reader to verify that $\sigma_X$ has the required properties. 
\end{proof}

\begin{corollary}\label{C:closeddiscrete}
Any definable, discrete subset in an \omin\ structure is closed and bounded.
\end{corollary}
\begin{proof}
Let $X\sub M^n$ be a definable, discrete subset in an \omin\ structure $\mathcal M$. By Theorem~\ref{T:projdisc}, all $\pi_i(X)$ are discrete, whence bounded and closed by \eqref{i:disc}. It is now easy to deduce from this that so is then $X$. 
\end{proof}

\begin{corollary}\label{C:imdisc}
In an \omin\ structure, the image under a definable map of a definable  discrete set is again discrete.
\end{corollary}
\begin{proof}
If $\mathcal M$ is \omin, and  $f\colon X\sub M^n\to M$ is definable with $X$ discrete, then the graph of $f$ is also discrete (as a subset of $M^{n+1}$). Since $f(X)$ is the projection of this graph, it is discrete by Theorem~\ref{T:projdisc}. 
\end{proof}

\begin{corollary}\label{C:discfib}
A   definable subset  $X\sub M^n$ in an \omin\ structure $\mathcal M$   is discrete \iff\ for some (equivalently, for all)  $I\sub\{1,\dots,n\}$, the projection $\pi_I(X)$  as well as each fiber $\fiber X{\tuple a}I$ is discrete.
\end{corollary}
 \begin{proof}
The converse implication is easy, and for the direct, suppose $X$ is discrete. We may assume, after renumbering, that  $I=\{1,\dots,e\}$. Since each fiber $\fiber X{\tuple a}I$ is homeomorphic to the subset $\tuple a\times \fiber X{\tuple a}I\sub X$, for $\tuple a\in M^e$, and since the latter  is discrete, so is the former. On the other hand, each $\pi_i(X)$, for $i\leq e$, is discrete by Theorem~\ref{T:projdisc}, and since these are just the projections of $\pi_I(X)$, the latter is also discrete, by the same theorem. 
\end{proof}

%
%
 Suppose $\mathcal M$ is an \omin\ expansion of an (Abelian, divisible) ordered group (we will discuss this situation in more detail in \S\ref{s:omingroup}). We define the absolute value $\norm a$ as the maximum of  $a$ and $-a$, for any $a\in M$. We call a map $f\colon X\to X$,  for $X\sub M$, \emph{contractive}, if 
\begin{equation}\label{eq:contract}
\norm{f(x)-f(y)}<\norm {x-y},
\end{equation}
 for all $x\neq y\in X$. We say that $f$ is \emph{weakly contractive}, if instead we have only a weak inequality in \eqref{eq:contract}. Recall that a \emph{fixed point}  of $f$ is a point $x\in X$ such that $f(x)=x$. If $f$ is contractive, it can have at most one fixed point.  
  
 \begin{theorem}[Fixed Point Theorem]\label{T:fixpt}
Let $f\colon D\to D$ be a definable map on a discrete, definable subset $D\sub M$ in an \omin\ expansion $\mathcal M$ of an ordered group. If $f$ is contractive, it has a unique fixed point. If $f$ is weakly contractive, then $f^2$ has a fixed point.
\end{theorem}
\begin{proof}
We treat both cases simultaneously. Assume $f$ does not have a fixed point. In particular, $f(l)>l$, where $l$ is the minimum of $D$. Hence the set of $x\in D$ such that $x<f(x)$ is non-empty, whence has a maximal element $u$. Clearly, $u<h$, where $h$ is the maximum of $D$, and hence $u$ has an immediate successor $v:=\sigma_D(u)$ by \eqref{i:succdisc}. By maximality, we must have $f(v)<v$. Hence  $v\leq f(u)$ and $u\leq f(v)$, and therefore
 $v-u\leq\norm{f(u)-f(v)}$, leading to a contradiction in the contractive case with \eqref{eq:contract}, showing that $f$ must have a fixed point, necessarily unique. In the weak contractive case, we must have  an equality in the latter inequality, whence also in the former two, that is to say,  $f(u)=v$ and $f(v)=u$. Hence,  $u$ and $v$ are fixed points of $f^2$.
\end{proof}

\section{Sets with non-empty interior}
Shortly, we will introduce the notion of dimension, and whereas the discrete sets are those with minimal dimension (=zero), the sets with non-empty interior will be those of maximal dimension.   Note that the non-empty definable subsets of $M$ are exactly of one of these two types by \eqref{i:dmin}.  

\begin{proposition}\label{P:nondisc2}
In an \omin\ structure $\mathcal M$, a definable subset $X\sub M^n$ has non-empty interior \iff\ the set of points $a\in M$ such that the fiber $\fiber Xa{}$ has non-empty interior is non-discrete.
\end{proposition}
\begin{proof}
If $X$ has non-empty interior, it contains an open box, and the assertion is clear. For the converse, note that, since we can pick definably the first open interval inside a definable non-discrete subset of $M$ by the properties proven in Proposition~\ref{P:Dedlim}, we may reduce to the case that $\pi(X)$ is an open interval and each fiber $\fiber Xa{}$ for $a\in \pi(X)$ is   an open box, where $\pi\colon M^n\to M$ is the projection onto the first coordinate. The proof for $n>2$ is practically identical to that for $n=2$, and so, for simplicity, we assume $n=2$. Let $l(a)$ and $h(a)$ be respectively the infimum and supremum of $X_a$, so that $l,h\colon \pi(X)\to M_\infty$ are definable maps. The subset of $\pi(X)$  where either function takes an infinite value is definable, whence it or its complement contains an open interval, so that we can either assume that $l$ is either finite everywhere or equal to $-\infty$ everywhere, and a similar dichotomy for $h$. The infinite cases can be treated by a similar argument, so we will only deal with the case that they are both finite (this is a practice we will follow often in our proofs). By Theorem~\ref{T:disctu}, there is a point $a\in\pi(X)$ at which both $l$ and $h$ are continuous. Fix some $c<l(a)<p<q<h(a)<d$,  so that by continuity, we can find $u<a<v$ so that $l(\oo uv)\sub  \oo cp$ and $h(\oo uv)\sub \oo qd$. I claim that $\oo uv\times \oo pq$ is entirely contained in $X$. Indeed, if $u<x<v$ and $p<y<q$, then from $c<l(x)<p<y<q<h(x)<d$, we get $y\in \fiber Xx{}$, that is to say, $(x,y)\in X$.
\end{proof}

By a simple inductive argument, we get the following analogue of Corollary~\ref{C:discfib}:
 
\begin{corollary}\label{C:openfib}
A   definable subset  $X\sub M^n$ in an \omin\ structure $\mathcal M$  has non-empty interior \iff\ for some (equivalently, for all)  $I\sub\{1,\dots,n\}$, the set of points $\tuple a$ for which $\fiber X{\tuple a}I$ has non-empty interior, has non-empty interior.\qed
\end{corollary}

\begin{corollary}\label{C:non-disc2}
In an \omin\ structure $\mathcal M$, a finite  union of definable subsets of $M^n$ has non-empty interior \iff\ one of the subsets has non-empty interior.
\end{corollary}
\begin{proof}
One direction is immediate, and to prove the other we may by induction reduce to the case of two definable subsets $X_1,X_2\sub M^n$ whose union $X:=X_1\cup X_2$ has non-empty interior. We induct on $n$, where the case $n=1$ follows from Corollary~\ref{C:uniondisc} and \eqref{i:dmin}. Let $W\sub M$ be the subset of   all points $a\in M$ for which the fiber $\fiber Xa{}\sub M^{n-1}$ has non-empty interior.  By Proposition~\ref{P:nondisc2}, the interior of $W$ is non-empty. Since $\fiber Xa{}=\fiber {(X_1)}a{}\cup \fiber {(X_2)}a{}$, our induction hypothesis implies that for $a\in W$, at least one of $\fiber {(X_i)}a{}$ has non-empty interior, in which case we put  $a$ in $W_i$.  In particular, $W=W_1\cup W_2$ so that at least one of the $W_i$ has non-empty interior, say, $W_1$. By Proposition~\ref{P:nondisc2}, this then implies that $X_1$ has non-empty interior.
\end{proof}

\section{Planar cells and arcs}
For the remainder of our analysis of multi-variable definable sets, apart from separate remarks, we restrict to subsets of the plane, that is to say, of $M^2$. Given an ordered structure $\mathcal M$ (soon to be assumed also \omin), let us define a $\emph{$2$-cell}$ in $M^2$ as a definable subset $C$ of the following form: suppose $I$ is an open interval, called the \emph{domain} of the  cell, and $f,g\colon I\to M$ are definable, continuous maps such that $f<g$ (meaning that $f(x)<g(x)$ for all $x\in I$). Let $C$ be the subset of all $(x,y)\in M^2$ with $x\in I$ and $f(x)\diamond_1 y\diamond_2g(x)$, where $\diamond_i$ is either no condition or a strict inequality (when we only have at most one inequality, we get an example of an \emph{unbounded} cell; the remaining ones are call \emph{bounded}, and in arguments we often only treat the latter case and leave the former with almost identical arguments to the reader).  Any $2$-cell is open.   By a \emph{$1$-cell} $C\sub M^2$, we mean either the graph of a continuous definable map $f$  with domain an open interval $I$, or  a Cartesian product $x\times I$. We call the former \emph{horizontal} and the latter \emph{vertical}. Finally, by a \emph{$0$-cell}, we mean a point. We may combine all these definitions into a single definition: a  cell $C$ is determined by elements $a<b$ and definable, continuous maps $f<g\colon M\to M$, as the set of all pairs $(x,y)$ such that $a\diamond_{1} x\diamond_{2} b$ and $f(x)\diamond_3y\diamond_4g(x)$, where each $\diamond_i$ is either no condition, equality or strict inequality. Moreover, if $C$ is non-empty, then  it is a $d$-cell, where $d$ is equal to two minus the number of equality signs among the $\diamond_i$. We sometimes use some suggestive notation like $C(I;f<g)$ to denote, for instance, the  cell given by $x\in I$ and $f(x)<y<g(x)$. If $C$ is a  cell with domain $I$ and $J\sub I$ is an open interval, then we call $C\cap (J\times M)$ the \emph{restriction} of $C$ to $J$. Any restriction of a cell is again a cell, and so is  any principal projection.

\begin{remark}\label{R:gencell}
For higher arity, we likewise define  cells inductively: we say that $C\sub M^n$ is a \emph{$d$-cell} if either $C$ is the graph of a definable, continuous function with domain some $d$-cell in $M^{n-1}$, or otherwise, is the region strictly between two such graphs with common domain some $(d-1)$-cell in $M^{n-1}$. As we shall see below in Remark~\ref{R:dimdef}, the $d$ in $d$-cell refers to the dimension of the cell. 
\end{remark}

\subsection{Arcs}
Assume again that we work in an \omin\ structure $\mathcal M$. 
Given a definable subset $X\sub M^n$, a point $P=(a,b)\in M^2$, and a definable map $h\colon Y\to M$ such that $a\in Y$ and $h^-(a)=b$, we will say that \emph{$P^-_h$ belongs to $X$}, if there exists  an open interval $\oo ua\sub I$  so that the graph of the restriction of $h$  lies inside $X$. By Theorem~\ref{T:disctu}, we may shrink $\oo ua$ so that $h$ is continuous on that interval, and so we could as well view this as a property of the horizontal $1$-cell $C$ defined by $h$. Note that $P$ lies in the closure of $C$. Moreover, we only need $a$ to lie in the closure of $Y$ to make this work. So, given a $1$-cell $C$ such that $P$ lies in its closure, we say that \emph{$P_C^-$ belongs to $X$} if $P_h^-$ does, where $h$ is the definable, continuous map determining $C$, in case $C$ is a horizontal cell, or if $b^-$ belongs to $\fiber Xa{}$ in case $C$ is a vertical cell. Of course, we can make a similar definition for $P_h^+$ or $P_C^+$. The following result essentially shows that viewed as a type, $P_C^-$ is complete:

\begin{lemma}\label{L:limcell}
Given a definable subset $X\sub M^2$ in an \omin\ structure $\mathcal M$, a $1$-cell $C\sub M^2$, and a point $P$ in the closure of $C$, either $P_C^-$ belongs to $X$ or it belongs to its complement.
\end{lemma}
\begin{proof}
Let $P=(a,b)$ be in the closure of $C$. If $C$ is a vertical cell, then $P_C^-$ belongs to $X$ \iff\ $b^-$ belongs to $\fiber Xa{}$, and so we are done in this case by \eqref{eq:event}. In the horizontal case, there exists  a definable, continuous map $h\colon \oo ua\to M$ whose graph is contained in $C$. By \eqref{eq:event}, either $a^-$ belongs to $\pi(X\cap C)$ or to its complement. In the former case, after increasing $u$   if necessary, we have $\oo ua\sub \pi(X\cap C)$, whence $(x,h(x))\in X$ for every $x\in \oo ua$. In the latter case, $\oo ua$ is disjoint from $\pi(X\cap C)$, and hence $(x,h(x))\notin X$ for every $x\in \oo ua$.  
\end{proof}

Using this, it is not hard to show that   the following is an equivalence relation (and in particular symmetric): given a point $P\in M^2$ and $1$-cells $V,W\sub M^2$ such that $P$ lies in each of their closures, we say that $V\equiv_{P^-}W$, if $P_V^-$ belongs to $W$. By a \emph{left arc} at $P$, we mean an $(\equiv_{P^-})$-equivalence class of $1$-cells whose closure contains $P$; and a similar definition for $V\equiv_{P^+}W$ and \emph{right arc}. It is now easy to see that $P_V^-$ belongs to some definable subset $X\sub M^2$  \iff\ $P_W^-$ belongs to it, for any $W\equiv_{P^-}V$, so that we may make sense of the expression \emph{$P_\alpha$ belongs to $X$}, for any left (or right) arc $\alpha$ at $P$. There are two unique equivalence classes containing a vertical cell, called respectively the \emph{lower} and \emph{upper vertical arc}; the remaining ones are called \emph{horizontal}. Given two left horizontal arcs $\alpha$ and $\beta$ at $P$, we can find a common domain $I=\oo ua$ and definable continuous functions $f$ and $g$ on $I$, such that $\alpha$ and $\beta$ are the respective equivalence classes of the graphs of $f$ and $g$. Let $I_-$, $I_=$ and $I_+$ be the subsets of all $x\in I$ such that  $f(x)$ is less than, equal to, or bigger than $g(x)$ respectively. If $\alpha\neq \beta$, then $a$ cannot be in the closure of $I_=$, so that upon shrinking even further, we may assume $I_=$ is empty. Hence $a^-$ belongs either to $I_-$ or $I_+$ and we express this by saying that $\alpha<_{P^-}\beta$ and $\alpha>_{P^-}\beta$ respectively. This yields a well-defined total order relation $<_{P^-}$ on the left horizontal arcs at a point $P$. To include the vertical arcs, we declare the lower one to be  smaller than any horizontal left arc and the upper one to be bigger than any.  

\begin{proposition}\label{P:infarc}
Let $X\sub M^2$ be a definable subset, and $P\in M^2$ a point in an \omin\ structure $\mathcal M$. The set of all left arcs $\alpha$ at $P$  such that $P_\alpha$ belongs to $X$ has an infimum $\beta$ (with respect to the order $<_{P^-}$). If $\beta$ is not vertical, then $P_\beta$ belongs to $\partial X$. 
\end{proposition}
\begin{proof}
Since a point is either interior, exterior or a boundary point, we may upon replacing $X$ by its complement, reduce to the case that $P=(a,b)$ is either interior or a boundary point. In the former case,  the lower vertical arc is clearly minimal,  so assume $P\in\partial X$. In what follows, $\alpha$ always denotes a left arc  at $P$. Consider the set $L_\emptyset$ of all $x<a$ such that $\fiber Xx{}\cap J$ is empty for some open interval $J$ containing $b$. If $a^-$ belongs to $L_\emptyset$, then no $P_\alpha$ belongs to $X$ so that the upper vertical arc is the minimum. So we may assume that the $\fiber Xx{}\cap J$ are non-empty for $x$ close to $a$ from the left. If $b^-$ belongs to $\fiber Xa{}$, then the lower vertical arc is the infimum, so assume $b^-$ belongs to $M\setminus{\fiber Xa{}}$. Hence we may shrink $J$ so that $J\cap (\fiber Xa{})_{<b}$ is empty. For each $x<a$, let $f(x)$ be the infimum of $\fiber Xx{}\cap J$. On a sufficiently small open interval $\oo ua$, the function $f$ is continuous, whence defines a $1$-cell $V$. Since $J\cap (\fiber Xa{})_{<b}=\emptyset$, the left limit $f^-(a)$ must be equal to $b$, showing that  $(a,b)$ lies in the closure of $V$, and hence the equivalence class of $V$ at $P^-$ is a left arc $\beta$. It is now easy to show that $\beta$ is the required infimum, and that it is contained in the boundary $\partial X$.
\end{proof}

\begin{corollary}\label{C:dimcl}
In an \omin\ structure $\mathcal M$, if   $C\sub M^2$  is a definable subset without interior, then so is its closure, that is to say, $C$ is nowhere dense.
\end{corollary}
\begin{proof}
Suppose $P=(a,b)$ is an interior point of the closure $\bar C$, so that there exists an open box $U\sub \bar C$ containing $P$. By Proposition~\ref{P:nondisc2}, the fibers $\fiber Cx{}$ for $x$ close to $a$ must be discrete. By Proposition~\ref{P:infarc},  the infimum $\alpha$ of all  left arcs at $P$ belonging to $C\setminus \fiber Cx{}$ exists, and by discreteness of the surrounding fibers, it must be a minimum, whence also belong to $C$. Similarly, the infimum $\beta$ of all left arcs at $P$ belonging to  $C$ and strictly bigger than $\alpha$ is also a minimum.  Choose an open interval $\oo ua$ such that $\alpha$ and $\beta$ are represented by the respective continuous, definable maps $f,g\colon \oo ua\to M$. Enlarging $u$ if needed, we may assume $f<g$, so that the $2$-cell $S:=C(\oo ua;f<g)$ is disjoint from $C$. 
 Since $S$ is open and $P$ lies in its closure,   $S\cap U$ is non-empty. Since $(S\cap U)\cap C=\emptyset$, no point of $S\cap U$ can lie in the closure of $C$, contradiction.
\end{proof}

%
%

\subsection{The Hardy structure of an \omin\ structure}
We now extend this to infinity in the obvious way: given two horizontal cells $V$ and $W$ with domain an  interval unbounded to the right, we say that $V\equiv_\infty W$ if their restrictions to some interval $\oo u\infty$ are equal. Let $\hardy M$ be the set of all  equivalence classes of cells defined on an open interval unbounded to the right. Note that any definable map $f\colon Y\to M$ whose domain is unbounded to the right yields an equivalence class in $\hardy  M$, denoted $[f]$, since $f$ is continuous by Theorem~\ref{T:disctu} on some open interval $\oo u\infty\sub Y$. Given a definable subset $X\sub M^2$, we can say,  as before, that \emph{$\infty_\alpha$ belongs to $X$}, if  $\infty^-$ belongs to the set of all $x\in Y$ such that $(x,f(x))\in X$, for some $f$ with arc $\alpha$. However, in this case we can do more and make $\hardy  M$ into an $L$-structure: if $\underline c$ is a constant symbol, then we interpret it in $\hardy M$ as the class of the constant function with value $c:=\underline c^{\mathcal M}$; if $\underline F$ is an $n$-ary function symbol, and $\alpha_1,\dots,\alpha_n\in \hardy M$, then $\underline F(\alpha_1,\dots,\alpha_n)$ is the class given by the definable map $\underline F(g_1,\dots,g_n)$, where the $g_i$ are definable functions with domain $I:=\oo u\infty$  such that $[g_i]=\alpha_i$; if $\underline R$ is an $n$-ary predicate symbol, then $\underline R(\alpha_1,\dots,\alpha_n)$ holds in $\hardy M$ \iff\  $\infty^-$ belongs to the set of all $x\in I$ such that $\underline R(g_1(x),\dots,g_n(x))$ holds in $\mathcal M$. 

\begin{definition}\label{D:Hardy}
We call this $L$-structure on $\hardy M$ the \emph{Hardy structure} of $\mathcal M$. In particular, by the same argument as above, $<$ interprets a total order on $\hardy M$, making it into a densely ordered structure without endpoints (note that the notion of vertical arc makes no sense in this context).
\end{definition}

 By induction on the complexity of formulae, we easily can show:

\begin{lemma}\label{L:memberHardy}
Given an \omin\ structure $\mathcal M$, let $\varphi\rij xn$ be a formula with parameters from $M$ and let $X\sub M^n$ be the set  defined by it. For given arcs $\alpha_1,\dots,\alpha_n\in \hardy M$, we have $\hardy M\models \varphi\rij\alpha n$ \iff\  there is a $u\in M$ such that $(g_1(x),\dots,g_n(x))\in X$, for all $x>u$, where each $g_i$ is some continuous function  defined on $\oo u\infty$ representing the arc  $\alpha_i$.\qed
\end{lemma}

Since a continuous function with values in a discrete set must be constant,   Lemma~\ref{L:memberHardy} yields:

\begin{corollary}\label{C:discHardy}
Given an \omin\ structure $\mathcal M$, if a discrete subset  $D\sub (\hardy M)^n$ is definable with parameters in $\mathcal M$, then $D\sub M^n$.\qed
\end{corollary}

\begin{theorem}\label{T:Hardy}
Given an \omin\ structure $\mathcal M$, there is a canonical elementary embedding $\mathcal M\to \hardy M$. In particular, $\hardy M$  is also \omin.
\end{theorem}
\begin{proof}
The map $M\to \hardy M$ sending an element $a\in M$ to the class of the corresponding constant function $x\mapsto a$, is easily seen to be an elementary embedding. 
\end{proof}

These two results together show that if $\mathcal M$ is \omin\ but not o-minimal, then    $(\mathcal M,\hardy M)$ is a Vaughtian pair (see, for instance, \cite[Proposition 9.3]{PillIntro}). In particular,  o-minimalism has Vaughtian pairs.

\begin{remark}\label{R:proto}
Note that we only really used the properties of $\Ded$, so that in the above, we may replace \omin\ by this weaker condition. 
We can think of $\hardy M$ as a sort of \emph{protoproduct}, in the meaning  of a ``controlled'' subring of an ultraproduct as studied in \cite[Chapter 9]{SchUlBook}. Namely, endowing the set $M$ with an ultrafilter containing all right unbounded open intervals, then $\hardy M$ consists of all elements in the ultrapower $\ul{\mathcal M}$  given by definable maps (whereas an arbitrary element is given by any map).

We also can define a \emph{standard part} operator, at least on the  subset $\hardyfin  M$ of all \emph{finite} elements, that is to say, the set of all arcs $\alpha$ at infinity represented by some definable, continuous  map $f\colon \oo u\infty \to M$ such that $f(\infty)\in M$ (see Remark~\ref{R:lim} for the definition). Indeed, the value of $f(\infty)$ only depends on $\alpha$, thus yielding a standard part map   $\hardyfin M\to M$. Note, however, that as $\mathcal M$ is not definable in $\hardy M$, neither is $\hardyfin M$. 
\end{remark} 

\section{Planar curves}
In this section, we fix, without further notice, an \omin\ structure $\mathcal M$. 

\subsection{Dimension}
Let us say that a non-empty definable subset $X\sub M^2$  has \emph{dimension zero} if it is discrete, and \emph{dimension two}, if it has non-empty interior. In the remaining case, we will put $\op{dim}(X)=1$ and call $X$ a (generalized) \emph{planar curve}.  We will assign to the empty set dimension $-\infty$, in order to make the following formula work (with the usual conventions that $-\infty +n=-\infty$): 

\begin{corollary}\label{C:fibdim}
Given   a definable subset $X\sub M^2$, let $F_e$ be the set of  all $a\in M$ for which the fiber $\fiber Xa{}$ has dimension $e$, for $e=0,1$. Then each $F_e$ is definable and the dimension of $X$ is equal to the maximum of all $e+\op{dim}(F_e)$.  
\end{corollary}
\begin{proof}
Being discrete and having interior are definable properties, whence so is being a planar curve, showing that each $F_e$ is definable. The formula then follows by  inspecting the various cases by means of Corollary~\ref{C:discfib} and Proposition~\ref{P:nondisc2}.
\end{proof}

\begin{remark}\label{R:dimdef}
There are several ways of extending this definition to larger arity, and the usual one is to define the dimension of a definable subset $X\sub M^n$ as the largest $d$ such that the image of $X$ under some projection $\pi\colon M^n\to M^d$ has non-empty interior. It follows that a $d$-cell has dimension $d$. 
\end{remark}

\begin{proposition}\label{P:dimulomin}
In an ultra-o-minimal structure $\mathcal M$, a definable set has dimension $e$ \iff\ it is an ultraproduct of $e$-dimensional definable sets. 
\end{proposition}
\begin{proof}
Suppose $\mathcal M$ is the ultraproduct of  o-minimal structures $\mathcal M_i$, and let $X=\varphi(\mathcal M)$ be a definable subset. By \los, $X$ is the ultraproduct of the definable sets $X_i:=\varphi(\mathcal M_i)$. The result now follows from the definability of dimension: we leave again the general case to the reader, but the planar case is clear from our definitions. Indeed, both being discrete and  having non-empty interior are first-order definable properties, and hence pass through the ultraproduct by \los. For instance, for $e=0$, we have   that  $X$ is discrete \iff\ almost all $X_i$ are, and in the o-minimal context the latter is equivalent with being  finite, that is to say, with having  dimension zero. In other words,  $X$ is discrete \iff\ it is \emph{ultra-finite} in the terminology of \S\ref{s:ominexp} below. 
\end{proof}

We may rephrase the previous result as a trichotomy theorem for planar subsets:

\begin{theorem}[Planar Trichotomy]\label{T:triplanar}
In an \omin\ structure $\mathcal M$, a definable planar subset $X\sub M^2$ satisfies exactly one of the following three conditions:
\begin{enumerate}
\item\label{i:0} $X$ is discrete, closed, and bounded;
\item\label{i:1} $X$  is nowhere dense, but at least one projection onto a coordinate axis has non-empty interior;
\item\label{i:2} $X$ has non-empty interior.
\end{enumerate}
\end{theorem}
\begin{proof}
We only need to show that \eqref{i:1} \iff\ $X$ is one-dimensional. The converse  is clear from Corollaries~\ref{C:fibdim} and \ref{C:dimcl}, and for the direct implication, we must show that a  definable subset   satisfying \eqref{i:1} cannot be discrete, and this follows from Theorem~\ref{T:projdisc}.
\end{proof}

Immediate from the definitions and Corollary~\ref{C:non-disc2}, we have:
\begin{corollary}\label{C:dimunion}
The dimension of a union $X_1\cup\dots \cup X_n\sub M^2$ of   definable subsets   is the maximum of the dimensions of the $X_i$.\qed
\end{corollary}

\subsection{Nodes}
Assume $C$ is a planar curve.  Let us call a point $P\in C$ a \emph{node}, if for every open box $B$ containing $P$, there is an open subbox $I\times J\sub B$ containing $P$ and some point $x\in I$ such that $\fiber Cx{}\cap J$ is not a singleton. We denote the set of nodes of $C$ by $\node C$. 
 We call a node an \emph{edge}, if in the above condition $\fiber Cx{}\cap J$ can be made empty. By an argument similar to the one proving Corollary~\ref{C:ctuconn}, one shows that  a function on an open  interval $h$ is continuous \iff\ its graph has no edges (since it is a graph, it cannot have any other type of nodes). Note that the closure of a $1$-cell has at most two edges: indeed, if  the cell $C$ is given as the graph of a definable, continuous function $h$ on an interval $\oo ab$, then $\bar C\setminus C$ consists of those points  among $(a,h^+(a))$ and $(b,h^-(b))$ that are finite (in the notation of Remark~\ref{R:lim}), and these are then the edges of $\bar C$. Isolated points are edges, and they form a discrete, closed, and bounded subset. Another special case of an edge is any point lying on an open interval inside a vertical fiber $\fiber Ca{}$. Let $\vertcomp C$ be the set of all such edges, called the \emph{vertical component} of $C$. Note that $\vertcomp C$ is equal to  the union of the interiors of all fibers, that is to say, $\vertcomp C=\bigcup_a(\fiber Ca{})^\circ$, and hence in particular is definable.   

\begin{proposition}\label{P:discnodes}
The set of nodes of a planar curve is the union of its  vertical  component and a  discrete set. 
\end{proposition}
\begin{proof}
Let $C\sub M^2$ be a planar curve. Replacing $C$ by $C\setminus\vertcomp C$, we may assume its vertical component is empty.   Assume towards a contradiction that  $N:=\node C$ is not discrete. Therefore, $\pi(N)$ cannot be discrete by Corollary~\ref{C:discfib}, and hence contains an open interval $I$. For each $x\in I$, let $h(x)$ be the minimal $y\in \fiber Cx{}$ such that $(x,y)\in N$. By Theorem~\ref{T:disctu}, we may shrink $I$ so that $h$ becomes a continuous function on $I$. In particular,  its graph $V$ is a $1$-cell contained in $N$. For each $x\in I$, let $l(x)$ and $u(x)$ be the  respective predecessor and successor in $\fiber Cx{}$ of $h(x)$ (if $h(x)$ is always an extremal element of $\fiber Cx{}$ then we can adjust the argument accordingly, and so we just assume that $l(x)<h(x)<u(x)$ always exist). Since $(x,h(x))$ is a node and $h$ is continuous, for $y<x$ sufficiently close to $x$, and $J$ an open interval such that $J\cap \fiber Cx{}=\{h(x)\}$, the intersection $J\cap \fiber Cy{}$ contains at least one other element besides $h(y)$, necessarily either $l(y)$ or $u(y)$.  By \eqref{eq:event}, either $l(y)$ belongs to all $J\cap \fiber Cy{}$, for all $y$ sufficiently close to the left of $x$, or otherwise $u(y)$ does. In particular, for a fixed $x\in I$, we have $h(x)=l^-(x)$ or $h(x)=u^-(x)$. Shrinking $I$ if necessary, \eqref{eq:event} then reduces to the case that one of these alternatives happens for every $x\in I$, say, $h(x)=l^-(x)$ for all $x\in I$. Shrinking $I$ even further, we may assume that  $l$   is continuous on $I$, and hence $l=h$ on $I$, contradiction.
\end{proof}

\begin{lemma}\label{L:nonnode}
  A point $P$ on a  planar curve $C$  is not a  node \iff\ there is some open box $B$ containing $P$ such that $C\cap B$ is a horizontal $1$-cell. On the other hand, $P$ is an edge \iff\ it does not belong to  any horizontal cell inside $C$.   
\end{lemma}
\begin{proof}
If $P=(a,b)\notin\node C$, there exist  open intervals $I$ and $J$ containing  respectively $a$ and $b$ such that $\fiber Cx{}\cap J$ is a singleton $\{f(x)\}$, for every $x\in I$, and this property is preserved for any subbox of $I\times J$ containing $P$. Hence $f\colon I\to J$ is a definable map with $f(a)=b$. Shrinking $I$ if necessary, we may assume by Theorem~\ref{T:disctu} that $f$ is continuous on $I$ with a possible exception at $a$. As already observed, $f$ is also continuous at $a$ lest $(a,f(a))$ be a node. Hence the graph of $f$ is a cell equal to $(I\times J)\cap C$. If $P$ is not a node, then by definition, no intersection with an open box around $P$ can be a cell. The second assertion is obvious.
\end{proof}

\begin{remark}\label{R:regular}
This allows us to generalize the notion to higher arity: let us say that a point $P$ on a definable subset $X\sub M^n$ is \emph{strongly $e$-regular}, for some $e\leq n$, if there exists an open box $B$ containing $P$ such that $B\cap X$ is an $e$-cell. When $n=2$, a point is strongly $2$-regular \iff\ it is interior, and strongly $0$-regular \iff\ it is isolated. The previous result then says that on a planar curve, a point is strongly $1$-regular \iff\ it is not a node. As with cells, this definition of regularity has a directional bias: nodes are really critical points with respect to projection onto the first coordinate.
To break this bias, just taking permutations of the variables does not give enough   transformations to turn some  point on a curve in a non-nodal position, as for instance the origin on the curve given by $(t,t)$ if $t\leq 0$ and $(-t,-2t)$ if $t\geq 0$.  However, if we assume that there is an underlying ordered group (see \S\ref{s:omingroup} below), then we say that  a point $x\in X\sub M^n$ is   \emph{$e$-regular}, if after a translation bringing $x$ to the origin $O$, we can find a rotation $\rho$ such that $\rho(O)$ is strongly $e$-regular in $\rho(X)$, where by a \emph{rotation} of $M^n$, we mean a linear map $\rho\colon M^n\to M^n$ given by an invertible matrix of determinant one over $\mathbb Q$ (as we shall see below, any \omin\ expansion of a group is divisible, whence admits a natural structure of a $\mathbb Q$-vector space). 
\end{remark}

\begin{proposition}\label{P:dimcl}
A definable subset $X\sub M^2$ has the same dimension as that of its closure $\bar X$, whereas the dimension of its frontier $\fr X$  is    strictly less.
\end{proposition}
\begin{proof}
If $X$ is discrete, then it is closed by Corollary~\ref{C:closeddiscrete}, and so $\fr X=\emptyset$, proving the assertion in this case. If $X$ has dimension one, then so does $\bar X$ by Corollary~\ref{C:dimcl}. Let $V:=\vertcomp C$ be the vertical component of $C$ and let $\pi(V)$ be its projection. Since $\pi(V)$ is discrete by Proposition~\ref{P:nondisc2}, the boundary   $\partial V$ is equal to the union of all $\partial (\fiber Xa{})$, whence is discrete by Corollary~\ref{C:discfib}. Hence, upon removing $V$ from $X$, we may reduce to the case that $X$ has no vertical components.  Suppose towards a contradiction that  $\fr X$ is a planar curve. By Proposition~\ref{P:discnodes}, the set of nodes on $\bar X$ and on $\fr X$ are both discrete sets, and so, there exists a $P\in \fr X$ which is not a node on $\fr X$ nor on $\bar X$. By Lemma~\ref{L:nonnode}, there exists an open box $B$ containing $P$ such that both $B\cap \fr X$ and  $B\cap \bar X$ are cells, and therefore the inclusion $B\cap \fr X\sub B\cap \bar X$ must be an equality. In particular,  $B\cap X$ is empty, contradicting that $P$ lies in the closure of $X$. 

Finally, if $X$ has dimension two, then so must $\bar X$. Let $Y:=X^\circ$ and $Z:=X\setminus Y$. Since $\bar X=\bar Y\cup\bar Z$, we have $\fr X=(\bar Y\setminus X)\cup (\bar Z\setminus X)$, so that it suffices to show that neither of these two differences has interior. The first one, $\bar Y\setminus X$, is equal to $\partial Y$ whence has no interior, being the boundary of an open set. By construction, $Z$ has no interior, and hence by Corollary~\ref{C:dimcl}, neither does its closure. 
\end{proof}

Recall that a \emph{constructible} subset is a finite Boolean combination of open subsets, and hence every one-variable definable subset in an \omin\ structure is constructible. This is still true in higher dimensions: by an easy induction on the dimension, and using that the closure is obtained by adjoining the frontier, Proposition~\ref{P:dimcl} yields:

\begin{corollary}\label{C:const}
In an \omin\ structure $\mathcal M$, any definable subset  is constructible. \qed
\end{corollary}

\begin{corollary}\label{C:bddim2}
The boundary of a two-dimensional definable subset of $M^2$ has dimension at most one. 
\end{corollary}
\begin{proof}
Let $X\sub M^2$ has dimension two. Its boundary $\partial X$ is the union of its frontier $\fr X$ and $X\setminus X^\circ$. The former has dimension at most one by Proposition~\ref{P:dimcl} and the latter has no interior. The result now follows from Corollary~\ref{C:non-disc2}.
\end{proof}

Recall that a subset in a topological space is called \emph{codense} if its complement is dense.

\begin{corollary}\label{C:codense}
If $Y$ is a codense definable subset of a non-empty definable subset $X\sub M^2$, then $\op{dim}(Y)<\op{dim}(X)$.
\end{corollary}
\begin{proof}
If $X$ is discrete, then it is closed by Corollary~\ref{C:closeddiscrete}, and hence its only codense subset is the empty set. If $X$ and $Y$ both have dimension two, then $Y^\circ$ is disjoint from the closure of $X\setminus Y$, contradicting that $Y$ is codense in $X$. So remains the case that $X$ is a curve. 
If $Y$ is codense in $X$, then it must be contained in the frontier of $X\setminus Y$, and the latter has dimension strictly less than one by Proposition~\ref{P:dimcl}. 
%
\end{proof}

\subsection{Quasi-cells}
To obtain a cell decomposition  as in the o-minimal case, we must generalize the notion of $1$-cell by the following equivalent conditions:

\begin{lemma}\label{L:supercell}
Let $S$ be a union of mutually intersecting $1$-cells in $M^2$. We call $S$ a (horizontal) \emph{$1$-quasi-cell} if it satisfies one of the following equivalent conditions:
\begin{enumerate}
\item\label{i:nonodes} $S$ has no nodes;
\item\label{i:locfct} $S$ is the graph of a continuous map $h\colon \pi(S)\to M$ which is   \emph{locally definable}, meaning that its restriction to any open interval in its domain is definable.
\end{enumerate}
Moreover, $\pi(S)$ is open and convex, and    $S$    is a $1$-cell \iff\  $\pi(S)$ is definable.
\end{lemma}
\begin{proof}
The 
implication  \implication{i:locfct}{i:nonodes} is easy, since the graph of a continuous function has no nodes. To show  \implication{i:nonodes}{i:locfct}, suppose $S$ has no nodes, so that in particular,   no vertical cell lies inside $S$. Fix $a_1,a_2\in\pi(X)$  and choose $1$-cells $C_1$ and $C_2$ containing $a_1$ and $a_2$ respectively. Let $I_k:=\pi(C_k)$ and let $h_k$ be the definable (continuous) function on $I_k$ whose graph is $C_k$. Let $I:=I_1\cup I_2$. Since $C_1\cap C_2$ is non-empty, so is $I_1\cap I_2$, showing that $I$ is an interval. Let $H$ be the subset of $I_1\cap I_2$ on which $h_1$ and $h_2$ agree, that is to say, $H=\pi(C_1\cap C_2)$. For $a\in H$ with common value $b$, since $(a,b)$ is not a node, there exist    open intervals  $U$  and $V$ containing respectively  $a$ and $b$ such that $\fiber Sx{}\cap V$ is a singleton, for all $x\in U$. Shrinking $U$ if necessary, continuity allows us to assume that $h_k(U)\sub V$, so that $(x,h_k(x))$ both lie in $\fiber Sx{}\cap V$ for $x\in U$, whence must be equal. This shows that $U\sub H$, and hence that $H$ is open. Let $a\in \partial H$. Since $H$ is open, $a$ does not belong to $H$ whereas either $a^-$ or $a^+$ does. If $a$ lies in $I_1\cap I_2$, then $a\notin H$ implies $h_1(a)\neq h_2(a)$, and by continuity, this remains the case on some open around $a$, contradicting that either $a^-$ or $a^+$ belongs to $H$. Hence $a\notin I_1\cap I_2$. This means that the only boundary points of $H$ are the endpoints of the interval $I_1\cap I_2$, proving that $h_1$ and $h_2$ agree on this interval. Let $h(x)$ be equal to $h_1(x)$ if $x\in I_1$ and to $h_2(x)$ otherwise. Since the graph of $h$ is then equal to $C_1\cup C_2\sub S$, whence contains no nodes, $h$ is continuous. It is not hard to see that $\pi(S)$ is open and convex. The last assertion then follows since $\pi(S)$ is   a disjoint union of open intervals by Theorem~\ref{T:onevar}, whence, being also convex,  a single open interval. 
\end{proof}

Although we should also entertain the notion of  vertical quasi-cells (see Remark~\ref{R:gensupercell} below), they do not occur in the analysis of planar subsets.  Given a curve $C$ without nodes and a quasi-cell $S\sub C$, we say that $S$ is \emph{optimal} in $C$, if    no quasi-cell inside $C$ strictly contains $S$. 

\begin{corollary}\label{C:supercell}
Any point on a planar curve  $C$ without nodes  lies on a (uniquely determined) optimal quasi-cell  in $C$. In particular, $C$ is a disjoint union of quasi-cells.
\end{corollary}
\begin{proof}
Fix $P\in C$. By Lemma~\ref{L:nonnode}, there exists a cell $V\sub C$ containing $P$.  Let $S$ be the union of all cells inside $C$ containing $P$. Since $S\sub C$,   has no nodes, it  is a quasi-cell by Lemma~\ref{L:supercell}. Suppose $S'\sub C$ is a quasi-cell containing $S$ and let $P'\in S'$. By  Lemma~\ref{L:supercell}, there exists a cell $V'\sub S'$ containing both $P$ and $P'$. By construction, we then have $V'\sub S$, whence $P'\in S$, showing that $S=S'$ is optimal.
\end{proof}

If the curve has nodes, then to preserve uniqueness of optimal quasi-cells, we have to amend this definition as follows: for an arbitrary planar curve, a horizontal $1$-quasi-cell $S$ is called \emph{optimal} if $S\sub C$ contains no node of $C$ and is maximal with this property.  

Finally, we define the notion of a \emph{$2$-quasi-cell} $S\sub M^2$ given as the region between two continuous, locally definable maps defined on an open, convex subset of $M$ (again called the \emph{domain} of the quasi-cell), or an unbounded variant as in the case of $2$-cells. More precisely, let $V\sub M$ be an open convex subset, $f,g\colon V\to M$ continuous and locally definable with $f<g$, then $S$ consist of all $(x,y)$ such that $x\in V$ and $f(x)\diamond y\diamond' g(x)$, with  $\diamond,\diamond'$ strict inequality or no condition. By definition of local definability, the restriction of a quasi-cell $S$ to an   open interval $I\sub V$, that is to say,  $S\cap (I\times M)$ is a cell,   and hence every $2$-quasi-cell is the union of $2$-cells and therefore open. Moreover, we can arrange for   all these cells in this union to contain a  given  fixed point of the quasi-cell. 

\begin{remark}\label{R:gensupercell}
The definition of an arbitrary $d$-quasi-cell is entirely similar: simply replace in the recursive definition from Remark~\ref{R:gencell} `cell' by `quasi-cell' and `definable, continuous map' by `locally definable, continuous' map everywhere.
\end{remark}

\subsection{Locally definable subsets}\label{s:locdef}
Quasi-cells are particular instances of locally definable subsets, which we now briefly study. In an arbitrary ordered structure $\mathcal M$, we call a subset $X\sub M^n$ \emph{locally definable} if for each point $P\in M_\infty^n$, there exists an open box $B$ containing $P$ such that $B\cap X$ is definable. It is important to include  in this definition also the \emph{infinite points} of $M_\infty^n$, that is to say, points with at least one coordinate equal to $\pm\infty$, where,  just as an example,  we mean by an \emph{open box around an infinite point}  like $(0,\infty,-\infty)$ one of the form $\oo uv\times\oo p\infty\times\oo{-\infty}q$, with $u<0<v$. It is also important to note that the definition applies to all points, not just to those belonging to $X$. In fact, the condition is void if $P$ is either an interior or an exterior point,   since then some open box is entirely contained in or entirely disjoint from $X$. So we only need to verify local definability at boundary points and at infinite points. Therefore, any clopen is locally definable. It is not hard to show that a finite Boolean combination of locally definable sets is again locally definable. Moreover, the interior, closure and exterior of a locally definable subset are again locally definable.

To see that a $1$-quasi-cell $S$ is locally definable, write it as the graph of a locally definable, continuous map  $f\colon V\to M$. We leave the infinite points again  to the reader, so fix $P=(a,b)$. If $a$ is exterior to $V$, then $P$ lies in the exterior of $S$,  and there is nothing to show. If $a$ lies in the interior of $V$, then there is some open interval $I\sub V$ containing $a$, and   $S\cap (I\times M)$ is the graph of $f$ restricted at $I$ by \eqref{i:locfct}, and the same is true if $a$ is a boundary point of $V$, say on the left, by choosing   $I=\oo aq\sub V$ for some $q>a$.  Using this, it is not hard to see that $2$-quasi-cells are also locally definable.

\begin{proposition}\label{P:locdefdisc}
In an \omin\ structure, a discrete set is locally definable \iff\ it is closed and bounded. 
\end{proposition}
\begin{proof}
Let $D\sub M^n$ be discrete. If it is not closed, then there is a $P\in\partial D$ not belonging to $D$. But then the intersection $D\cap B$ with any open box $B$ containing $P$ will have $P$ in its closure, so that $D\cap B$ is not closed, whence cannot be definable by Corollary~\ref{C:closeddiscrete}. Similarly, if $D$ is not bounded, say, in the first coordinate on the right, then its intersection with any open box of the form $\oo p\infty\times B'$ will still be unbounded, whence not definable by Corollary~\ref{C:closeddiscrete}. Suppose therefore $D$ is closed and bounded. To check local definability at a boundary point $P$, as it belongs to $D$ by closedness, there is an open box $B$ such that $D\cap B=\{P\}$. To check at an infinite point, we can find an open box around $P$ which is disjoint from $D$.
\end{proof}

\begin{corollary}\label{C:locdefbd}
In an \omin\ structure $\mathcal M$,  the topological boundary of a  locally definable subset $Y\sub M$ is a discrete, closed, bounded set.
\end{corollary}
\begin{proof}
If the locally definable set $\partial Y$ has non-empty interior, it would contain an open interval $I$ and we may shrink this so that $F:=I\cap Y$ is definable.  Since $\partial F=I\cap \partial Y=I$, we get a contradiction with \eqref{i:secbound}. Hence $\partial Y$ has no interior, and so, for $b\in \partial Y$ and an open interval $I$ containing $b$ such that $I\cap \partial Y$ is definable, the latter set, having no interior, must be discrete by \eqref{i:dmin}, and hence, shrinking $I$ further if necessary, $I\cap \partial Y=\{b\}$. Hence $\partial Y$   is discrete, whence bounded and closed by Proposition~\ref{P:locdefdisc}.  
\end{proof}

 Given an arbitrary ordered structure $\mathcal M$, let $\locdef M$ be the structure generated by the locally definable subsets of $\mathcal M$ (formally, we have a language with an $n$-ary predicate $\tt X$ for any locally definable subset $X\sub M^n$, and we interpret $\tt X(\locdef M)$ as the subset $X$). Since the class of locally closed subsets is closed under projection, fibers, and finite Boolean combinations, the definable subsets of $\locdef M$ are precisely the locally definable subsets of $\mathcal M$. 

\begin{corollary}\label{C:omincore}
If $\mathcal M$ is an \omin\ structure, then   $\locdef M$ is locally o-minimal.  
\end{corollary}
\begin{proof}
Given a one-variable definable subset  of $\locdef M$, whence a locally definable subset $Y\sub M$, and a point $a\in M_\infty$, we may choose an open interval $I$ around $a$ such that $Y\cap I$ is definable. Since $a^-$ belongs to $Y\cap I$ or to its complement, the same is true with respect to $Y$, proving \eqref{eq:event}. 
\end{proof}

Since bounded clopens are locally definable but have no infimum, axiom~\eqref{eq:infsup} usually fails and  $\locdef M$ is in general not \omin.

\section{Planar cell decomposition}
 In o-minimality, \emph{cell decomposition} is the property that we can partition any given definable subset $X$ into a disjoint union of cells. Every point is a $0$-cell but writing $X$ as a union of its points should not qualify as a cell decomposition. Slightly less worse, if $X$ is planar, then each fiber $\fiber Xa{}$ is a disjoint union of intervals and points, so that we can partition $X$ into points and vertical cells. Of course, in the o-minimal context these pathologies are avoided  by demanding  the partition be finite. In the \omin\ case, however, we can no longer enforce finiteness, and so to exclude any unwanted partitions, we must impose some weaker restrictions. Moreover, at present, I do not see how to avoid having to use quasi-cells (but see \S\ref{s:tame} below). 

Let us introduce the following terminology, which we give only for the planar case (but can easily be extended to larger arity, see Remark~\ref{R:celldec} below). First we extend the definition of  dimension to arbitrary planar subsets  (which is not necessarily well-behaved if the subset is not definable) by the same characterization: a non-empty subset $B\sub M^2$ has \emph{dimension} $2$, if it has non-empty interior; \emph{dimension} $1$, if it has empty interior but is non-discrete; and \emph{dimension zero} if it is discrete.  We can also define the \emph{local dimension} $\op{dim}_P(B)$ of $B$ at a point $P\in M^2$ as the minimal dimension of $B\cap U$ where $U$ runs over all open boxes containing $P$. Note that $\op{dim}_P(B)\geq 0$ if $P\in\bar B$. It follows that the dimension of $B$ is the maximum of its local dimensions at all points. It is not hard to see that the dimension of $B$ is the largest $e$ for which it contains an $e$-cell. 
In particular, a $2$-(quasi-)cell has dimension $2$, whereas a $1$-(quasi-)cell has dimension one. More generally, by Corollaries~\ref{C:uniondisc} and \ref{C:non-disc2} and the local nature of dimension, we showed that:

\begin{lemma}\label{L:uniondim}
The dimension of a finite union of $e_i$-quasi-cells is equal to the maximum of the $e_i$.\qed
\end{lemma}

Given a collection $\mathcal B$ of (not necessarily definable) subsets of   $M^2$, we say that a definable subset $X\sub M^2$ has a \emph{$\mathcal B$-decomposition}, if there exists a partition $X=\bigsqcup_{i\in I} B_i$ with all $B_i\in\mathcal B$, with the additional property that if $X^{(e)}$ denotes the union of all $e$-dimensional $B_i$ in this partition, then $X^{(e)}$ is definable and has dimension at most $e$, for $e=0,1,2$ (whence of dimension  $e$ \iff\ it is non-empty). Put simply, in a decomposition  there cannot be too many lower dimensional subsets. By a \emph{cell} decomposition (respectively, a \emph{quasi-cell} decomposition) we mean a $\mathcal B$-decomposition where  $\mathcal B$  is the collection of all (quasi-)cells.  By Lemma~\ref{L:uniondim}, any finite partition into (quasi-)cells is a cell decomposition (there can be no quasi-cell in a finite decomposition since each subset in the partition is then definable).

\begin{remark}\label{R:celldec}
For higher arities, we define the dimension of a subset $B\sub M^n$ to be the largest $d$ such that it contains a $d$-cell (in case $B$ is not definable, this might be different from the largest $d$ such that the projection of $B$ onto some $M^d$ has non-empty interior, but both are equal in the definable case). The definition of $\mathcal B$-decomposition for a definable set $X\sub M^n$ now easily generalizes: it is a partition of $X$ into sets from $\mathcal B$ such that the union of all sets in this partition of a fixed dimension is a definable subset of dimension at most $d$.
\end{remark}

\begin{theorem}\label{T:plane}
In an \omin\   structure  $\mathcal M$, any planar definable subset $X\sub M^2$ has a quasi-cell decomposition.
\end{theorem}
\begin{proof}
There is nothing to show if $X$ has dimension zero. If $X$ is a curve,  then $\vertcomp X$ is  a disjoint union of vertical cells (see the proof of Proposition~\ref{P:dimcl}). So after removing it from $X$, we may assume $X$ has no vertical components. In that case,  $\node X$   is    discrete   by Proposition~\ref{P:discnodes}, and so after removing it, we may assume $X$ has no nodes, and so we are done by Corollary~\ref{C:supercell}. 

So remains the case that $X$ is $2$-dimensional. Let $C:=\partial X$ be its boundary. Since $C$ has dimension at most one by Corollary~\ref{C:bddim2}, and so can be decomposed into disjoint quasi-cells by what we just argued, we may assume,  after removing it, that $X$ is moreover open. The projection $\pi(N)$ of the set $N:=\node C$ of all nodes  is discrete by  Proposition~\ref{P:discnodes} and Theorem~\ref{T:projdisc}, and therefore $X\cap (\pi(N)\times M)$ can be partitioned into vertical cells. So remains to deal with points $(a,b)\in X$ such that $a\notin \pi(N)$. Since $X$ is open, $b$ is an interior point of $\fiber Xa{}$. Let  
  $l$ and $h$ be respectively the maximum of $(\fiber Ca{})_{\leq b}$ and the minimum of $(\fiber Ca{})_{\geq b}$, so that the open interval $\oo {l}{h}$ lies inside $\fiber Xa{}$ and contains $b$ (we leave the case that one of these endpoints is infinite to the reader). By choice,  neither    $(a,l)$ nor $(a,h)$ is a node of $C$, so that by Corollary~\ref{C:supercell}, there are (unique) optimal $1$-quasi-cells  $L,H\sub C$ containing $(a,l)$ and $(a,h)$ respectively, say, given as the graphs of  locally definable, continuous maps $f\colon V\to M$ and $g\colon W\to M$. Consider all open intervals $I\sub V\cap W$ containing $a$ such that the $2$-cell $C(I;\restrict {f}I<\restrict {g}I)$ lies entirely inside $X$, and let $Z\sub M$ be the union of all these intervals. Hence $Z$ is open and convex, and  $C(Z;\restrict {f}Z<\restrict {g}Z)$ is an (optimal) $2$-quasi-cell inside $X$, by Lemma~\ref{L:supercell}.  To show that this construction produces a disjoint union of quasi-cells, we need to show that if $(a',b')$ is any point in $S$, then the above procedure yields exactly the same quasi-cell  containing $(a',b')$. Indeed, by convexity, we can find an open interval $I\sub V$ containing $a$ and $a'$. Since the intersections of $F$ and $G$ with $I\times M$ are $1$-cells,  $C(I;\restrict {f}I<\restrict {g}I)=S\cap (I\times M)$ is a $2$-cell contained in $X$, whence must lie inside $S$ by construction. 
  
To show that this is a decomposition, we induct on the dimension $d$ of  $X$, where the case $d=0$ is trivial. In the above, at various stages, we had to remove some subsets of $X$ of dimension strictly less than $d$, and  partition each separately. Since each of these finitely many exceptional sets was definable,  so is their union and by Lemma~\ref{L:uniondim}, has dimension strictly less than $d$. Hence the complement $X^{(d)}$, consisting of all $d$-quasi-cells in the partition, is also definable.   After removing $X^{(d)}$, we are left with a definable subset of strictly less dimension,  and so we are done by induction.
\end{proof}

The proof gives in fact some stronger results, where for the sake of brevity, we will view any point as a $0$-quasi-cell:

\begin{remark}\label{R:plane}
 Keeping track of the various (quasi-)cells, we actually showed that we may partition $X$ in  quasi-cells $S_i$, such that each $\bar S_i\cap X$ is a disjoint union of $S_i$ and   some of the other  $S_j$. 
\end{remark}

%

\begin{remark}\label{R:border}
For a definable subset $Y\sub M$, define a \emph{left outer point} (respectively, a \emph{right outer point}) as a point $b\in Y$ such that $(\partial Y)_{\leq b}$ (respectively, $(\partial Y)_{\geq b}$) is finite. For a subset $X\sub M^2$, we call $(a,b)$ an \emph{outer point}, if $b$ is an outer point   in $\fiber Xa{}$. We can now choose a partition $X=\bigsqcup_i S_i$ in quasi-cells $S_i$  such that any quasi-cell $S_i$ containing an outer point is a cell. 

Indeed, for a fixed $n$, let $f_n\colon M\to X$ be defined by letting $f_n(x)$ be  the  $n$-th left outer point in $\partial (\fiber Xx{})$ (note that this is a definable function by  \eqref{i:succdisc}). The set where $f_n$ takes finite values is definable, and after removing a discrete set, $f_n$ is continuous on it, whence defines a cell above each open interval in its domain. So, given a (left) outer point $(a,b)$, if $n$ is the cardinality of $\partial (\fiber Xa{})_{\leq b}$, then $(a,b)$ lies either on the graph of $f_n$  or between the graphs of $f_n$ and $f_{n+1}$, and therefore, on a cell. 
 \end{remark}
 
 \begin{corollary}\label{C:compactcell}
In an \omin\ structure $\mathcal M$,  any definable subset of $M^2$ lying inside   a compact subset   admits a finite partition into cells.
\end{corollary}
\begin{proof}
We induct on the dimension $d$ of the definable subset $X\sub M^2$ contained in the compact set $K$. By Corollary~\ref{C:closeddiscrete}, if $d=0$, then $X$ is  closed, whence compact, and discrete, whence finite. So assume $X$ is a curve. Let $D$ be the $a\in M$ for which $\fiber Xa{}$ has non-empty interior, so that $D$ is a discrete subset of the compact $\pi_1(K)$, whence finite. By Proposition~\ref{P:locomin}, since each $\fiber Xa{}$ lies in the compact set $ \pi_2(K)$, it is a finite union of intervals, and hence $\vertcomp X$, being the union of all $\fiber Xa{}$ for $a\in D$, is a finite union of cells. Removing it, we reduced to the case that $X$ has no vertical components. In particular, each fiber is discrete whence finite, and so is $\node X$ by Proposition~\ref{P:discnodes}. Upon removing the latter, we may assume $X$ has no nodes.  Apply the quasi-cell decomposition procedure from Theorem~\ref{T:plane} to $X$, which we will call here the \emph{canonical quasi-cell decomposition}. Since all its fibers are finite, every point is an outer   point in the terminology of Remark~\ref{R:border}, and hence belongs to a cell. In other words, the partition does not contain quasi-cells and so is a cell decomposition. Suppose there are infinitely many $1$-cells in this partition. Since each cell inside $K$ must have exactly two boundary points, and $\partial X$ is finite, there must be at least two of these boundary points common to infinitely many cells. But then, by continuity, any vertical line in between these two points would intersect each of these cells in a different point, contradicting that the fibers are all finite.

So remains the case that $d=2$. 
We follow again the procedure, and its notation,  from Theorem~\ref{T:plane} to produce the canonical quasi-cell decomposition of $X$. All the parts of $X$ we will throw out in this procedure are lower dimensional and hence have finite cell decomposition by induction. So we may reduce to the case that $X$ is open and that $(a,b)\in X$ has the property that no node of $\partial X$ lies in  the closure of $\fiber Xa{}$. The definable maps $f\colon V\to M$ and $g\colon W\to M$ as in the proof of Theorem~\ref{T:plane} are uniquely determined by $(a,b)$ and $X$,  and their graphs are part of a cell decomposition of $\partial X$. In particular, by induction, there are only finitely many choices for $f$ and $g$. Since   the $2$-cell $C(Z;\restrict {f}Z<\restrict {g}Z)$  in the canonical cell decomposition containing $(a,b)$ is uniquely determined  by the property that $Z$ is the largest open interval  in $V\cap W$  such that the $2$-cell lies inside $X$,   there are only finitely many possibilities for it.
\end{proof}

\begin{remark}\label{R:compactcell}
The analogue of Proposition~\ref{P:locomin} also holds: if $X\sub M^2$ is definable and $K\sub M^2$ is compact and open, then $X\cap K$ admits a finite cell decomposition. Indeed, choose for each $a\in K$ and open box $U_a$ containing $a$ and contained in $K$. Since each $X\cap U_a$ is definable and contained in $K$, it admits a finite cell decomposition by the previous result. But by compactness, $K$ is the union of finitely many of the $U_a$, proving the claim.
\end{remark}

Both the Corollary and the Remark extend to higher arities, the detail of which we leave to the reader.

\section{Tameness}\label{s:tame}

The quasi-cell decomposition version given by Theorem~\ref{T:plane} is not very useful in applications. Moreover, the non-definable nature of quasi-cells is a serious obstacle. Perhaps quasi-cells never occur, but in the absence of a proof of this, we make the following definitions, for $\mathcal M$   any ordered structure. Let us call a definable map $c\colon X\to Y$ \emph{pre-cellular}, if every fiber $\inverse cy$ is a cell. Note that the non-empty fibers of $c$ then constitute a partition of $X$ into cells. Injective maps are pre-cellular, but the resulting partition in cells is clearly not a decomposition if $X$ has positive dimension. To guarantee that we get a cell decomposition, we require moreover that the image of $c$ be discrete, bounded, and closed, and we call such a map simply \emph{cellular}. In particular, we may assume, if we wish to do so, that the    cellular map $c\colon X\to D$ is surjective, where $D$ is discrete, bounded, and closed.  
%
%

Assume now that  $\mathcal M$ is \omin\ (or, merely a model of $\Ded$), and let $c\colon X\to D$ be   cellular. The collection $X^{(e)}$ of all fibers $\inverse cy$  of dimension $e$ is a definable subset, for each $e$, since we can express in a first-order way whether a fiber $\inverse cy$ has dimension $e$ (for instance, if $X$ is planar, then  having  interior or  being discrete are elementary properties).  If   $X^{(e)}$ is non-empty, then its dimension is equal to $e$ by Corollary~\ref{C:fibdim} and the fact that  $D$ is discrete, showing that we have indeed a cell decomposition. Note that in particular, the graph of $c$    and $X$ have the same dimension. 
We have the following converse: given a definable subset $X\sub M^k$ and a cellular map $c\colon X\to Y$ such that $X=\bigsqcup_{y\in c(X)}\inverse cy$  is a cell decomposition, then $c$ is   cellular, that is to say, $c(X)$ is discrete. It suffices to show that  each $c(X^{(e)})$ is discrete, for $e\leq k$,
and this follows from
  Corollary~\ref{C:fibdim}. It is not clear whether every definable subset admits a   cell decomposition of this type, and so we make the following definition:

\begin{definition}\label{D:tame}
A definable subset $X\sub M^n$  in an ordered structure $\mathcal M$ is called \emph{\tame} if it is the domain of a   cellular map. If every definable subset in $\mathcal M$ is \tame, then we call $\mathcal M$ \emph\tame.
\end{definition}

\begin{lemma}\label{L:tameDed}
Any \tame\ structure is a model of $\Ded$.
\end{lemma}
\begin{proof}
By Theorem~\ref{T:onevar}, it suffices to show that any definable subset $Y\sub M$ is a disjoint union of open intervals and a single discrete, bounded, and closed subset.
Let $c\colon Y\to D$ be   cellular, with $D$ discrete, bounded, and closed. Hence, each fiber $\inverse ca$ must be a one-variable cell, that is to say, either  a point or an open interval.  Let $E$ be the subset of all $a\in D$ for which $\inverse ca$ is a point. Hence $Y\setminus\inverse cE$ is a disjoint union of open intervals, so that upon removing them, we may assume $E=D$, so that  $c$ is a bijection. Since $D$ is discrete,   $\inv c$ is continuous. Therefore, $Y=\op{Im}(\inv c)$ is closed and bounded by \cite[Prop. 1.10]{MillIVT} (which we may invoke since \eqref{eq:infsup} holds).    Moreover, by the same result, any closed subset is again mapped to a closed subset, showing that $\inv c$, whence also $c$,  is a homeomorphism. In particular, $Y$ is discrete, as we needed to show. 
\end{proof}

I do not know whether a \tame\ structure is always \omin; in fact, I do not know whether \eqref{eq:DPP}---an \omin\ property to be defined below in \S\ref{s:gr}---holds in a \tame\ structure. We now investigate when an \omin\ structure is \tame. However, since so far we only used the axioms of $\Ded$ to derive properties of \omin\ structures, everything works within this weaker theory. So, for the remainder of this section, unless noted elsewhere, $\mathcal M$ is a model of $\Ded$. 
Clearly any cell is \tame. Since a (principal) fiber of a cell is again a cell, the same holds for \tame\ subsets. Since a principal projection  of a cell   is   a cell, the collection of  \tame\ subsets  is closed under principal projections (we will generalize this in Corollary~\ref{C:symm} below). Any  finite cell decomposition is easily seen to be given by a cellular map, and hence in particular, any o-minimal structure is \tame. 

%
%

\begin{proposition}\label{P:strongtame}
Suppose $\mathcal M$ and $\tilde{\mathcal M}$ are elementary equivalent ordered structures. If $\mathcal M$ is \tame, then so is $\tilde{\mathcal M}$. 
\end{proposition} 
\begin{proof}
Since both structures have isomorphic ultraproducts, we only need to show that \tame{ness} is preserved under elementary substructures and extensions. The former is easy, so assume $\mathcal M$ is a \tame\ elementary substructure of $\tilde{\mathcal M}$ and let 
$\tilde X$ be a definable subset in $\tilde M^n$. Since \tame{ness} is preserved under fibers, we may assume that $\tilde X$ is definable without parameters, say $\tilde X=\varphi(\tilde{\mathcal M})$. By assumption, there exists a   cellular map $c\colon \varphi(\mathcal M)\to D$, that is to say,   formulae $\gamma$ and $\delta$, with $\gamma(\mathcal M)$ the graph of a map all of whose fibers are cells of dimension at most $n$ and whose image is the discrete, closed, bounded subset $\delta(\mathcal M)$. Since all  this is   first-order, it must also hold in $\tilde{\mathcal M}$, so that $\gamma(\tilde{\mathcal M})$ is the graph of  a   cellular map $\tilde c\colon \tilde X\to \delta(\tilde{\mathcal M})$.
\end{proof} 

By Theorem~\ref{T:Hardy}, the associated Hardy structure of a \tame\ structure is therefore again \tame. An interesting question is whether any  ultraproduct of \tame\ structures is \tame. If so, then   any \omin\ structure is \tame, since it is  elementary equivalent by Corollary~\ref{C:ulomin} with an ultra-o-minimal structure, and the latter would then be \tame, whence so would  the former be by  
Proposition~\ref{P:strongtame}. A priori, the analogue of Lemma~\ref{L:ominred}  will fail, that is to say, \tame{ness} may not be preserved under reducts, as some cellular maps may no longer be definable, but we have:

\begin{proposition}\label{P:tamered}
In a reduct of a \tame\ structure, every definable set admits a cell decomposition.
\end{proposition}
\begin{proof}
Let $\mathcal M$ be a \tame\ structure, $\bar{\mathcal M}$ some reduct, and $X$ an $\bar{\mathcal M}$-definable subset. By assumption,  there exists an ${\mathcal M}$-definable cellular map $c\colon X\to D$, the fibers of which yield an ${\mathcal M}$-cell decomposition of $X$. We will need to show how we can turn this into an $\bar{\mathcal M}$-cell decomposition. As always, we   only treat   the planar case, $X\sub M^2$. There is nothing to show if $X$ is discrete, so assume it is a curve. We already argued that its vertical component $\vertcomp X$ admits a cell decomposition, and so we may remove it. The remaining set of nodes is discrete, and hence may be removed as well, so that we are left with the case that $X$ has no nodes. By Corollary~\ref{C:supercell}, every point of $X$ lies on a unique optimal quasi-cell. Hence if   $C:=\inverse cd$ is one of the cells in the above decomposition, then it is  contained in a    unique $\bar{\mathcal M}$-quasi-cell $S$. As $C$ is then the restriction of $S$ to $I$,  it is  $\bar{\mathcal M}$-definable by Lemma~\ref{L:supercell}.

If $X$ has dimension two, we may assume it is open after removing its boundary, as we already dealt with curves. By Theorem~\ref{T:plane}, there exists an $\bar{\mathcal M}$-quasi-cell decomposition of $X$. Following that proof, we may assume, after removing all points lying on a vertical line containing a node of $\partial X$, that any quasi-cell $S$ in this decomposition is open, and its  boundary  consists of quasi-cells of $\partial X$. By the one-dimensional case, the latter decompose into $\bar{\mathcal M}$-cells, whence so does $S$. 
\end{proof}

\begin{remark}\label{R:tamered}
We have  the following puzzling fact that at least one among the following three statements holds:
\begin{enumerate}
\item\label{i:tamenonomin} there is a \tame\ structure which is not \omin;
\item\label{i:nontameomin} there is an \omin\ structure which is not \tame;
\item\label{i:tamered} any reduct of a \tame\ structure is again \tame.
\end{enumerate}
Indeed, suppose both \eqref{i:tamenonomin} and \eqref{i:tamered} fail. So, by the latter, there is a \tame\ structure $\mathcal M$ with a non-\tame\ reduct $\bar{\mathcal M}$, and by the former, $\mathcal M$ is \omin, whence so is $\bar{\mathcal M}$ by Lemma~\ref{L:ominred}. Hence $\bar{\mathcal M}$ is \omin\ but not \tame. Note that \eqref{i:tamenonomin} implies that $\Ded$ is not equal to $\ominth{}$, and by the discussion preceding the proposition, \eqref{i:nontameomin} implies that  \tame{ness} is not preserved under ultraproducts.
\end{remark}


\begin{lemma}\label{L:onetame}
In a model $\mathcal M$ of $\Ded$, every one-variable definable subset  is \tame. 
\end{lemma}
\begin{proof}
Most proofs involving \tame{ness} will require some coding of disjoint unions, and as  we will gloss over this  issue below,  let me do the proof in detail here. For ease of discussion, let us assume   $Y\sub M$ is bounded (the unbounded case is only slightly more complicated and left to the reader). Assume $0$ and $1$ are distinct elements in $M$. Define $c\colon Y\to M^2$ by letting $c(y)$ be equal to $(y,0)$, in case $y\in \partial Y$; and equal to $(x,1)$ where $x$ is the maximum of $(\partial Y)_{<y}$, in the remaining case. The fiber $\inverse c{d,e}$ is either a point in $\partial Y$ (when  $e=0$), or the interval $\oo d{\sigma_{\partial Y}(d)}\sub Y$. Since its image  is contained in $\partial Y\times\{0,1\}$, the map $c$ is cellular by \eqref{i:disc}.
\end{proof}

\begin{remark}\label{R:onetame}
As it will be of use later, note that by the above argument, we can   refine the cell decomposition given by $c$  as follows:  for any discrete subset $D$ containing $\partial Y$, we can construct a cellular map $c_D\colon Y\to M^2$ whose cells have endpoints in $D$.
\end{remark}

To facilitate the coding of certain disjoint unions, we make the following definitions. Given a definable subset $X\sub M^n$, we define its \emph{valence} as the smallest $d$ such that $X\preceq M^d$, that is to say, the smallest $d$ for which there exists an injective, definable map $X\to M^d$ (see the paragraph preceding Theorem~\ref{T:Euldim} below). We say that \emph{definable discrete subsets are univalent}, if every discrete subset $D\sub M^n$ has valence one. The advantage of the assumption that definable discrete subsets are univalent is  that in the definition of definable cell decomposition, we may always take $n=1$.  If $\mathcal M$ expands an ordered field, then its definable discrete subsets are univalent. Indeed, by induction, it suffices to show that if $D\sub M^{n+1}$ is discrete, then there is a definable, injective map $g\colon D\to M^n$.  The set of lines connecting two points of $D$ is again a discrete set (in the corresponding projective space) and hence we can find a hyperplane   which is non-orthogonal to any of these lines. But then the restriction to $D$ of the  projection   onto this hyperplane is injective.

\begin{proposition}\label{P:defcelldiscunion}
Let $\mathcal M$ be a model of $\Ded$ whose definable discrete subsets are univalent (e.g., an expansion of an ordered field) and let $g\colon X\to M^n$ be a definable map with discrete image. If every fiber $\inverse ga$ is \tame, then so is $X$. 
\end{proposition}
\begin{proof}
Let $A:=g(X)$.  By assumption, there exists for each $a\in A$, a 
  cellular map   $c_a\colon \inverse ga\to D_a$ with $D_a\sub M$ discrete. Let  $D$ be the union of    all $\{a\}\times D_a\sub M^{n+1}$, for $a\in A$.   It follows from Corollary~\ref{C:discfib} that $D$ is discrete. Define $c\colon X\to D$ by the rule $c(x)=(g(x),c_{g(x)}(x))$. To see that this is cellular, note that the fiber over a point $(a,d)\in D$ is equal to $\inverse {c_a}d$ for $d\in D_a$, whence is by assumption a cell.
\end{proof}

\begin{remark}\label{R:defcelldiscunion}
The condition on the univalence of discrete sets can be relaxed: let us say that a   subset  $X\sub M^k$ is \emph{$e$-\tame}, if there is a   cellular map $g\colon X\to M^e$. Hence under the assumption that all definable discrete subsets are univalent, \tame\ is the same as $1$-\tame. If now every fiber $\inverse ga$ is merely $e$-\tame, the above proof still goes through to show  that $X$ is \tame.

It is also not necessary to assume that $\mathcal M$ is a model of $\Ded$, provided we impose that the image of $g$ be also closed and bounded. Indeed, one can directly prove that $D$ as in the proof is discrete, closed and bounded, without any appeal to Corollary~\ref{C:discfib}.
\end{remark}

\begin{theorem}\label{T:defcellBool}
In a model  of $\Ded$, the collection of \tame\ subsets   is closed under (finite) Boolean combinations. 
\end{theorem}
\begin{proof}
Let $\mathcal M\models\Ded$. We only treat the case that the definable discrete subsets are univalent, where the general case follows again by a more careful analysis, based upon  Remark~\ref{R:defcelldiscunion}. Moreover,  we restrict for simplicity to the case of planar subsets, and leave the general case to the reader (by an induction   on the arity). 
 Since the  complement of a cell $V\sub M^2$ is   a finite union of cells, it is \tame. For instance, if $V=C(\oo ab;f<g)$, then its complement consists of the four $2$-cells $\oo{-\infty}{a}\times M$, $\oo {b}\infty\times M$, $C(I;-\infty<f)$ and $C(I;g<\infty)$, and the four $1$-cells, $a\times M$, $b\times M$ and the graphs of $f$ and $g$. Since any union can be written as a disjoint union by taking complements, an application of Proposition~\ref{P:defcelldiscunion} then reduces to showing that the   intersection of two cells $V_1$ and $V_2$ in $M^2$  is \tame.  
This is trivial if either one is discrete, whence a singleton. Suppose $V_1$ is a $1$-cell, given by the definable, continuous map $f_1\colon I_1\to M$. Let $Y$   be the subset  of all $x\in I_1$ such that $(x,f_1(x))$   belongs to $V_2$. Choose a  cellular map $c\colon Y\to D$ (by Lemma~\ref{L:onetame}, or, for higher arities, by induction). Its composition with the (bijective) projection $V_1\cap V_2\to Y$ is then also  cellular.

Suppose next that $V_i=C(I_i;f_i<g_i)$ are both $2$-cells, assumed once more for simplicity  to be bounded.  Let $I:=I_1\cap I_2$ and for $x\in I$,  let $f(x)$ be the maximum of $f_1(x)$ and $f_2(x)$, and let $g(x)$ be the minimum of $g_1(x)$ and $g_2(x)$. Note that   $f$ and $g$ are   continuous on $I$. Let $Y$ consist of all $x\in I$ for which $f(x)<g(x)$, and  let $c\colon Y\to D$ be  cellular. The composition of $c$ with the projection $V_1\cap V_2\to Y$ is again cellular, since its fibers are the cells $C(\inverse ca;f<g)$. 
%
%
\end{proof}

\begin{example}\label{E:nontamenondc}
It is important in this result that the structure is already a model of $\Ded$. For instance, let $D$ be the subset of the ultrapower $\ul{\real}$ of the reals (viewed as an ordered field) consisting of all elements of the form $n$ or $\diag-n$, for $n\in\nat$. Note that $D$ is closed, bounded, and discrete, and hence \tame. However, $(\ul{\real},D)$ is not \tame, since $\nat=D_{<\diag/2}$ is definable in it but fails to satisfy \eqref{eq:infsup}, and so, $(\ul{\real},D)$ is not even a model of $\Ded$. In particular, in the terminology from below, $D$ is not \ofin\ by Corollary~\ref{C:ofindisc}.
\end{example}

It is not hard to show that the product of two cells is again a cell. Therefore, the product of two \tame\ subsets is again \tame. Similarly, the fiber of a cell is again a cell, and hence if $X\sub M^n$ is \tame, then so is each fiber $\fiber X{\tuple a}{}$. Together with Theorem~\ref{T:defcellBool} and the fact that a principal projection of a \tame\ subset is again \tame, we showed that the collection of \tame\ subsets determines a first-order structure on $M$ (in the sense of \cite[Chapt.~1, 2.1]{vdDomin},  with a predicate for every \tame\ subset of $M^n$). Calling this induced structure on $M$ the \emph{\tame\ reduct} of $\mathcal M$ and denoting it $\tamered M$, is justified by: 

\begin{corollary}\label{C:tamered}
If $\mathcal M\models\Ded$, then   $\tamered M$ is \tame, whence in particular a model of $\Ded$.
\end{corollary} 
\begin{proof}
The definable subsets of $\tamered M$ are precisely the \tame\ definable subsets of $\mathcal M$, so that in particular, $\mathcal M$ and $\tamered M$ have the same cells. So remains to show that if $c\colon X\to D$ is cellular in $\mathcal M$, then it is also cellular in $\tamered M$. Discrete sets are \tame\ by definition, and the graph $\Gamma(c)$ of $c$ is \tame\ by applying Proposition~\ref{P:defcelldiscunion} to the projection $\Gamma(c)\to D$. In particular, $c$ is $\tamered M$-definable, and since its fibers are cells, we are done. The last assertion then follows from Lemma~\ref{L:tameDed}.
\end{proof}

\begin{proposition}\label{P:tameredomin}
If $\mathcal M$ is \omin, then so is   $\tamered M$.
\end{proposition}
\begin{proof}
Let $\bar L$ be the language with predicates for the \tame\ subsets of $\mathcal M$, so that $\tamered M$ is an $\bar L$-structure.   Viewing $\mathcal M$ as a structure in the language having a predicate for every definable subset of $\mathcal M$  yields again a  \tame\ structure, since we added no new definable subsets (see Lemma~\ref{L:expdef} below). Therefore,   upon replacing $L$ by the latter language, we assume from the start that $\bar L\sub L$, and the result follows from Lemma~\ref{L:ominred}. 
\end{proof}

%
%
 We also have the following joint cell decomposition:
 
 \begin{corollary}\label{C:jointtame}
Given \tame\ subsets $Y_1,\dots,Y_n$ of a \tame\ subset $X$ in $\mathcal M\models\Ded$, there exists a cellular map  $c\colon X\to D$, such that for each $i$, the restriction of $c$ to $Y_i$ is also cellular.
\end{corollary}
\begin{proof}
Since any Boolean combination of \tame\ subsets is again \tame\ by Theorem~\ref{T:defcellBool}, we may reduce first to the case that all $Y_i$ are disjoint, and then by induction, that we have a single \tame\ subset $Y\sub X$. Since $X\setminus Y$ is \tame\ too, we have cellular maps $d\colon Y\to D$ and $d'\colon X\setminus Y\to D'$, and their disjoint union is then the desired cellular map.
\end{proof}
 
%
%
%

We call a definable map \emph{\tame}, if its graph is. Note that its domain then must also be \tame.  
As already observed in the previous proof,   cellular maps are \tame.  To characterize \tame\ maps, we make the following observation/definition: given a   cellular map $c\colon X\sub M^n\to D$, for $e\leq n$, let 
$X^{(e)}_c=X^{(e)}$ be the union of all $e$-dimensional $\inverse ca$. Since dimension is definable, so is each $X^{(e)}$, and hence the restriction of $c$ to $X^{(e)}$ is also cellular, proving in particular that each $X^{(e)}$ is \tame.

\begin{theorem}\label{T:tamemap}
In a model of $\Ded$,  a definable map $f\colon X\to M^k$  is \tame\ \iff\  $X$ is \tame, and the restriction of  $f$ to the set of its discontinuities   is also \tame. In particular, a definable, continuous map with \tame\ domain is \tame. 
\end{theorem} 
\begin{proof}
If $f$ is \tame, then $f$ is $\tamered M$-definable, and hence so is its set of discontinuities $X'$, proving that $X'$ is \tame. Since the graph of $\restrict f{X'}$ is $\Gamma(f)\cap (X'\times M^k)$, the restricted map is again \tame\  by Theorem~\ref{T:defcellBool}. For  the converse,   $U:=X\setminus X'$ is \tame\ by Theorem~\ref{T:defcellBool}, so that we have a cellular map $c\colon U\to D$. Since $f$ is continuous on $U$, the composition of $c$ with the principal projection $\Gamma(\restrict fU)\to U$ is also cellular, showing that $\Gamma(\restrict fU)$ is \tame. Since by assumption,   the graph of the restriction to $X'$ is   \tame, so is   $\Gamma(f)=\Gamma(\restrict fU)\cup \Gamma(\restrict f{X'})$   by   Theorem~\ref{T:defcellBool}, showing that $f$ is \tame.
\end{proof}

A \tame\ map is $\tamered M$-definable, and hence so its its image, proving:

\begin{corollary}\label{C:symm}
In a model $\mathcal M$ of $\Ded$, if  the domain of  a definable, continuous map is \tame, then so is its image.   More generally, the image of a \tame\ subset under a \tame\ map is again \tame.\qed
\end{corollary}

Let us call a definable map $f\colon X\sub M^n\to M^k$ \emph{almost continuous}, if its set of discontinuities is discrete. By Theorem~\ref{T:disctu}, any one-variable definable map in a model of $\Ded$ is almost continuous. Given a definable map $f\colon X\to M^k$, let us inductively define $D_i(f)\sub X$, by setting $D_0(f):=X$, and by setting $D_i(f)$, for $i>0$,  equal to the set of discontinuities of the restriction of $f$ to $D_{i-1}(f)$. By Remark~\ref{R:genctu}, each $D_i(f)$ has strictly lesser dimension than $D_{i-1}(f)$, and hence $D_n(f)$ is empty for $n$ bigger than the dimension of $X$.  Hence $f$ is (almost) continuous if $D_1(f)$ is empty (respectively, discrete). Since the domain of a \tame\ function is \tame,    an easy inductive argument using Theorem~\ref{T:tamemap} immediately yields:

\begin{corollary}\label{C:almostctu}
In a model of $\Ded$,   an almost continuous (e.g., a one-variable) map with \tame\ domain is   \tame.
 In particular, a definable map $f$ is  \tame\ \iff\ all  $D_i(f)$ are  \tame.
\qed
\end{corollary}


 Let us say that an ordered structure is \emph{almost continuous}, if apart from a binary predicate denoting the order, all other symbols represent almost continuous functions.

\begin{corollary}\label{C:tame}
If  is $\mathcal M$ an almost continuous model of $\Ded$, then $\mathcal M$ is \tame. 
\end{corollary} 
\begin{proof}
Since the collection of \tame\ subsets is closed under Boolean operations, projections, and products by Theorem~\ref{T:defcellBool}, we only have to verify that the ones defined by unnested atomic formulae are \tame. Since by assumption the only predicate is the inequality sign, and the set it defines is a cell, we only have to look at formulae of the form $f(x)=g(x)$ or $f(x)<g(x)$, with $f, g$ function symbols. Since $f$ and $g$ are total functions representing almost continuous maps, their graphs are \tame\ by Corollary~\ref{C:almostctu}, whence so is their intersection by Theorem~\ref{T:defcellBool}. The projection of the latter is the set defined by $f(x)=g(x)$, proving that is a \tame\ subset. Let $F$ and $G$ be the subsets of $M^{n+2}$ of all $(a,f(a),c)$ and all $(a,b,g(a))$ respectively, with $a\in M^n$ and $b,c\in M$. Since these are just products of the respective graphs and $M$, both are \tame, and so is the subset $E$ of all $(a,b,c)\in M^{n+2}$ with $a<b$. Therefore, by another application of Theorem~\ref{T:defcellBool}, the intersection $F\cap G\cap E$ is \tame, and so is its projection, which is just the set defined by the relation $f(a)<g(a)$. 
%
\end{proof}

\begin{remark}\label{R:tame}
More generally, by the same argument, if $\mathcal M\models\Ded$ is an expansion of a \tame\ structure by \tame\ functions and by predicates defining \tame\ subsets, then $\mathcal M$ itself is \tame.  
\end{remark}

\section{\Omin\ expansions of groups}\label{s:omingroup}
In this section, we study \omin\ $L$-expansions of   (ordered) groups. As in the o-minimal case, they are definably simple, whence in particular divisible and Abelian:

\begin{theorem}[Definable simpleness]\label{T:Dedlimgr}
Any \omin\ expansion of a group  is divisible and commutative, and has no proper,    definable  subgroups.
\end{theorem}
\begin{proof}
Let $G$ be an \omin\ expansion of a group. Since we do not assume $G$ to be commutative, we will use the multiplicative notation (although we will continue to write $x^-$ and $x^+$). To show divisibility, fix some $p\in\nat$ and let $f\colon G\to G$ be the $p$-th power map $x\mapsto x^p$. Since $f$ is  definable and continuous, its image  is definably connected, whence an interval by \eqref{i:defconn}. Since the $p$-th powers are cofinal, this image interval must be the whole set, proving divisibility. 
%
%
%
%
In particular, since ordered groups are torsion-free, any element $c\in G$ has by assumption a unique $p$-th root, denoted $\sqrt[p] c$.  

We start with showing that every definable subgroup $1\neq H\subsetneq G$ is  convex. If $H$ were discrete, it would admit a maximum $b>1$ by \eqref{i:disc}. Since then $b<b^2\in H$, we get a contradiction. Let us next show that  $H$ cannot have any isolated point. Indeed, by \eqref{i:dmin}, there are $a<b$ with $\cc ab\sub H$. If $h\in H$,  then $h\inv aH=H$ contains the closed interval $\cc h{h\inv ab}$, showing that $h$ is not isolated. Let $b$ be the infimum of $(M\setminus H)_{>1}$. To establish convexity, it suffices to show that $H_{>b}$ is empty, so towards a contradiction, suppose not and let $c$ be its infimum. Since $b^+$ belongs to $M\setminus H$, we must have $c>b$ and $\oo bc$ is disjoint from $H$. Since neither $1$ nor $c$ can be isolated by our previous observation,  $1<b$ and $c^+$ belongs to $H$. In particular, $\co 1b$ and $\oo ce$, for some $e>c$, are contained in $H$. Choose $d>1$ strictly smaller than $b$ and $e\inv c$. Since $d$ is in $\oo 1b$, it belongs to  $H$, and since $c<dc<e$, so does $dc$. Hence $c\in H$. Since $\sqrt b$ is less than $b$, it too belongs to $H$, whence so does $b$. Hence $c\cc 1b=\cc c{cb}\sub H$. Since $b<b^2\in H$, we  must have $c\leq b^2$. Now, take any $u\in \oo bc$. Hence $u\notin H$, but, on the other hand,  from $b^2<bu<bc$, we get $c<bu<bc$, showing that $bu$ whence also $u$ belongs to $H$, contradiction.

So $H$ is convex, whence an interval with endpoints $b'$ and $b$, the respective infimum and supremum of $H$. Since $H\neq G$, we must have $b\in G$. Since $\sqrt b<b<b^2$, the former   lies in $H$, whence so do $b$ and $b^2$, contradiction.  The commutativity of $G$ follows now immediately by observing that  any   centralizer, being definable and non-trivial,  must be the whole group. 
\end{proof}

\begin{remark}\label{R:ordgr}
Since an \omin\ structure  has the IVT, the fact that an \omin\ expansion of a group  is Abelian and divisible then already follows from  \cite[Proposition 2.2]{MillIVT}.
\end{remark}

From now on, we will write \omin\ expansions of groups   additively.

\begin{lemma}\label{L:Dedlimdisc}
Let $\mathcal M$ be an \omin\  expansion of a group. For any definable subset $Y\sub M$, the topological boundary $\partial Y$ is uniformly discrete, that is to say, there exists $a>0$ such that $a\leq y-x$ for any pair of topological boundary points $x<y$.
\end{lemma}
\begin{proof}
By \eqref{i:secbound}, the topological boundary $\partial Y$ is discrete. Consider the 
the definable map $\partial Y\to M\colon b\mapsto \sigma(b)-b$, where $\sigma=\sigma_{\partial Y}$ is the successor function of $\partial Y$ given by  \eqref{i:succdisc}. By Lemma~\ref{L:imdisc}, its image is also discrete, whence has a minimum $a>0$ by \eqref{i:disc}, satisfying therefore the required property. 
\end{proof}

\begin{remark}\label{R:spread}
It follows from the above proof, that there is a maximal $a$ so that $a\leq y-x$ for all $x<y\in\partial Y$, which we call the \emph{spread} $a_{\text{spr}}(Y)$ (and, in particular, there is some $b\in \partial Y$ so that $a+b\in\partial Y$). 
\end{remark}

\begin{theorem}\label{T:ominth}
Every Archimedean, \omin\ expansion of a group is o-minimal.
\end{theorem} 
\begin{proof}
Let $\mathcal M$ be an \omin\ expansion of a group and let $Y\sub M$ be definable. By \eqref{i:intint}, it suffices to show that $\partial Y$ is finite, and so, upon replacing $Y$ by its boundary, we may already assume that $Y$ is discrete. Towards a contradiction, suppose $Y$ is not finite and let $a:=a_{\text{spr}}(Y)>0$ be its spread. Let $H$ be the (non-empty) set of all $x\in M$ such that $Y_{<x}$ is infinite. In particular, $H$ is (the upper part of) a Dedekind cut. Since $a>0$ and $M$ is Archimedean, $H$ is strictly contained in $H-\frac a2$. Hence,  their difference must contain infinitely many, whence at least two elements    $x<y$ belonging to $Y$. However, this leads to the contradiction  $y-x<a$.
\end{proof}

Recall that a first-order structure $\mathcal M$ is said to have \emph{definable Skolem functions}, if for every definable map $f\colon X\to Y$, there is a \emph{definable section} $g\colon f(X)\to X$, where the latter means that $g$ is definable and $f\after g$ is the identity. Assume now that $\mathcal M$ is an expansion of   an ordered structure. Let $X\sub M^3$ be the subset of all $(a,b,x)$ such that either $a<x<b$, or   $x<a=b$, or $b<a<x$. Note that every open interval occurs exactly once as a fiber of $X$ under the principal projection $\pi\colon M^3\to M^2$. A special case of definable Skolem functions is for this projection to have a definable section $\bary\cdot$, in which case we say that $\mathcal M$ has \emph{definable barycenters} (with the understanding that we assign $\bary M:=c$ to be some point $c\in M$ fixed once and for all). Instead of writing $\bary{a,b}$, we may also write $\bary I$ where $I$ is the open interval $\oo ab$.
Any expansion of a divisible, ordered Abelian group has definable barycenters, namely let $\bary I:=(a+b)/2$ be the midpoint of the bounded open interval $I=\oo ab$, and let it be equal to respectively $a-c$ and $a+c$ if $I$ is equal to $\oo{-\infty} a$ and $\oo a\infty$ respectively, where $c$ is a fixed positive element.

\begin{lemma}\label{L:defskol}
If an \omin\ structure has definable barycenters (e.g., an expansion of an ordered group), then it has definable Skolem functions. 
\end{lemma}
\begin{proof}
Since a section of a definable map $f\colon X\sub M^n\to Y\sub M^m$ is also a section of the projection of its graph onto $M^m$, we may reduce to the case of a projection of a definable subset. By induction on $m$, we then easily reduce to the case that $f=\restrict \pi X$ where $\pi\colon M^{n+1}\to M^n$ and $X\sub M^{n+1}$. For   $\tuple a\in \pi(X)$, let $l\in M_\infty$ be the infimum of the fiber $\fiber X{\tuple a}{}$. If $l\in \fiber X{\tuple a}{}$, 
we put $g({\tuple a}):=(\tuple a,l)$. In the remaining case, $l^+$   belongs to $\fiber X{\tuple a}{}$, and hence  the infimum $b$ of $(M\setminus \fiber X{\tuple a}{})_{>l}$ is a boundary point strictly bigger than $l$. In particular,  $I:=\oo lb$ lies inside $\fiber X{\tuple a}{}$, whence so does $\bary I$. Putting $g(\tuple a):=(\tuple a,l)$, it is now easy to verify that $g$  is a definable section of $\restrict \pi X$. 
\end{proof}

\subsection{Paths}
The following definition can be made in any ordered structure $\mathcal M$. By a \emph{path}  in $M^2$, we mean the image $\Gamma=g(I)$ of a closed interval $I=\cc ab$ under a definable, continuous, injective map $g\colon I\to M^2$. We call $g(a)$ and $g(b)$ the \emph{endpoints} (provided $a$ and/or $b$ are finite); we refer to $\Gamma$ minus its endpoints as an \emph{open path}. A note of caution:   an open path $\Gamma$ therefore  is the image of an open interval $\oo ab$ under a definable, continuous, injective map $g$ with the additional property that its \emph{endpoints} $g(a^+)$ and $g(b^-)$ are different and do not lie on $\Gamma$.  In an \omin\ structure, paths, open or closed,  are one-dimensional, and by continuity,   definably connected.  By Corollary~\ref{C:symm}, every path, being the projection of a graph, is \tame.

To obtain some further properties, we need to assume some additional structure: let us  fix an \omin\ expansion  $\mathcal M$  of a  group. Most   results are proven in the same way as in the o-minimal case (see \cite[Chapt. 6, \S1]{vdDomin}). We will write $\norm x$ for the maximum of all $\norm {x_i}$ if $x=\rij xn$. Clearly, a path is a planar curve, and in fact it is smooth in the following sense (see Remark~\ref{R:regular} for the definition of regular point):

\begin{proposition}\label{P:pathsmooth}
In an \omin\ expansion of a  group,  any point on an open path is $1$-regular.
\end{proposition}
\begin{proof}
Let $\Gamma=g(I)$ be a path, where $g\colon I\to M^2$. Recall that, after applying a translation, the origin $O\in \Gamma$ is $1$-regular if $\rho(O)$ is not a node on $\rho(\Gamma)$ for some translation given by an invertible $2\times2$-matrix over $\mathbb Q$. 
 Write $g(t)=(g_1(t),g_2(t))$. By the Monotonicity Theorem~\ref{T:disctu}, we may subdivide $I$ in open intervals and a discrete subset $D$ so that $g_1$ and $g_2$ are monotone on each interval. The image of these intervals is then easily seen to consist of non-nodes. So remains the points of the form $g(d)$ with $d\in D$. Checking all cases, one sees that there always exists a rotation $\rho$ such that $\rho \after g$ is strictly increasing  in either component on some open interval $\oo ud$ with $u<d$.  One then easily checks that $\rho(g(d))$ is not a node on $\rho(\Gamma)$.   
\end{proof}

\begin{lemma}[Path selection]\label{L:pathsel}
Given a definable subset $X$ and a point $x$ in its frontier $\fr X$, there exists an open path in $X$  with one of its  endpoints equal to $x$.
\end{lemma}
\begin{proof}
Let $H$ be the set of all $\norm {x-y}$,  for $y\in X$. Since $x\in \bar X$,  the type $0^+$ belongs to $H$, so that $\oo 0u\sub H$ for some $u>0$. Hence for every $0<t<u$, there exists $y_t\in X$ such that $\norm{x-y_t}=t$. By Lemma~\ref{L:defskol}, there exists therefore a definable map $g\colon \oo 0u\to X$ such that $\norm {x-g(t)}=t$, for all $t\in \oo0u$. By Theorem~\ref{T:disctu}, by taking $u$ sufficiently small, we may assume $g$ is continuous, so that its image is an open path. By construction, $g(0^+)=x\notin X$, proving the assertion.
\end{proof}
 
\begin{theorem}\label{T:closbdd}
In an \omin\ expansion of a  group, the image of a closed and bounded definable subset under a definable, continuous map is closed and bounded. 
\end{theorem}
\begin{proof}
We restrict our proof once more to $n=2$, so that we have a definable, continuous map $f\colon X\sub M^2\to M^2$ with $X$ bounded and closed. If $f(X)$ were unbounded, we could find, for each $t\in M$, some $x_t\in X$ with $\norm{f(x_t)}>t$. By Lemma~\ref{L:defskol}, we get a definable map $p\colon M\to X$ such that $\norm{f(p(t))}>t$, for all $t\in M$. Since $X$ is bounded and closed, the Monotonicity Theorem~\ref{T:disctu} applied to each component of $p$ shows that $z:=p(\infty^-)$ belongs to $X$. By continuity, $f(z)=(f\after p)(\infty^-)\in M^2$, which is impossible since $\norm{f(p(t))}>t$, for all $t\in M$.

Since $f$ is continuous, it is \tame\ by Theorem~\ref{T:tamemap}, whence $\tamered M$-definable. Therefore, so is its reverse graph $\Gamma^*(f)$, proving that the latter is a \tame\ subset of  $ M^2$. Hence, there exists a cellular map $c\colon\Gamma^*(f)\to D$ with $D$ discrete.  Since $\Gamma^*(f)$ is closed, it is the union of the closures  $\overline{\inverse ca}$ for $a\in D$. By a similar use of Lemma~\ref{L:pathsel} as in the proof of \cite[Chapt.~6, Lemma 1.7]{vdDomin}, each projection  $\pi(\overline{\inverse ca})$ is closed. Moreover, it is not hard to see that   each $x\in M^2$ is contained in an open box $U$  intersecting  at most three among the   $\pi(\overline{\inverse ca})$, and so their union, that is to say, $f(X)=\pi(\Gamma^*(f))$, is closed by Lemma~\ref{L:closunion} below.  
\end{proof}

\begin{lemma}[Folklore]\label{L:closunion}
In a metric space,  given a locally finite collection of closed subsets, that is to say, such that   each point   admits an open neighborhood  which intersects only finitely many of its members, then their union is again closed.
\end{lemma}
\begin{proof}
Let $X$ be the union of closed subsets $X_i$ with the given property and let $x$ be a point in the closure of $X$. Let $x_n\in X$ be a sequence converging to $x$ and let $U$ be an open containing $x$  intersecting only finitely many $X_i$. In particular, one of these $X_i$ then contains a subsequence of $x_n$. Since $X_i$ is closed, $x\in X_i\sub X$. 
\end{proof}

Without the boundedness assumption,  Theorem~\ref{T:closbdd} even fails for the closure of a cell under projection: consider the graph of $1/x$ in $\real_{>0}$, which is closed but whose projection is open. Theorem~\ref{T:closbdd} has the usual corollaries (see \cite[Chapt. 6, \S1]{vdDomin}): let $f\colon X\to M^n$ be a definable, continuous map  with  closed, bounded domain $X$, then (i) a subset    of $f(X)$ is closed \iff\ its preimage in $X$ is; (ii) a definable map $f(X)\to M^n$ is continuous \iff\ its composition with $f$ is; (iii) if $n=1$, then $f$  attains its maximum and minimum; (iv)  if $f$ is injective, then it is a  homeomorphism   onto its image $f(X)$ (in particular, any such   map is open).

\section{The \gr\ of an \omin\ structure}\label{s:gr}
Given any first-order structure $\mathcal M$ in a language $L$, we define its \emph\gr\ $\grot M$ as follows.  Given two formulae $\varphi(x)$ and $\psi(y)$ in $L(M)$ (that is to say,  allowing parameters), with $x=\rij xn$ and $y=\rij ym$, we say that $\varphi$ and $\psi$ are \emph{$\mathcal M$-definably isomorphic} (or simply \emph{definably isomorphic} if $\mathcal M$ is understood), if there exists a definable bijection   $f\colon \varphi(\mathcal M)\to \psi(\mathcal M)$. Let $\grot M$ be the quotient of the free Abelian group generated by ${\mathcal M}$-definable  isomorphism classes $\isoclass\varphi$ of formulae $\varphi\in L(M)$ modulo the \emph{scissor relations}
\begin{equation}\label{eq:sciss}
\isoclass\varphi+ \isoclass\psi-\isoclass{\varphi\en\psi}-\isoclass{\varphi\of\psi}\tag{sciss}
\end{equation}
where $\varphi,\psi$ range over all pairs of formulae in the same free variables. See for instance \cite{KraEul,KraSca} for more details. 

We will   write $\class \varphi$ or $\class Y$ for the image of the formula $\varphi$, or the set $Y$ defined by it, in $\grot M$. Since we can always replace a definable subset with a definable copy that is disjoint from it, the scissor relations can be simplified, by only requiring them for disjoint unions: $\class {X\sqcup Y}=\class X+\class Y$. In particular, combining all terms with a positive sign as well as all terms with a negative sign by taking disjoint unions, we see that every element in the \gr\ is of the form $\class X-\class Y$, for some definable subsets $X$ and $Y$. To make $\grot M$ into a ring, we define the product of two classes $\class\varphi$ and $\class\psi$ as the class of the \emph{product} $\varphi(x)\en\psi (y)$ where $x$ and $y$ are disjoint sets of variables. One checks that this is well-defined and that the class of a point is the unit for multiplication, therefore denoted $1$. Note that in terms of definable subsets, the product corresponds to the Cartesian product and the scissor relation to the usual inclusion/exclusion relation. 
 Variants are obtained by restricting the class of formulae/definable subsets. For our purposes, we will only do this for discrete subsets. Call a formula \emph{discrete} if it defines a discrete subset. In an \omin\ structure, discrete formulae are closed under Boolean combinations and products, and if two discrete definable subsets are definably isomorphic, then the graph of this isomorphism is also given by a discrete formula. Therefore, the \gr\ on discrete formulae is well-defined and  since definably discrete is equivalent with having dimension zero, we will   denote this  by $\grotth M0$. 
We have a canonical \homo\ $\grotth M0\to \grot M$ with image the subring generated by classes of discrete formulae.
The following is useful when dealing with \gr{s}:

\begin{lemma}\label{L:grequal}
Two definable subsets $X$ and $Y$ in a first-order structure $\mathcal M$ have the same class in $\grot M$ \iff\ there exists a definable subset $Z$ such that $X\sqcup Z$ and $Y\sqcup Z$ are definably isomorphic. 
\end{lemma}
\begin{proof}
One direction is immediately, for if $X\sqcup Z$ and $Y\sqcup Z$ are definably isomorphic, then $\class X+\class Z=\class {X\sqcup Z}=\class{Y\sqcup Z}=\class Y+\class Z$ in $\grot M$, from which it follows $\class X=\class Y$. Conversely, if $\class X=\class Y$, then by definition of scissor relations, there exist mutually disjoint, definable subsets $A_i,B_i,C_i,D_i\sub M^{n_i}$ such that
$$
\isoclass X+\sum_i\isoclass{A_i}+ \isoclass{B_i}-\isoclass{A_i\sqcup B_i}=  \isoclass Y+\sum_i\isoclass{C_i}+ \isoclass{D_i}-\isoclass{C_i\sqcup D_i}
$$
in the free Abelian group on isomorphism classes. Bringing the terms with negative signs to the other side, we get an expression in which each term on the left hand side must also occur on the right hand side, that is to say, the collection of all isomorphism classes $\{\isoclass X, \isoclass {A_i}, \isoclass {B_i}, \isoclass {C_i\sqcup D_i}\}$ is the same as the collection of all isomorphism classes  $\{\isoclass Y, \isoclass {C_i}, \isoclass {D_i}, \isoclass {A_i\sqcup B_i}\}$. By properties of disjoint union, we therefore get $\isoclass{X\sqcup Z}=\isoclass{Y\sqcup Z}$, where   $Z$ is the disjoint union of all  definable subsets $A_i,B_i,C_i,D_i$.  
\end{proof}

If $\mathcal M$ is an expansion of an ordered, divisible  Abelian group, then we have the following classes of open intervals. If $I=\oo ab$, then $I$ is definably isomorphic to $\oo 0{b-a}$ via the translation $x\mapsto x-a$. Moreover, $\oo 0a$ is definably isomorphic to $\oo 0{2a}$ via the map $x\mapsto 2x$. Hence the class $\mathbbm i$ of $\oo 2a$ is by \eqref{eq:sciss} equal to the sum of the classes of $\oo 0a$, $\{a\}$, and $\oo a{2a}$. In other words, $\mathbbm i=2\mathbbm i+1$, whence $\mathbbm i=-1$ (the additive inverse of $1$). Let $\mathbbm h$ be the class of the unbounded interval $\oo 0\infty$. By translation and/or the involution $x\mapsto -x$, any half  unbounded interval is definably isomorphic with $\oo 0\infty$. Finally, we put $\lef:=\class M$ (the so-called \emph{Lefschetz class}). Since $M$ is the disjoint union of $\oo{-\infty}0$, $\{0\}$, and $\oo0\infty$, we get 
\begin{equation}\label{eq:eulint}
\lef=2 \mathbbm h+1.\tag{lef}
\end{equation} 
If $M$ is moreover an ordered field, then taking the reciprocal makes $\oo 01$ and $\oo 1\infty$ definably isomorphic, so that $\mathbbm h=\mathbbm i=-1$, and hence also $\lef=-1$.

 Under the assumption of an underlying ordered structure, whence a topology, we can also strengthen the definition by calling two definable subset \emph{definably homeomorphic}, if there exists a definable (continuous) homeomorphism between them, and then build the \gr, called the \emph{strict \gr} of $\mathcal M$ and denoted $\grots M$,   on the free Abelian group generated by homeomorphism classes of definable subsets. Note that there is a canonical surjective \homo\ $\grots M\to \grot M$. In the o-minimal case, the monotonicity theorem implies that both variants are  equal, but this might fail in the \omin\ case, since  cell decompositions are no longer finite (but see Corollary~\ref{C:tamestrict} below). In fact, in the o-minimal case, 
the \gr\ is extremely simple, as   observed by Denef and Loeser (see \cite[Chap. 4, \S2]{vdDomin} for more details):

\begin{proposition}\label{P:omingr}
The \gr\ of an o-minimal expansion of an ordered field is canonically isomorphic to the ring of integers $\mathbb Z$.
\end{proposition}
\begin{proof}
By the previous discussion, the class of any open interval is equal to $-1$. The graph of a function   is definably isomorphic with its domain, and so the class of any $1$-cell is equal to $-1$. Since a bounded planar $2$-cell lies in between two $1$-cells, it is definably isomorphic to an open box, and by definition of the multiplication in $\grot M$, therefore its class is equal to $\lef^2=1$. The unbounded case is analogous, and so is the case that the $2$-cell lies in a higher Cartesian product. This argument easily extends to show that the class of a $d$-cell in $\grot M$ is equal to $\lef^d=(-1)^d$. By Cell Decomposition, every definable subset  is a finite union of cells, and hence its class in $\grot M$ is an integer (multiple of $1$). 
\end{proof}

We denote the canonical \homo\ $\grot M\to \zet$ by $\euler M\cdot$ and call it the   \emph{Euler \ch} of $\mathcal M$. Inspired by  \cite{CluEd}, we define the \emph{Euler measure} of a definable subset $X$ in an o-minimal structure $\mathcal M$ the pair $\eulerm MX:=(\op{dim}(X),\euler MX)\in (\nat\cup\{-\infty\})\times\zet$, where we view the latter set in its lexicographical ordering.

 In an arbitrary first-order structure, let us  say, for definable  subsets $X$ and $Y$,   that  $ X\preceq Y$ \iff\ there exists a definable injection $X\to Y$. In general, this relation, even up to definable isomorphism,  will fail to be symmetric (take for instance in the reals the sets $X=\cc 01$ and $Y=X\cup\{3/2\}$, where $x\mapsto x/2$ sends $Y$ inside $X$), and therefore is in general only a partial pre-order. As we will discuss below in \S\ref{s:disc}, it does induce a partial order on isomorphism classes of discrete, definable subsets in an \omin\ structure, and we can even make it total by extending the class of isomorphisms (Theorem~\ref{T:totorddisc}). In the o-minimal case, $\preceq$ is a total pre-order  by the following (folklore) result.  In some sense, the rest of the paper is an attempt to extend  this result to the \omin\ case.

\begin{theorem}\label{T:Euldim}
In an o-minimal expansion of an ordered field, two definable sets $X$ and $Y$ are definably isomorphic \iff\ $\eulerm MX=\eulerm MY$.  Moreover,   $X\preceq Y$ \iff\ $\op{dim}(X)\leq \op{dim}(Y)$ with the additional condition that  $\euler  MX\leq \euler  MY$ whenever both are zero-dimensional. 
\end{theorem}
\begin{proof}
The first statement is proven in \cite[Chap. 8, 2.11]{vdDomin}. So, suppose 
$X\preceq Y$. Since $X$ is definably isomorphic with a subset of $Y$, its dimension is at most that of $Y$. If both are zero-dimensional, that is to say, finite, then the pigeonhole principle gives $\euler  MX=\norm X\leq\norm Y =\euler  MY$.

Conversely, assume $\op{dim}(X)\leq \op{dim}(Y)$. If both are finite, the assertion is clear by the same argument, so assume they are both positive dimensional. 
Without loss of generality, by adding a cell of the correct dimension, we may  then assume that they have both the same dimension $d\geq 1$. 
Let $e:= \euler  MY-\euler  MX$ and let $F$ consist of $e$ points disjoint from $X$ if $e$ is positive and of $-e$ open intervals disjoint from $X$ if $e$ is negative. Since $\euler{{}}F=e$, the Euler measure of $X\sqcup F$ and $Y$ are the same, and hence they are definably isomorphic by the first assertion, from which it follows that $X\preceq Y$. 
%
%
\end{proof}

Let ${\mathcal M}$ be an ultra-o-minimal structure, say, realized as the ultraproduct of o-minimal structures  $\mathcal M_i$. We define the \emph{ultra-Euler \ch} $\euler{{ M}}\cdot$ as follows. Let $Y\sub  M^n$ be a definable subset, say given by a formula $\varphi(x,{\tuple b})$ with $ {\tuple b}$ a tuple of parameters realized as the ultraproduct of  tuples $\tuple b_i$ in each $M_i$. Let $Y_i:=\varphi(\mathcal M_i,\tuple b_i)$, so that $ Y$ is the ultraproduct of the $Y_i$, and let $\euler{{ M}}{ Y}$ now be the ultraproduct of the $\euler{M_i}{Y_i}$, viewed as an element of $\ul\zet$. If $ X$ is definably isomorphic with $ Y$, via a definable bijection with graph $ G$, choose as above definable subsets $X_i$ and $G_i$ in $\mathcal M_i$ with ultraproduct equal to $ X$ and $ G$ respectively. By \los, almost each $G_i$ is the graph of a definable bijection between $X_i$ and $Y_i$, and therefore  $\euler{M_i}{X_i}=\euler{M_i}{Y_i}$ for almost all $i$, showing that $\euler{{ M}}{ X}=\euler{{ M}}{ Y}$. Similarly, we define the \emph{ultra-Euler measure} $\eulerm{ M}X:=(\op{dim}( X),\euler MX)$. 
Since the ultra-Euler \ch\ is easily seen to be also compatible with the scissor relations~\eqref{eq:sciss}, we showed (compare with Theorem~\ref{T:Drank} below):

\begin{corollary}\label{C:ulomingrot}
For an ultra-o-minimal structure $ {\mathcal M}$, we have a canonical \homo\ $\grot{{  M}}\to \ul\zet$.\qed
\end{corollary}

\subsection{The Discrete Pigeonhole Principle} 
Before we proceed, we identify another \omin\ property, that is to say, a first-order property of o-minimal structures:

\begin{proposition}[Discrete Pigeonhole Principle]\label{P:DPP}
Given an \omin\ structure $\mathcal M$, if a definable map $f\colon Y\to Y$, for some $Y\sub M^n$, is injective and its image is co-discrete, meaning that $Y\setminus f(Y)$ is discrete, then it is a bijection. In particular, any definable map from a discrete subset $D$ to itself is injective \iff\ it is surjective. 
\end{proposition} 
\begin{proof}
For each formula $\varphi(x,y,\tuple z)$, we can express in a first-order way that if $\varphi(x,y,\tuple c)$, for some tuple $\tuple c$ of parameters, defines the graph of an injective map $f\colon Y\to Y$ then 
\begin{equation}\label{eq:DPP}
\text{$Y\setminus f(Y)$  discrete implies $Y=f(Y)$.}\tag{DPP}
\end{equation}
Remains to show that \eqref{eq:DPP} holds in any o-minimal structure $\mathcal M$. Indeed, if $D=Y\setminus f(Y)$, then $\euler M Y=\euler M{f(Y)}+\euler MD$. Since $f$ is injective, $Y$ and $f(Y)$ are definably isomorphic, whence have the same Euler \ch,  and so $\euler MD=0$. But a discrete subset in an o-minimal structure is finite and its Euler \ch\ is then just its cardinality, showing that $D=\emptyset$. One direction in the last assertion is immediate, and for the converse, assume $f\colon D\to D$ is surjective. For each $x\in D$, define $g(x)$ as the (lexicographical) minimum of $\inverse fx$, so that $g\colon D\to D$ is an injective map, whence surjective by the above, and therefore necessarily the inverse of $f$.  
\end{proof}

At present, I do not know how to derive \eqref{eq:DPP} from $\Ded$.

\begin{corollary}\label{C:ominicgr}
An \omin\  structure is o-minimal \iff\ its \gr\ is isomorphic to $\mathbb Z$. 
\end{corollary}
\begin{proof}
If $\mathcal M$ is \omin, but not o-minimal, then there exists at least one definable, discrete, infinite set $D$. By assumption,  $\class D=n$ for some integer $n$. After removing $n$ points, if $n$ is positive, or adding $-n$ points, if negative, we may suppose $\class D=0$. By Lemma~\ref{L:grequal}, there exists a definable subset $X$ such that $X$ and $X\sqcup D$ are definably isomorphic. By \eqref{eq:DPP}, this forces $D=\emptyset$, contradiction.
\end{proof}

\begin{corollary}\label{C:discmon}
A monotone map $f\colon D\to D$ on  a definable, discrete subset $D$ in an \omin\ structure $\mathcal M$ is either constant or an involution.
\end{corollary}
\begin{proof}
Suppose $f$ is  non-constant and hence $f^2$ is strictly increasing. So  upon replacing $f$ by its square, we may already assume  that $f$ is increasing, and we need to show that it is then the identity.  Since $f$ is injective, it is  bijective by Proposition~\ref{P:DPP}.   Let $h$ be the maximum of $D$, and suppose $f(d)=h$. If $d<h$, then $h=f(d)<f(h)\in D $, contradiction, showing that $f(h)=h$. If $f$ is not the identity, then the set $Q$ of all $d\in D$ for which $f(d)\neq d$ is non-empty, whence has a maximum, say, $u<h$. In particular, if $v:=\sigma_D(u)$ is its immediate successor, then $f(u)<f(v)=v$, since $v\notin Q$, whence $f(u)<u$, since $u\in Q$. Since $u=f(a)$ for some $a\neq u$, then either $a<u$ or $v\leq a$, and hence $u=f(a)<f(u)<u$ or $v=f(v)\leq f(a)=u$, a contradiction either way.
\end{proof}

\begin{remark}\label{R:discmon}
Note that the map sending $h$ to the minimum of $D$, and equal to $\sigma_D$ otherwise is a definable permutation of $D$, but it obviously fails to be monotone.
The map $x\mapsto \diag-x$ on $D=(\ul\nat)_{\leq \diag}$ as in Example~\ref{E:ulomin} is a strictly decreasing involution. It is not hard to see that if an involution exists, it must be unique: indeed, if $f$ and $g$ are both decreasing, let $a$ be the maximal element at which they disagree (it cannot be $h$ since $f(h)=l=g(h)$), and assume $f(a)<g(a)$. Since $f(\sigma(a))=g(\sigma(a))<f(a)<g(a)$, it is now easy to see that $f(a)$ does not lie in the image of $g$, contradicting that $g$ must be a bijection by \eqref{eq:DPP}.
\end{remark}

\begin{proposition}\label{P:onevar}
In an \omin\ expansion $\mathcal M$ of an ordered field, there exists for every definable subset $Y\sub M$, two definable, discrete subsets $D, E\sub Y$ such that $\class Y=\class D-\class E$ in $\grot M$.
\end{proposition}
\begin{proof}
Since  the boundary $\partial Y$ is discrete, we may remove it and assume $Y$ is open, whence a disjoint union of open intervals by Theorem~\ref{T:onevar}. 
Let  us introduce some notation that will be useful later too, assuming $Y$ is open. For $y\in Y$, let $l(y)$ and $h(y)$ be respectively the maximum of $(\partial Y)_{<y}$   and the minimum of  $(\partial Y)_{>y}$ (allowing $\pm\infty$). Hence $\oo {l(y)}{h(y)}$ is the maximal interval in $Y$ containing $y$, and we denote its midpoint by  $m(y):=\bary{\oo {l(y)}{h(y)}}$. Let  $L(Y)$, $M(Y)$ and  $R(Y)$ consist respectively of all $y$ less than, equal to, or greater than $m(y)$.   Removing a maximal unbounded interval from $Y$ if necessary (whose class is equal to $-1$ as already observed above), we may assume $Y$ is bounded, so that $l(y)$ and $h(y)$ are always finite.
Since the maps $f_Y\colon L(Y)\to Y\colon y\mapsto 2y-l(y)$ and $g_Y\colon L(Y)\to R(Y)\colon y\mapsto y+m(y)$ are bijections,   $\class Y=\class{L(Y)}=\class{R(Y)}$.  Since the  scissor relations yield 
   $\class Y=\class {L(Y)}+\class {M(Y)}+\class {R(Y)}$, we get $\class Y=-\class {M(Y)}$. By construction $M(Y)$ is discrete, and so we are done.
\end{proof}

The proof gives the following more general result: given any discrete subset $D_0\sub Y$, we can find disjoint discrete subsets $D,E\sub Y$ such that $D_0\sub D$ and $\class Y=\class D-\class E$. Indeed, let $D:=D_0\cup (Y\cap \partial Y)$ and $E:=M(Y\setminus D)$. 
If $\mathcal M$ merely expands an ordered group, then we have to also include the class $\mathbbm h$, that is to say, in that case we can write $\class Y=e\mathbbm h+\class D-\class E$, where $e\in\{0,1,2\}$ is the number of unbounded sides of $Y$. For higher arities, we need to make a \tame{ness} assumption:
%

\begin{corollary}\label{C:defcell}
Let $X\sub M^2$ be a definable subset in an \omin\ expansion $\mathcal M$ of an ordered field. If $\partial X$ is \tame, then there exist definable, discrete subsets $D,E\sub X$ such that $\class X=\class D-\class E$ in $\grot M$. 
In fact, the class of any \tame\ subset in  $\grot M$ is of the form $\class D-\class E$, for some definable discrete subsets   $D,E\sub M$.
\end{corollary}  
\begin{proof}
There is nothing to show if $X$ is discrete. Assume next that it has dimension one. Let   $V:=\vertcomp X$ be the vertical component of $X$. Since $\pi(V)$ is discrete, as we argued before, we can carry out the argument in the proof of Proposition~\ref{P:onevar} on each fiber separately to write $\class V$ as the difference of two discrete classes (we leave the details to the reader, but compare with  the two-dimensional case below). Removing $V$ from $X$, we may assume $X$ has no vertical components. In particular, the set $N:=\node X$ of nodes of $X$ is discrete by Proposition~\ref{P:discnodes}. Removing it, we may assume $X$ has no nodes, so that every point lies on a unique optimal quasi-cell by Corollary~\ref{C:supercell}. However, by assumption, $X$ is  \tame, and hence there exists a cellular map $c\colon X\to D$. Given $x\in X$, let $I_x$ be the domain $\pi(\inverse c{c(x)})$ of the unique cell $\inverse c{c(x)}$ containing $x$. Let  $L(X)$, $M(X)$, and $R(X)$ consist respectively of all $x\in X$ such that 
$\pi(x)$ lies in $L(I_x)$, $M(I_x)$, and $R(I_x)$ respectively  (in the notation of the proof of  Proposition~\ref{P:onevar}). Define $f_X\colon L(X)\to X$ and $g_X\colon L(X)\to R(X)$ by sending $x$ to the unique point on $\inverse c{c(x)}$ lying above respectively $f_{I_x}(\pi(x))$ and $g_{I_x}(\pi(x))$, showing that $X$, $L(X)$, and $R(X)$ are definably isomorphic. Since  $M(X)$ is discrete and  $\class X=\class {L(X)}+\class {R(X)}+\class {M(X)}$,   we are done in this case.

If $X$ has dimension two, its boundary has dimension at most one, and so we have already dealt with it by the previous case. Upon removing it,  we may assume $X$ is open. This time, we let  $L(X)$, $M(X)$ and $R(X)$ be the union of respectively all $L(\fiber Xa{})$,  $M(\fiber Xa{})$, and  $R(\fiber Xa{})$, for all $a\in\pi(X)$. The maps $(a,b)\mapsto f_{\fiber Xa{}}(b)$ and  $(a,b)\mapsto g_{\fiber Xa{}}(b)$ put $L(X)$ in definable bijection with respectively  $X$ and $R(X)$ (with an obvious adjustment left to the reader if the fiber $\fiber Xa{}$ is unbounded), and hence $\class X=-\class{M(X)}$. Since $M(X)$ has dimension at most one by Proposition~\ref{P:nondisc2}, we are done by induction.
 Without providing the details, we can extend this argument to higher dimensions, proving the last claim, where we also must use that definable discrete subsets are univalent in an ordered field.   
\end{proof}

\begin{remark}\label{R:defcell}
Inspecting the above proof, we actually proved the following: if $c\colon X\to D$ is a cellular surjective map, then
\begin{equation}\label{eq:eulerdisc}
\class X=\sum_{e=0}^d (-1)^e\class {D_e}
\end{equation}
where $D_e=c(X^{(e)})$ consist of all $a\in D$ with $e$-dimensional fiber $\inverse ca$, and  where $d$ is the dimension of $X$. We may reduce to the case that all fibers have the same dimension, and the assertion is then  clear in the one-dimensional case, since the restriction of $c$ to $M(X)$ is a bijection. Repeating the argument therefore to $X$, we get $\class X=-\class{M(X)}=\class {M(M(X))}$, and  now $M(M(X))$ is definably isomorphic with $D$ via $c$. Higher dimensions follow similarly by induction.
\end{remark} 

In particular, if $\mathcal M$ is a  \tame\ expansion of an ordered field, then its \gr\ is generated by the definable discrete subsets of $M$, and the canonical \homo\ $\grotth M0\to \grot M$ is surjective. Inspecting the above proof, we see that all isomorphisms involved are in fact homeomorphisms, and so the result also holds in the strict \gr\ $\grots M$. Since any function with discrete domain is continuous, we showed:

\begin{corollary}\label{C:tamestrict}
For a \tame, \omin\ expansion of an ordered field, its  \gr\ and its strict \gr\ coincide. \qed
\end{corollary}

\subsection{The partial order on $\disc M$}\label{s:disc}
 Let $\disc M$ denote   the collection of isomorphism classes of definable, discrete subsets in an \omin\ structure $\mathcal M$. Recall that  $ X\preceq Y$ if there exists a definable injection $X\to Y$. Theorem~\ref{T:Euldim} suggests that this relation is not very useful in higher dimensions, so we study it only on $\disc M$. Assume   $D$ and $E$ are discrete, definable subsets with $ D\preceq  E$ and $ E\preceq  D$. Hence there are definable injections $D\to E$ and $E\to D$. By Proposition~\ref{P:DPP}, both compositions are bijections, showing that $D$ and $E$ are definably isomorphic. Since transitivity is trivial, we showed that we get a partial order on $\disc M$. To obtain a partial order on the zero-dimensional \gr\ $\grotth M0$, we define $\class D\preceq \class E$, if there exists a definable, discrete subset $A$ such that $D\sqcup A\preceq E\sqcup A$. To show that this well-defined, assume $\class D=\class{D'}$ and $\class E=\class{E'}$. By Lemma~\ref{L:grequal}, there exist definable, discrete subsets $F$ and $G$ such that $D\sqcup F\iso D'\sqcup F$ and $E\sqcup G\iso E'\sqcup G$. Since $D\sqcup A\preceq E\sqcup A$, we get
 $$
 D'\sqcup F\sqcup G\sqcup A\iso D\sqcup F\sqcup G\sqcup A\preceq E\sqcup F\sqcup G\sqcup A\iso E\sqcup F\sqcup G\sqcup A
 $$
as required. We then extend this to a partial ordering on $\grotth M0$ by linearity.  
In the o-minimal case, $\grotth M0$ is just $\zet$ in its natural ordering, but we will give some examples  where the order is not total (though, see Theorem~\ref{T:totorddisc} below). Let us first  prove a comparison result in a special case. In an expansion of an ordered group, we call  a definable, discrete set $D$   \emph{equidistant}, if  the map $a\mapsto \sigma_D(a)-a$ is constant on all non-maximal elements of $D$, where $\sigma_D$ is the successor function. In other words, $D$ is equidistant, if  for any non-maximal $a\in D$, also $a+\rho\in D$, where $\rho$ is the spread of $D$.

\begin{proposition}\label{P:equidist}
In an \omin\ expansion $\mathcal M$ of an ordered field, any two definable equidistant   subset{s} of $M$  are comparable. 
\end{proposition}
\begin{proof}
Let $D,E\sub M$ be definable  equidistant subsets. Since they are bounded by \eqref{i:disc}, we may assume after a translation that both have minimum equal to $0$, and then after taking a scaling, that both have spread $1$. Let $m$ be the maximum of (the non-empty) $D\cap E$. If $m$ is non-maximal in  either set, then $m+1$ lies both in $D$ and in $E$ by assumption, contradiction. Hence $m$ is the maximum, say, of $D$, and therefore $D\sub E$, whence $  D\preceq   E$.
\end{proof}

More generally, given a definable, discrete subset $D\sub M$ in an \omin\ expansion $\mathcal M$ of an ordered field, define the \emph{derivative} $D'$ of $D$ as the set of all differences $\sigma_D(a)-a$, where $a$ runs over all non-maximal elements of $D$. Hence an equidistant set is one whose derivative is a singleton. Since we have a surjective map $D\setminus\{\max D\}\to D'\colon a\mapsto \sigma_D(a)-a$, it follows from the next  lemma that $D'\preceq D$. 

\begin{lemma}\label{L:imageorder}
In an \omin\ structure $\mathcal M$, if   $g\colon X\to M^k$ is a definable map, then $ {g(D)}\preceq  D$, for every discrete, definable subset $D\sub X$. 
\end{lemma}
\begin{proof}
This follows by considering the injective map $ g(D)\to D$  sending $a$ to the minimum of $\inverse ga$.
\end{proof}

In general, this partial order will not be total.  Since $D\preceq E$ implies $\class D\preceq \class E$, but not necessarily the converse, the former being total implies that the latter is too, but again, the converse is not clear. To construct examples, let us introduce the following notation. 

\begin{example}[Discrete Overspill]\label{E:discoverspill}
Given a    sequence $\tuple a=(a_n)$ of real numbers,  let  $\overspill{\tuple a}$  be the ultraproduct of the $\mathcal R_n$, where each $\mathcal R_n$ is the expansion of the real field with a unary predicate $\tt D$ interpreting the first $n$ elements $a_1,\dots,a_n$  in   the sequence. Since each $\mathcal R_n$ is o-minimal, $\overspill{\tuple a}$  is \omin. Moreover, ${\tuple a}$ is the ``finite'' part of the set $D_{\tuple a}:=\tt D(\overspill {\tuple a})$ defined by $\tt D$,   that is to say, 
$$
D_{\tuple a}\cap \real=\{a_1,a_2,\dots,\}.
$$
so that we refer to  $\overspill{\tuple a}$  as the structure obtained from  $\tuple a$ by \emph{discrete overspill} (for a related construction, see also \S\ref{s:Taylor} below).
\end{example}

In this notation, Example~\ref{E:ulomin} is the discrete overspill $\overspill\nat$ of $\nat$ listed in its natural order. I do not know whether $\preceq$ is total on it. Any countable subset can be enumerated, including $\mathbb Q$, although this enumeration might not be order preserving. Nonetheless, we get a structure $\overspill{\tuple q}$ with $D_{\tuple q}\cap \real=\mathbb Q$ (the non-standard elements of $D_{\tuple q}$ form a proper subset of $\ul{\mathbb Q}$ and are harder to describe as they depend on the choice of enumeration). 
 We can repeat this construction with more than one   sequence, taking one unary predicate for each. Any structure obtained by discrete overspill is \tame\ by  Remark~\ref{R:tame}. 
%

\begin{example}\label{E:nontotal}
 Now, if we take two unary predicates, representing, say, the sequence of prime numbers $\tuple p$ and the sequence of powers of two $\tuple t$, then in $\overspill{\tuple p,\tuple t}$, it seems very unlikely that the discrete sets $D_{\tuple p}$ and $D_{\tuple t}$  are comparable. For if they were, they would have to be definably isomorphic as they have the same ultra-Euler \ch\ (equal to $\diag$, the ultraproduct of the diagonal sequence $(n)_n$). In fact, if, instead, $\tuple p$ enumerates a computable set in $\nat$  and $\tuple t$  a non-computable one, then such a definable isomorphism restricted to $\nat$ would be a computable isomorphism between these two sets, which is of course impossible. It is easy to combine these two unary sets into a single one, by letting $a_{2n}:=p_n$  and   $a_{2n-1}:=-t_n$, so that  then $D_{\tuple a}\cap (\ul {\real})_{\leq 0}=D_{\tuple t}$ and $D_{\tuple a}\cap (\ul {\real})_{\geq 0}=D_{\tuple p}$, giving an example of a single discrete overspill $\overspill{\tuple a}$ with non-total order.  Nonetheless, as we will see in Theorem~\ref{T:totorddisc} below, the failure is due to a missing isomorphism, and so,  in an appropriate expansion (no longer computable of course), both sets become definably isomorphic.
 \end{example}


\begin{theorem}[Euler O-minimality Criterion]\label{T:ominul}
An ultra-o-minimal structure   $\ul{\mathcal M}$, given as the ultraproduct of o-minimal structures  $\mathcal M_i$, is o-minimal, if for each formula $\varphi$ without parameters, there exists an $N_\varphi\in\nat$ such that $\norm{\euler{M_i}\varphi}\leq N_\varphi$.
\end{theorem}
\begin{proof}
Let $\ul Y\sub \ul M$ be definable, say, given as the fiber of a $\emptyset$-definable subset $\ul X\sub \ul M^{1+n}$ over a tuple $\ul{\tuple b}$. Let $X_i\sub M_i^{1+n}$ be the corresponding $\emptyset$-definable subset, and choose $\tuple b_i$ in $M_i$ with ultraproduct $\ul{\tuple b}$, so that $\ul Y$ is the ultraproduct of the $Y_i:=\fiber{X_i}{\tuple b_i}{}$. By the proof of Theorem~\ref{T:plane} (which in the o-minimal case does yield a finite cell decomposition), we can decompose each $X_i$ as a disjoint  union of $\emptyset$-definable subsets  $X_i^{(e)}$ consisting of the union of all $e$-cells in a cell decomposition of $X_i$. In fact, this proof can be carried out in the theory $\Ded$, so that it holds in any \omin\ structure $\mathcal M$ uniformly. For instance, assuming $n=2$ and    $X=\varphi(\mathcal M)$,  then $X^{(2)}$ consists exactly of all interior points that do not lie on a vertical fiber containing some node of $\partial X$, whereas $X^{(0)}$ consists of all nodes of $\partial X$ that belong to $X$, and $X^{(1)}$ of all remaining points. And so we can find formulae $\varphi^{(e)}$ that define in each model $\mathcal M\models\Ded$ the sets $X^{(e)}$, for $e\leq n+1$. Now, since each $X_i^{(e)}$ is a disjoint union of $e$-cells, its Euler \ch\ is equal to $(-1)^eN_{i,e}$, where $N_{i,e}$ is    the number of $e$-cells in the decomposition. By assumption (applied to the formula $\varphi^{(e)}$), this Euler \ch\ is bounded in absolute value, whence so are the $N_{i,e}$, that is to say, there exist $N_e\in\nat$ such that $N_{i,e}<N_e$ for all $i$. But then the  fiber $\fiber{X_i^{(e)}}{\tuple b_i}{}$ admits a decomposition in at most $N_e$ cells. Since the union of the latter for all $e$ is just $Y_i$, we showed that there is a uniform bound on the number  of cells, that is to say, open intervals and points,    in a decomposition of $Y_i$. Since this is now first-order expressible, $\ul Y$ too is a finite union of intervals, as we needed to show. 
\end{proof}

\section{Expansions of \omin\ structures}\label{s:ominexp}
Since an expansion by  definable sets does not alter the collection of definable sets, we immediately have:

\begin{lemma}\label{L:expdef}
If $\mathcal M$ is \omin\ and $X\sub M^n$ is definable, then $(\mathcal M,X)$ is again \omin. \qed
\end{lemma}

So we ask in more generality, what properties does a subset of an \omin\ structure need to have in order for the expansion to be again \omin? Let us call such a subset \emph\omin\ (or, more correctly, \emph{$\mathcal M$-\omin} as this depends on the surrounding structure), where we just proved that definable subsets are. 

\begin{corollary}\label{C:imagomin}
The image of an \omin\ subset   under a definable map is again \omin, and so is its complement, its closure, its boundary, and its interior. More generally, any set definable from an \omin\ set is again \omin.
\end{corollary}
\begin{proof}
It suffices to prove the last assertion. 
Let $X$ be an \omin\ subset of an \omin\ structure $\mathcal M$.   Since $(\mathcal M,X)$ is \omin, any set definable in $(\mathcal M,X)$ is \omin\  (in the expansion, whence also in the reduct) by Lemma~\ref{L:expdef}. 
\end{proof}

To define a weaker isomorphism relation, we introduce the following notation. Let $X$ be a definable subset in a structure $\mathcal M$, say, defined by the formula (with parameters) $\varphi$, that is to say, $X=\varphi(\mathcal M)$.  If $\mathcal N$ is an elementary extension of $\mathcal N$, then we set $\defext NX:=\varphi(\mathcal N)$, and call it the \emph{definitional extension} of $X$ in $\mathcal N$.

Let us call two definable subsets $X$ and $Y$  of an \omin\ structure $\mathcal M$ \emph{\omin{ally} isomorphic}, denoted $X\equiv Y$, if their definitional extensions have the same ultra-Euler measure in every ultra-o-minimal elementary extension $\mathcal M\preceq \mathcal N$, that is to say, if  $\eulerm N{\defext NX}=\eulerm N{\defext NY}$. It is easy to see that this constitutes an equivalence relation on definable subsets.


\begin{proposition}\label{P:ominiso}
In an \omin\ expansion  $\mathcal M$ of an ordered field, if two definable subsets $X$ and $Y$ are \omin{ally} isomorphic, then there exists an \omin\ expansion of $\mathcal M$ in which they become definably isomorphic.
\end{proposition}
\begin{proof}
Suppose $X$ and $Y$ are \omin{ally} definable, and let $\mathcal N$ be some ultra-o-minimal elementary extension of $\mathcal M$, given as the ultraproduct of o-minimal structures $\mathcal N_i$. Let $X_i$ and $Y_i$ be  $\mathcal N_i$-definable subsets   with respective ultraproducts   $\defext N X$ and $\defext N Y$. Since  by  Proposition~\ref{P:dimulomin} dimension is definable, $\defext N X$ and $\defext N Y$ have the same dimension, whence so do almost each $X_i$ and $Y_i$ by \los. By assumption, they have also the same Euler \ch\ for almost all $i$, so that they are definably isomorphic by   Theorem~\ref{T:Euldim}.  Hence, there exists for almost all $i$, a definable isomorphism  $f_i\colon X_i\to Y_i$. Let $\ul\Gamma$ be the ultraproduct of the graphs $\Gamma(f_i)$, so that $(\mathcal N,\ul\Gamma)$, whence also $(\mathcal M,\Gamma)$, is \omin, where $\Gamma$ is the restriction of $\ul\Gamma$ to $\mathcal M$. Moreover, by \los, $\ul \Gamma$ is the graph of a bijection $\defext N X\to \defext N Y$, and hence its restriction $\Gamma$ is the graph of a bijection $X\to Y$, proving that $X$ and $Y$ are definably isomorphic in $(\mathcal M,\Gamma)$. 
%
\end{proof}

I do not know whether the converse is also true: if $X$ and $Y$ are definably isomorphic in some \omin\ expansion $\mathcal M'$, are they \omin{ally} isomorphic?  They   will have the same Euler \ch\ in any (reduct of an) ultra-o-minimal elementary extension of $\mathcal M'$ by essentially the same argument, but what about ultra-o-minimal elementary extensions of $\mathcal M$ that are not such reducts? A related question is in case $\mathcal M$ itself is already ultra-o-minimal, if two sets have the same Euler \ch, do their definitional extensions also have  the same Euler \ch\ in an ultra-o-minimal elementary extension? This would follow if Euler \ch\ was definable, but at the moment, we can only prove a weaker version (see Theorem~\ref{T:defrk}).  Before we address these issues, we prove a result  yielding non-trivial examples of \omin{ally} isomorphic sets that need not be definably isomorphic.

\begin{corollary}\label{C:ominiso}
In an \omin\ expansion  $\mathcal M$ of an ordered field, if two definable subsets $X$ and $Y$ have the same dimension and the same class in $\grot M$, then they are \omin{ally} isomorphic. 
\end{corollary}
\begin{proof}
By Lemma~\ref{L:grequal}, there exists a definable subset $Z$ such that $X\sqcup Z$ and $Y\sqcup Z$ are definably isomorphic. 
Let $\mathcal M\preceq\mathcal N$ be an ultra-o-minimal elementary extension. Hence $\defext NX\sqcup \defext N Z$ and $\defext N Y\sqcup \defext N Z$ are definably isomorphic, and therefore 
$$
\euler N{\defext N X} +\euler N{ \defext N Z}=\euler N{\defext N X\sqcup \defext N Z}=\euler N{\defext N Y\sqcup \defext N Z}=\euler N{\defext N Y} +\euler N{ \defext N Z}
$$
showing that $\defext N X$ and $\defext N Y$ have the same ultra-Euler \ch, as we needed to show.
\end{proof}

\subsection{Contexts and \visoism{s}}
To overcome the difficulties alluded to above, we must make our definitions context-dependable in the following sense.  
Given an \omin\ structure $\mathcal M$, by a \emph{context} for $\mathcal M$, we mean an ultra-o-minimal structure $\mathcal N$ that contains $\mathcal M$ as an elementary substructure (which always exists by Corollary~\ref{C:ulomin}). An expansion $\mathcal M'$ of $\mathcal M$ is then called \emph{permissible} (with respect to the context $\mathcal N$), if $\mathcal N$ can be  expanded to a context $\mathcal N'$, that is to say, $\mathcal M'\preceq\mathcal N'$ and $\mathcal N'$ is again ultra-o-minimal. If $\mathcal M$ itself is ultra-o-minimal, then we may take it as its own context, but even in this case, not every expansion will be permissible, as it may fail to be an ultraproduct.

From now on, we fix an \omin\ structure $\mathcal M$ and a context $\mathcal N$. We  define a (context-dependable)   \emph{Euler \ch} $\euler M\cdot$ (or, simply $\chi$) on $\mathcal M$ by restricting the ultra-Euler \ch\ of $\mathcal N$, that is to say,  by setting $\euler {{}}X:=\euler N{\defext NX}$, for any definable subset   of  $\mathcal M$, and we define similarly its \emph{Euler measure} $\eulerm {{}}X:=(\op{dim}(X),\euler{{}}X)$. We say that two definable subsets are \emph{\viso}, if there exists a permissible expansion of $\mathcal M$ in which they become definably isomorphic.  In particular,  two definable subsets that are \omin{ally} isomorphic are also \viso, but the converse is unclear. We can now prove an \omin\ analogue of Theorem~\ref{T:Euldim}.

\begin{theorem}\label{T:viso}
In an \omin\ expansion $\mathcal M$ of an ordered field, two definable subsets  are \viso\ \iff\ they have the same Euler measure.
\end{theorem} 
\begin{proof}
One direction is proven in the same way as Proposition~\ref{P:ominiso}, so assume $X$ and $Y$ are \viso\ definable subsets. By assumption, $\mathcal M\preceq\mathcal N$ expands into   \omin\ structures $\mathcal M'\preceq \mathcal N'$,  with $\mathcal N'$ again ultra-o-minimal, such that $X$ and $Y$ are $\mathcal M'$-definably isomorphic. Let $\mathcal N'$ be the ultraproduct of  o-minimal structures $\mathcal N_i'$. Since $\defext {N'}X$ and $\defext{N'}Y$ are  definably isomorphic, so are almost all $X_i$ and $Y_i$, where $X_i$ and $Y_i$ are $\mathcal N_i'$-definable subsets with respective ultraproducts $\defext {N'}X$ and $\defext{N'}Y$. In particular,  $X_i$ and $Y_i$ have the same Euler measure for almost all $i$, by Theorem~\ref{T:Euldim}. Hence $\defext {N'}X$ and $\defext{N'}Y$ have the same ultra-Euler measure, by Proposition~\ref{P:dimulomin}. Since both invariants remain the same in the reduct $\mathcal N$, elementarity then yields $\eulerm{{}}X=\eulerm{{}}Y$.
\end{proof} 

\begin{remark}\label{R:totorddisc}
Any two  subsets   given by discrete overspill  with respect to    non-repeating sequences  (see Example~\ref{E:discoverspill}) are \viso, since they both have Euler \ch\ $\diag$. 
\end{remark}

\subsection{O-finitism}
In the terminology of \cite{ForTame},  the definable  discrete sets in an \omin\ structure are exactly the \emph{pseudo-finite} sets. As we already mentioned in the introduction,  in the \omin\ context, discrete sets play the role of finite sets, and so we briefly discuss the first-order aspects of this assertion.  

Given  a (non-empty) collection of $L$-structures $\mathfrak K$, and a subset $X\sub M^n$ in some $L$-structure $\mathcal M$, we say that $X$ is \emph{$\mathfrak K$-finitistic}, if $(\mathcal M, X)$ satisfies every $L(\tt U)$-sentence $\sigma$ which  holds in every expansion $(\mathcal K,F)$ of a structure $\mathcal K\in\mathfrak K$ by a finite set $F\sub K^n$.  In case $\mathfrak K$ is the collection of o-minimal structures,  we call $X\sub M^n$ \emph{\ofin}. Applying the definition just to $L$-sentences $\sigma$ (not containing the predicate $\tt U$, so that $(\mathcal K,F)\models\sigma$ \iff\  $\mathcal K\models\sigma$), we see that $\mathcal M$ is then necessarily \omin. Put differently, an \ofin\ set in an \omin\ structure is a model of \emph{o-finitism}, that is to say, of the theory of a finite set in an o-minimal structure.   By Propositions~\ref{P:Dedlim}, \ref{P:locdefdisc}, and \ref{P:DPP}, we have:

\begin{corollary}\label{C:ofindisc}
In an \omin\ structure, an \ofin\ set is discrete, closed,  bounded, and locally definable, every non-empty intersection with an open interval has a maximum and a minimum, and every injective, definable self-map on it is an isomorphism. \qed
\end{corollary}
It seems unlikely that these properties characterize fully o-finitism. A complete axiomatization of o-finitism would be of interest in view of the following results.

\begin{theorem}\label{T:ofinexp}
A subset $X$ of an \omin\ structure $\mathcal M$ is \ofin\ \iff\ it is discrete and  \omin. In particular, any definable, discrete subset in an \omin\ structure is \ofin.
\end{theorem}
\begin{proof}
Assume first that $X$ is \ofin, whence discrete by Corollary~\ref{C:ofindisc}. We have to show that given an $L(\tt U)$-sentence $\sigma$ holding true in every o-minimal $L(\tt U)$-structure, then $(\mathcal M,X)\models\sigma$. Let $\mathcal K$ be an o-minimal structure and let $F\sub K^n$ be a finite subset. Hence $(\mathcal K,F)$ is also o-minimal and therefore satisfies $\sigma$. Since this holds for all such expansions, $\sigma$ is true in  $(\mathcal M,X)$ by o-finitism, as we needed to show.

Conversely, suppose $X\sub M^n$ is discrete and \omin, that is to say, $(\mathcal M,X)$ is \omin. To show that $X$ is \ofin, let $\sigma$ be a sentence true in every expansion $(\mathcal K,F)$ of an o-minimal structure $\mathcal K$ by a finite subset $F\sub K^n$. Consider the disjunction $\sigma'$ of $\sigma$ with the sentence expressing that  the set defined by $\tt U$ is not discrete. Hence $\sigma'$  is true in any o-minimal expansion $(\mathcal K,Y)$. Since $X$ is \omin, this means that $(\mathcal M,X)\models \sigma'$, and since $X$ is discrete, this in turn implies that $\sigma$ is true in $(\mathcal M,X)$, as we needed to show. The last assertion then follows from Lemma~\ref{L:expdef}.
\end{proof}

Let us call a subset of an ultra-o-minimal structure \emph{ultra-finite}, if it is the ultraproduct of finite subsets (such a set may fail to be definable, since the definition in each component  may not be uniform). An ultra-finite set  is \ofin. As for the converse, if an \ofin\ set is an ultraproduct (a so-called \emph{internal} set), then it must be \ofin, but what for external sets? In any case, we have:

\begin{theorem}\label{T:ulominic}
A subset $X\sub M^k$ of an \omin\ structure $\mathcal M$ is \ofin\ \iff\ there exists an elementary extension $\mathcal M\preceq \mathcal N$ with $\mathcal N$ ultra-o-minimal and an ultra-finite subset $Y\sub N^k$, such that $X=Y\cap M^k$.
\end{theorem}
\begin{proof}
Suppose $\mathcal N$ and $Y$ have the stated properties, and let   $\mathcal N_i$  be     o-minimal structures and $Y_i\sub N_i^k$  finite subsets, so that $\mathcal N$ and $Y$ are their respective ultraproducts. Since $(\mathcal N_i,Y_i)$ is again o-minimal, their ultraproduct $(\mathcal N,Y)$ is \omin. Since $(\mathcal M,X)$ is then an elementary substructure, the latter is also \omin. Moreover, since  $Y$ is   discrete, so must $X$ be, and hence $X$ is \ofin\ by Theorem~\ref{T:ofinexp}. Conversely, by the same theorem, if $X$ is \ofin, then $(\mathcal M,X)$ is \omin. Hence there exists an  elementary extension  $(\mathcal N,Y)$ which is ultra-o-minimal as an $L(\tt U)$-structure by Corollary~\ref{C:ulomin}. Write  $(\mathcal N,Y)$   as an ultraproduct of o-minimal structures $(\mathcal N_i,Y_i)$. Since $X$ is discrete, so must $Y$ be by elementarity, whence so are almost all $Y_i$ by \los. The latter means that almost all are in fact finite, showing that $Y$ is ultra-finite, and the assertion follows since $X=Y\cap M^k$. 
\end{proof}

%

Next, we give a criterion for a subset $Y\sub M$ to be \omin. 
By Theorem~\ref{T:onevar}, its boundary $\partial Y$ should be discrete, and $Y^\circ=Y\setminus\partial Y$ should be a disjoint union of open intervals. Given an arbitrary set $Y\sub M$, define its \emph{enhanced boundary} $\enh Y$ as the set consisting of the following pairs: $(d,0)$ if $d\in  Y$,   $(d,1)$ if $d^+$ belongs to $Y$, and $(d,-1)$ if $d^-$ belongs to $Y$, where $d$ runs over all boundary points of $Y$. An enhanced boundary cannot have fibers with three points and its projection is the ordinary boundary $\partial Y$. If $Y$ is \omin, then $\enh Y$ must satisfy some extra conditions: it must be bounded, discrete and closed, and, by \eqref{eq:event}, if $(d,1)$ belongs to it, then so must $(d',-1)$, where $d'$ is the immediate successor of $d$ in $\partial Y$.

\begin{theorem}\label{T:enhbd}
A one-variable subset $Y$ in an \omin\ structure $\mathcal M$ is \omin\ \iff\ its enhanced boundary $\enh Y$ is \ofin\ and its interior is a disjoint union of open intervals. 
\end{theorem}
\begin{proof}
Suppose $Y$ is \omin, so that $Y^\circ$ is a disjoint union of open intervals. Since $\enh Y$ is definable from $Y$, it too is \omin\ by  Corollary~\ref{C:imagomin}, whence \ofin\ by Theorem~\ref{T:ofinexp}. To prove the converse, let $D:=\partial Y=\pi(\enh Y)$, a bounded, closed, discrete set, and let $l$ be its minimum. Define $X\sub M$ as  the set of all $x\in M$  such that one of the following three conditions holds
\begin{enumerate}
\item\label{i:incase} $x\in D$;
\item\label{i:leftcase} $x>l$ and $(d,1)\in\enh Y$, where $d=\max D_{<x}$;
\item\label{i:unbcase}   $x<l$  and $(l,-1)\in\enh Y$.
\end{enumerate}
Since $X$ is definable from $\enh Y$, it is \omin\ by Corollary~\ref{C:imagomin}. Remains to show that $X=Y$. It follows from \eqref{i:incase} that  $X\cap \partial Y=Y\cap \partial Y$, so that it suffices to show that $X^\circ=Y^\circ$. Therefore, we may as well assume from the start that $Y$ is open.  Write $Y=\sqcup_n I_n$ as a disjoint union of open intervals, and let $\oo ab$ one of the $I_n$ (we leave the unbounded case  to the reader, for which one needs \eqref{i:unbcase}). In particular, $a\in D$ and $a^+$ belongs to $Y$, so that $(a,1)\in\enh Y$. By \eqref{i:leftcase}, the entire interval $\oo ab$ lies in $ X$, whence so does the whole of $Y$. Conversely, if $x\in X$, let $d:=\max D_{<x}$, so that $(d,1)\in\enh Y$. Hence $d^+$ belongs to $Y$, and so $d$ must be an endpoint of one of the $I_n$. The other endpoint must be bigger than $d$, and hence bigger than $x$, showing that  $x\in I_n\sub Y$.
\end{proof}
%
%
\subsection{The virtual \gr}
We fix again an \omin\ structure $\mathcal M$ and a context $\mathcal N$.
We can use   \visoism{s} instead of definable isomorphisms in the definition of the zero-dimensional or the full  \gr, that is to say,  the quotient modulo the scissor relations of the free Abelian group on \visoism\   classes of respectively all discrete, definable subsets,  and of all definable subsets  yield the \emph{virtual \gr{s}} $\grotomin M0$ and $\grotomin M{}$ respectively.
We have   surjective \homo{s} $\grotth M0\to \grotomin M0$ and $\grot M\to \grotomin M{}$. 

\begin{corollary}\label{C:omingr}
Given an \omin\ expansion  $\mathcal M$ of an ordered field, there exist embeddings  $\grotomin M0{}\sub \grotomin M{}\into \ul\zet$, where $\ul\zet$ is the ring of non-standard  integers in the given context.
\end{corollary}
\begin{proof}
Since the Euler \ch\ vanishes on any  scissor relation, it induces by  
Theorem~\ref{T:viso}  a \homo\ $\chi\colon  \grotomin M{}\to \ul\zet$. By the same result, its restriction to $\grotomin M0$ is injective. To see that $\chi$ is everywhere injective, assume $\euler{{}}X=\euler{{}}Y$ for some definable subsets $X$ and $Y$. If they have the same dimension, then they are \viso, again by Theorem~\ref{T:viso}. So assume $X$ has dimension $d\geq 1$ and $Y$ has lesser dimension. Let $U$ be the difference of a $d$-dimensional box minus a $(d-1)$-dimensional subbox, so that  in particular $\class U$ vanishes, whence also $\euler{{}}U$. As $X$ and $Y\sqcup U$ now have the same   Euler measure, they are \viso\ by Theorem~\ref{T:viso}, and hence $\class X=\class Y+\class U=\class Y$ in $\grotomin M{}$, as we needed to show. The injectivity of  $\grotomin M0\to \grotomin M{}$ is then also clear. 
\end{proof}

In particular, if $\mathcal M$ is moreover \tame, then we have an equality of   virtual  \gr{s} $\grotomin M0= \grotomin M{}$ by Corollary~\ref{C:defcell}. We can also extend the partial order on $\disc M$ to a total order on $\discomin M$, the set of \visoism\ classes of  definable, discrete subsets. First, given definable subsets $X$ and $Y$, we say that $X\leq Y$ (or, if we want to emphasize the context, $X\leq_{\mathcal N}Y$), if  $X\preceq_{\mathcal M'}Y$  in some  permissible \omin\ expansion $\mathcal M'$ of $\mathcal M$. Clearly, if $X\preceq Y$, then $X\leq Y$. The following two results are the \omin\ analogues of Theorem~\ref{T:Euldim}.
%

\begin{theorem}\label{T:totorddisc}
In an \omin\ structure $\mathcal M$,  two definable, discrete subsets $F$ and $G$   satisfy $F\leq G$ \iff\   $\euler{{}}{ F}\leq \euler{{}}{ G}$.  In particular, $\leq$ is a total order on $\discomin M$.
\end{theorem}
\begin{proof}
Suppose first that $\euler{{}}{  F}\leq \euler{{}}{ G}$. Write $\mathcal N$ as the ultraproduct of o-minimal structures $\mathcal N_i$, and let $F_i$ and $G_i$ be finite sets with respective ultraproducts the definitional extensions $\defext NF$ and $\defext NG$ of $F$ and $G$ respectively. Since $\euler N{\defext NF}\leq \euler N{\defext NG}$, the cardinality of $F_i$ is at most that of $G_i$, for almost all $i$. In particular, there exists an injective map $F_i\to G_i$ for almost all $i$.   Let $\ul\Gamma$ be the ultraproduct of the graphs of these maps $F_i\to G_i$. Hence $\ul\Gamma$ is ultra-finite and therefore its restriction $\Gamma$ to $\mathcal M$ is \ofin\ by Theorem~\ref{T:ulominic}, whence \omin\ by Theorem~\ref{T:ofinexp}. By \los\ and elementarity, $\Gamma$ is the graph of an injective map $F\to G$, showing that $F\preceq_{(\mathcal M,\Gamma)} G$. Since $(\mathcal M,\Gamma)$ is permissible,   $F\leq G$. The converse goes along the same lines: suppose $F\preceq_{\mathcal M'}G$, for some permissible \omin\ expansion $\mathcal M'$ of $\mathcal M$. By definition, there is an ultra-o-minimal expansion    $\mathcal N'$ of $\mathcal N$ with $\mathcal M'\preceq \mathcal N'$. Since  $\defext N F\preceq_{\mathcal N'}\defext N G$, we have $\euler{{}} F=\euler {N'}{\defext N F}\leq \euler {N'}{\defext N G}=\euler{{}}G$. 
%
\end{proof}

\begin{proposition}\label{P:preorddim}
In an \omin\ expansion $\mathcal M$ of an ordered field, we have $X\leq Y$ \iff\ $\op{dim}(X)\leq \op{dim}(Y)$, for $X$ and $Y$ definable subsets with $\op{dim}(Y)>0$.
\end{proposition} 
\begin{proof}
The direct implication is clear. For the converse, by definability of dimension, we may pass to the context of $\mathcal M$ and therefore already assume $\mathcal M$ is ultra-o-minimal, given as the ultraproduct of o-minimal structures $\mathcal M_i$. Let $X_i$ and $Y_i$ be definable subsets in $\mathcal M_i$ with respective ultraproducts $X$ and $Y$. By \los, $\op{dim}(X_i)\leq \op{dim}(Y_i)$, and hence $X_i\preceq Y_i$, by Theorem~\ref{T:Euldim}, for almost all $i$. Let $f_i\colon X_i\to Y_i$ be a definable injection and let $\ul\Gamma$ be the ultraproduct of the graphs $\Gamma(f_i)$. Since each $(\mathcal M_i,\Gamma(f_i))$ is again o-minimal, $(\mathcal M,\ul\Gamma)$ is ultra-o-minimal and hence in particular a permissible expansion. Since $\ul \Gamma$ is the graph of an injective map by \los, $X\preceq_{(\mathcal M,\ul\Gamma)}Y$, as we needed to show.
\end{proof} 

In particular, any definable, discrete subset is virtually univalent.

\begin{corollary}\label{C:OPP}[Virtual Pigeonhole Principle]
Given   an \omin\   structure $\mathcal M$,  two definable,  discrete  subsets $D$ and $E$         are \viso\ \iff, for some definable subset $X$,  the sets $D\sqcup X$ and $E\sqcup X$ are \viso, \iff\ $\class D=\class E$ in $\grotomin M{}$.  
\end{corollary}
\begin{proof}
One direction in the first equivalence is immediate, so assume  $D\sqcup X$ and $E\sqcup X$ are \viso. Passing to a permissible \omin\ expansion, we may assume that they are already definably isomorphic, say, by an isomorphism $f\colon D\sqcup X\to E\sqcup X$. By totality (Theorem~\ref{T:totorddisc}), we may assume that  $E\leq  D$, and hence after taking another permissible \omin\ expansion, and replacing $E$ with an isomorphic image, we may even assume that $E\sub D$. Therefore,  the composition of $f$ and the inclusion $E\sqcup X\sub D\sqcup X$ is a map with co-discrete image, and hence is surjective by \eqref{eq:DPP}. However, this can only be the case  if  $E=D$, as we needed to show. The last equivalence is now just Lemma~\ref{L:grequal}.
\end{proof}

\begin{corollary}\label{C:ordomingr}
For an \omin\ structure $\mathcal M$, its zero-dimensional, virtual \gr\ $\grotomin M0$  is an ordered ring with respect to $\leq$.
\end{corollary}
\begin{proof}
Every element in $\grotomin M0$ is of the form $\class A-\class B$, for some definable, discrete subsets $A$ and $B$ in the \omin\ structure $\mathcal M$. Therefore, for definable, discrete subsets $A_i$  and $B_i$, with $i=1,2$, we set $\class {A_1}-\class{B_1}\leq \class{A_2}-\class{B_2}$ \iff\ 
\begin{equation}\label{eq:ord}
A_1\sqcup B_2\leq A_2\sqcup B_1.
\end{equation} 
To see that this is well-defined, suppose $\class {A_i}-\class{B_i}=\class {A'_i}-\class{B'_i}$, for $i=1,2$ and definable, discrete subsets $A'_i$  and $B'_i$. Therefore,  $\class{A_i\sqcup B_i'}=\class{A_i'\sqcup B_i}$, whence $A_i\sqcup B_i'$ and $A_i'\sqcup B_i$ are \viso\ by Corollary~\ref{C:OPP}. We have to show that assuming \eqref{eq:ord}, the same inequality holds for the accented sets. Taking the disjoint union with $B_1'\sqcup B_2'$ on both sides of \eqref{eq:ord}, yields inequalities
\begin{align*}
(A_1\sqcup B_1')\sqcup B_2\sqcup B_2'&\leq (A_2\sqcup B_2')\sqcup B_1\sqcup B_1'\\
(A_1'\sqcup B_1)\sqcup B_2\sqcup B_2'&\leq (A_2'\sqcup B_2)\sqcup B_1\sqcup B_1'\\
(A_1'\sqcup B'_2)\sqcup (B_1\sqcup B_2)&\leq (A_2'\sqcup B_1')\sqcup (B_1\sqcup B_2)
\end{align*}
which by another application of Corollary~\ref{C:OPP} then gives $A_1'\sqcup B'_2\leq A_2'\sqcup B_1'$, as we needed to show. It is now easy to check that $\leq$ makes $\grotomin M0$ into a totally ordered ring.
\end{proof}

\begin{corollary}\label{C:ofincut}
Every \ofin\ subset defines a   cut in $\discomin M$. In particular, we can put a total pre-order on the collection of \ofin\ subsets.
\end{corollary}
\begin{proof}
Let $F$ be an \ofin\ subset of an \omin\ structure $\mathcal M$ and let $D\in\discomin M$ be arbitrary. Since $(\mathcal M,F)$ is \omin\ by Theorem~\ref{T:ofinexp}, we can compare $D$ and $F$ in $\disc{M,F}$ by Theorem~\ref{T:totorddisc}. If $G$ is another \ofin\ subset, then we set $F\leq G$ \iff\ the lower cut in $\disc M$ determined by $F$ is contained in the lower cut of $G$.
\end{proof}

A note of caution: even if $F\leq G$ and $G\leq F$,  for $F$ and $G$ \ofin\ subsets, they need not be \viso. For instance,  taking $D$ as in Example~\ref{E:ulomin}, it is an \ofin\ subset of $\ul{\real}$, and since $\disc{{\ul{\real}}}$ is just $\nat$ by o-minimality, its cut is  $\infty$. However, $D\setminus\{\diag\}$ determines the same cut, whence $D\leq D\setminus\{\diag\}\leq D$, but we know that they cannot be definably isomorphic  in any \omin\ expansion by \eqref{eq:DPP}. In fact, it is not clear whether  two given \ofin\ subsets  live in a common \omin\ expansion, and therefore can be compared directly.  This is also why we cannot (yet?) define a \gr\  on \ofin\ subsets. 

\subsection{Discretely valued Euler \ch{s}}
In order to calculate the zero-dimensional virtual \gr, we introduce a new type of Euler \ch. Fix an \omin\ structure $\mathcal M$ and a context $\mathcal N$, and let $D$ be a definable, discrete subset. In this section, we will always view $D$ in its lexicographical order $\leq_{\text{lex}}$ (or, when there is no risk for confusion, simply denoted $\leq$).

\begin{corollary}\label{C:defsub}
Any definable subset of a definable, discrete subset $D$ in an \omin\ structure $\mathcal M$ is \viso\ to an initial segment $D_{\leq a}$. 
\end{corollary}
\begin{proof}
The set of initial segments is a maximal chain in $\discomin M$, since any two consecutive subsets in this chain differ by a single point. Hence, any definable subset $E\sub D$ must be a member of this chain up to \visoism. 
\end{proof}

Clearly, such an $a$ must be unique, and so, given a non-empty definable subset $E\sub D$, we let $\discrk DE$   be the unique $a$ such that $E$ is \viso\ with $D_{\leq a}$. We add a new symbol $\varnothing$ to $D$ and  set $\discrk D\emptyset:=\varnothing$. For  definable subsets $E_1,E_2\sub D$, we have $E_1\leq E_2$ \iff\ $\discrk {D}{E_1}\leq \discrk {D}{E_2}$. Given  a definable map $g$ with domain $D$,   we can define by Lemma~\ref{L:imageorder} its \emph{rank} as $\rk {}g:=\discrk {D}{g(D)}$. A map is constant \iff\ its rank is minimal (that is to say, equal to the minimum of its domain). By   \eqref{eq:DPP}, we immediately have:

\begin{corollary}\label{C:defdiscmaprk}
In an \omin\ structure, a definable map with discrete domain  is injective \iff\ its rank is maximal (that is to say, equal to the maximum of its domain).\qed
\end{corollary}

Let $D\sub M$ be definable and discrete, with minimal element $l$ and maximal element $h$. For each $n$, we  view the Cartesian power $D^n$ as a definable subset of $D^{n+1}$ via the map $\tuple a\mapsto (l,\tuple a)$. We also need to take into consideration the empty set, and so we define $\varnothing$ to be lower than any element in any $D^n$, and we  let $D^\infty$ be the direct limit of the ordered sets $D^n\cup \{\varnothing\}$. Under this identification, the elements of $D^n\cup \{\varnothing\}$  form an initial segment   in $D^{n+1}\cup \{\varnothing\}$ with respect to the lexicographical ordering. In particular, if $E\sub D^n$ is a non-empty definable subset, then $\discrk {D^n}E=\discrk {D^{n+1}}{E'}$, where $E'$ is the image  of $E$ in $D^{n+1}$. After identification therefore, we will view $\discrk {D^n}E$ simply as an element of $D^\infty$, and we just denote it $\discrk DE$.  More generally, given an arbitrary definable subset $X\sub M^n$, we define its \emph{$D$-valued Euler \ch} (or, simply \emph{Euler \ch}) $\discrk DX:=\discrk {D^n}{X\cap D^n}$. 

We define an addition and a multiplication on $D^\infty$ as follows. First, let us define the \emph{disjoint union} $A\sqcup B$ of two definable subsets $A,B\sub M^n$ as the definable subset in $M^{n+1}$ consisting of all $(a,l)$ and $(b,h)$ with $a\in A$ and $b\in B$. For $a\in D^\infty$, we set $a\oplus\varnothing=\varnothing\oplus a=a$ and  $a\otimes\varnothing=\varnothing\otimes a=\varnothing$. For the general case, assume $a,b\in D^n$,   and    let $a\oplus b$ be the Euler \ch\ of the disjoint union $(D^n)_{\leq a}\sqcup(D^n)_{\leq b}\sub D^{n+1}$, and let   $a\otimes b$ be the Euler \ch\ of the Cartesian product $(D^n)_{\leq a}\times (D^n)_{\leq b}\sub D^{2n}$. One verifies that both operations are independent of the choice of $n$, making $D^\infty$ into a commutative semi-ring, where the zero for $\oplus$ is $\varnothing$, and where the unit for $\otimes$ is  $l$, the minimum of $D$. We even can define a subtraction: if $a\leq b$ in $D^\infty$, then we define $b\ominus a$ as the Euler \ch\ of $D^n\cap \oc a b$, where $n$ is sufficiently large so that $a,b\in D^n$. This allows us to define  the Grothendieck group generated by $(D^\infty, \oplus)$, defined as  all pairs $(x,y)$ with $x,y\in D^\infty$ up to the equivalence $(x,y)\sim (x',y')$ \iff\ $x\oplus y'=x'\oplus y$; the induced  commutative ring will be denoted   $\rankgr D$, and called the ring of \emph{$D$-integers}.

To turn this into a genuine Euler \ch, recall the construction of the \emph{induced structure} $\ind D$ on a subset $D\sub M$ of a first-order structure: for each definable subset $X\sub M^n$, we have a predicate defining in $\ind D$ the subset $M\cap D^n$. If $\mathcal M$ is an ordered structure, then so is $\ind D$. If $D$ is definable, then we have an induced \homo\ of \gr{s} $\grot{\ind D}\to \grotth M0$. If instead of definable isomorphism, we take \visoism, we get the virtual variant $\grotomin{\ind D}{}\to \grotomin M0$. By the Virtual Pigeonhole Principle (Corollary~\ref{C:OPP}), this latter \homo\ is injective. To discuss when they are isomorphic, let us call $D$   \emph{power dominant}, if for every definable, discrete subset $A$, there is some $n$ such that $A\leq D^n$.

\begin{proposition}\label{P:powdom}
In an \omin\ structure $\mathcal M$, a definable, discrete subset $D\sub M$  is power dominant \iff\ $
\grotomin {\ind D}{}\iso \grotomin M0$.
\end{proposition}
\begin{proof}
Suppose first that $D$ is power dominant and let $A$ be an arbitrary definable, discrete subset. By assumption, there exists an $n$ and a definable subset $B\sub D^n$, such that $A$ is \viso\ with $B$. Hence $\class A=\class B$ in $\grotomin M0$, proving that it lies in the image of $\grotomin{\ind D}{}\to \grotomin M0$. 

Conversely, assume that the latter map is surjective, and let $A$ be an arbitrary definable, discrete subset. Hence, there exists an $n$ and definable subsets  $E,F\sub D^n$ such that $\class A=\class E-\class F$ in $\grotomin M0$. By the   Virtual Pigeonhole Principle (Corollary~\ref{C:OPP}), this means that there is a \visoism\ $A\sqcup F\to E$. Hence the composition $A\sub A\sqcup F\to E\sub D^n$, shows that $A\leq D^n$.
\end{proof}

To study the existence of power dominant sets, let us say, for $D$ and $E$ discrete, definable subsets, that  $D\lll E$, if $D^n\leq E$ for all $n$. If neither $D\lll E$ nor $E\lll D$, then $D$ and $E$ are mutually power bounded, that is to say,  there  exist $m$ and $n$ such that $D\leq E^m$ and $E\leq D^m$, and we write $D\approx E$. Hence $\lll$ induces a total order relation on the set $\powarch M$ of $\approx$-classes of definable, discrete subsets of $\mathcal M$. The class of the empty set is the minimal element of $\powarch M$,  the class of a singleton is the next smallest element, and the class of a two-element set is the next (and consists of all finite sets). For an o-minimal structure, these are the only three classes, whereas for a proper \omin\ structure, there must be at least one more class, of some infinite set. I do not know whether $\powarch M$ is always  discretely ordered. In any case, it follows easily from the definitions that a class is maximal in $\powarch M$ \iff\ it is the class of a power dominant set. Thus,  the existence of a power dominant set corresponds to $\powarch M$ having a maximal element, which is especially interesting in view of Proposition~\ref{P:powdom} and its applications below.  I conjecture that $D$ as in Example~\ref{E:ulomin} is power dominant (and a similar property for any set obtained by discrete overspill). This would follow from the following growth conjecture in an o-minimal $L$-expansion $\mathcal R$ of $\real$: does there exist, for every formula $\varphi$ in the language $L(\tt U)$, some  $n\in\nat$, such that for any finite subset $F$, the set $\varphi(\mathcal R,F)$ defined by interpreting the unary predicate $\tt U$ by $F$, if finite, has cardinality at most $\norm F^n$.  Likewise, I conjecture that the following always produces a power dominant set: let $\mathcal M$ be o-minimal and let $D$ be \ofin, then $D$ is power dominant in the (\omin) expansion $(\mathcal M,D)$.    

\begin{theorem}\label{T:Drank}
In an \omin\ structure $\mathcal M$, every definable, discrete subset $D\sub M$  induces a ring isomorphism $\grotomin {\ind D}{}\iso \rankgr D$ by sending the class of a definable subset to its $D$-valued Euler \ch.
\end{theorem}
\begin{proof}
We already observed that the ring operations on $\rankgr D$ are invariant under \visoism. It is now easy to see that they also respect the scissor relations~\eqref{eq:sciss} in the \gr\ of $\ind D$. Surjectivity follows since every element in $\rankgr D$ is of the form $a\ominus b$ for some $n$ and some $a,b\in D^n$, and hence is the image of $\class{(D^n)_{\leq a}}-\class{(D^n)_{\leq b}}$. To calculate the kernel, we can write a general element   as $\class E-\class F$, with $E,F$ definable subsets in $\ind D$. Such an element lies in the kernel if $\discrk DE=\discrk DF$, which means that $E$ and $F$ are \viso,  whence $\class E=\class F$ in $\grotomin {\ind D}{}$.
\end{proof}

Summarizing, we have the following diagram of \homo{s} among the various \gr{s}, for $\mathcal M$ an \omin\ expansion of an ordered field:
\begin{equation} \label{eq:gr}
\xymatrix{
\grot{\ind D}\ar@{->>}[r]\ar[d]&\grotomin{\ind D}{}\ar^(.6){\sim}[r]\ar@{^{(}->}^i[d]&\rankgr D  \\
\grotth{M}0\ar@{->>}[r]\ar[d]&\grotomin{M}{0}\ar@{^{(}->}^j[d]\\
\grot{M}\ar@{->>}[r]&\grotomin{M}{}
}
\end{equation}
with $i$ an isomorphism if $D$ is power dominant by Proposition~\ref{P:powdom}, and with $j$ an isomorphism if $\mathcal M$ is \tame, by Corollary~\ref{C:defcell}, that is to say, we proved:

\begin{corollary}\label{C:gromintame}
If $\mathcal M$ is a \tame, \omin\ expansion of an ordered field admitting a definable, power dominant subset $D$, then its \omin\ \gr\ $\grotomin M{}$ is isomorphic to the ring of $D$-integers $\rankgr D$.\qed
\end{corollary}

  If we would allow classes of \ofin\  subsets in $\powarch M$, then there never is a maximal element: let $D$ be any definable, discrete subset (or even any \ofin\ subset). Take an ultra-o-minimal elementary extension $\mathcal N$, and choose $D_i\sub N_i$ such that their   ultraproduct is $\defext N D$. Assuming univalence, let $A_i\sub N_i$ be isomorphic with $D_i^i$ and let $\ul A\sub N$ be their ultraproduct. By Theorem~\ref{T:ulominic}, the restriction $\ul A\cap M$ is \ofin\ and satisfies by \los\ $D^n\leq A$ for all $n$, that is to say, $D\lll A$.  

\begin{theorem}[O-minimalism of Euler \ch{s}]\label{T:defrk}
Let $D\sub M$ be a definable, discrete subset of an \omin\ structure $\mathcal M$, and let $X\sub M^{n+k}$ be any definable subset. For each $e\in D^n$, the set of parameters $\tuple a\in M^n$ such that $\discrk D{\fiber X{\tuple a}{}}= e$ is \omin.
\end{theorem}
\begin{proof}
If $\tuple a$ does not belong to $D^n$, then the fiber $\fiber X{\tuple a}{}$ is empty, whence has Euler \ch\ $\varnothing$. As these $\tuple a$ form a definable subset, we may therefore replace $X$ by $X\cap D^{n+k}$ and assume already that $X$ is a definable subset of $D^n$.  
Let $\mathcal N$ be the context and write it as the ultraproduct of  o-minimal structures $\mathcal N_i$. Choose $D_i\sub N_i$, $e_i\in D_i^n$ and  $X_i\sub D_i^{n+k}$  with respective ultraproducts $\defext ND$, $e$,  and $\defext N X$.  For each $i$, let $F_i\sub N_i^n$ be the (finite) set of parameters for which the fiber has   the same cardinality as $(D_i^n)_{\leq e_i}$. Hence, for each $\tuple a\in F_i$, there exists  a bijection $f_{\tuple a}\colon \fiber{X_i}{\tuple a}{}\to (D_i^n)_{\leq e_i}$. Let $H_i\sub N_i^{3n}$ be the union of all $\{\tuple a\}\times \Gamma(f_{\tuple a})$, where $\tuple a$ runs over all tuples in $F_i$. Let $\ul F\sub N^n$ and $\ul  H\sub N^{3n}$ be their  ultraproduct, so that both sets are ultra-finite. By \los, for each $\tuple a\in \ul F$, the fiber $\fiber {\ul H}{\tuple a}{}$ is the graph of a bijection $\fiber{\defext N X}{\tuple a}{}\to \left((\defext N D)^n\right)_{\leq e}$. Therefore,     $F:=\ul F\cap M^n$ consists precisely of those $\tuple a\in M^n$ for which the fiber $\fiber X{\tuple a}{}$ has $D$-valued Euler \ch\ $e$ in the expansion $(\mathcal M,\ul H\cap M^{3n})$ whence in $\mathcal M$, as the former is \omin\ by Theorem~\ref{T:ulominic}. For the same reason,  $F$ is   \ofin, whence  \omin\ by Theorem~\ref{T:ofinexp}, so that  we are done.  
\end{proof}

\begin{remark}\label{R:ofinarch}
In everything in this section on Euler \ch{s}, we may, by passing to a suitable permissible expansion, even assume that $D$ is only \ofin. 
\end{remark}

\subsection{Archimedean reducts}
As before, let $D$ be  definable and  discrete with  respective minimum $l$ and maximum $h$. By \eqref{i:succdisc}, we have a successor function $\sigma:=\sigma_D$, defined on $D\setminus\{h\}$, with inverse $\inv\sigma$ defined on $D\setminus\{l\}$. Let us write $e\ll d$, if $\sigma^n(e)<d$, for all $n\in\nat$. If neither $d\ll e$ nor $e\ll d$, then $\sigma^n(d)=e$ for some $n\in\zet$, and we write $d\sim_D e$. The set of $\sim_D$-equivalence classes is totally ordered by $\ll$, and is called the \emph{Archimedean reduct} $\op{Arch}(D)$ of $D$.

\begin{theorem}\label{T:archred}
The Archimedean reduct $\op{Arch}(D)$ of a  definable, discrete subset $D$ in an \omin\ structure $\mathcal M$ is dense. 
\end{theorem}
\begin{proof}
This is clear if $D$ is finite, since then there is only one Archimedean class, so assume it is infinite. If $\op{Arch}(D)$  is not dense, there would exist $l\ll h$ in $D$  so that for no $d\in D$ we have $l\ll d\ll h$. Therefore, upon replacing $D$ with $D\cap \cc lh$, we may assume that  $\op{Arch}(D)$ consists of exactly two classes, those of $l$ and $h$ . By Corollary~\ref{C:ulomin} (or, Theorem~\ref{T:ulominic}), we can  embed $\mathcal M$ elementary in an ultra-o-minimal structure $\mathcal N$ so that $D$ is the restriction of a (definable) ultra-finite set $F$ in $\mathcal N$. Let $\mathcal N_i$ and $F_i$ be respectively o-minimal structures and finite subsets in these with ultraproduct equal to $\mathcal N$ and $F$ respectively. For each $i$, let $f_i\colon F_i\to F_i$ be the map reversing the (lexicographical) order and let $\ul\Gamma$ be the ultraproduct of the graphs of the $f_i$. Since this is an ultra-finite set, its restriction $\Gamma$ to $\mathcal M$ is an \ofin\ set by Theorem~\ref{T:ulominic}. By \los, $\Gamma$ is the graph of the order reversing permutation $f\colon D\to D$. In particular, $f$ is definable in the \omin\ expansion $(\mathcal M,\Gamma)$ and  maps any element in the class of $l$ to an element in the class of $h$ and vice versa.
 By definability, there is a maximal $a\in D$ such that $f(a)\geq  a$. In particular, $f(a')<a'$, where $a'$ is the successor of $a$ in $D$. A moment's reflection then shows that then either $f(a)=a$ or $f(a)=a'$, which contradicts that no element is $\sim_D$-equivalent with its image.  
\end{proof}

\begin{remark}\label{R:archred}
Similarly, given $D,E\in\discomin M$, we can define $D\ll E$ if for every finite subset $F$, we have $D\cup F\leq E$. If neither $D\ll E$ nor $E\ll D$, then we say that $D$ and $E$ have the same \emph{virtual Archimedean class}, and write $D\sim E$. This is equivalent with the existence of finite subsets $F$ and $G$ such that $D\cup F$ and $E\cup G$ are \viso. The induced order $\ll $ on virtual Archimedean  classes is dense: indeed, suppose $D\ll  E$ and let $d:=\discrk  ED$ and $h:=\discrk  EE$ (i.e., the maximum of $E$).  By Theorem~\ref{T:archred}, since $d\ll h$, there is some $a\in D$ with $d\ll a\ll h$. It follows that $D\ll E_{\leq a}\ll E$. 
\end{remark}

\section{\Taylor\ sets}\label{s:Taylor}
In this section, we work in an   expansion   of $\real$ and its ultrapower $\ul\real$, and we introduce some notation and terminology tailored to this situation. Recall that an element in $\ul\real$ is called \emph{infinitesimal} if its norm is smaller than  $1/n$, for all positive $n$. The \emph{standard part} of  $\alpha\in\ul\real$, denoted $\cp\alpha$, is the supremum of all $r\in\real$ with $r\leq \alpha$; if $\cp\alpha$ is not infinite (that is to say, if $\alpha$ is \emph{bounded}), then $\cp\alpha-\alpha$ is infinitesimal and $\cp\alpha$ is the unique real number with this property.  If $\alpha$ is a tuple $\rij\alpha k$, then we define $\cp\alpha$ coordinate-wise as $(\cp{\alpha_1},\dots,\cp{\alpha_k})$. For a   subset $X\sub\real^k$, we write $\cp X$ for the set of all   $\cp \alpha$ where  $\alpha$ runs over all bounded elements of $\ul X$   (so that $\pm\infty$ never belongs to $\cp X$), and, following the ideology from \cite[\S8]{SchUlBook}, we call $\cp X$ the \emph{catapower} of $X$. We note the following simple result from non-standard analysis:

\begin{lemma}\label{L:cpclos}
The catapower of a subset $X\sub\real^k$ is equal to its closure $\bar X$.
\end{lemma}
\begin{proof}
Suppose $\alpha\in \ul X$ is bounded, given as the ultraproduct of elements $\tuple a_n\in X$. Hence the ultraproduct of the sequence $\cp \alpha-\tuple a_n$ is an infinitesimal, showing that $\tuple a_n$ are arbitrary close to $\cp \alpha$ for almost all $n$. Put differently, there exists a subsequence of $(\tuple a_n)_n$ which converges to $\cp \alpha$, proving that $\cp \alpha$ lies in the closure of $X$. Conversely, if $\tuple b$ lies in closure of $X$, then we can find a sequence $\tuple b_n\in X$ converging to it, and by the same argument, $\cp{(\ul {\tuple b})}=\tuple b$, where $\ul {\tuple b}$ is the ultraproduct of the sequence $\tuple b_n$.
\end{proof}

In  \cite[Chapter 9]{SchUlBook}, we also introduce the notion of  a \emph{protopower}; since it was catered to deal with an additional ring  structure, which is not needed here, we will use only the following  simplified version:   for $X\sub\real^k$, let 
$\trunc Xn$ be the set of points in $X$ whose coordinates have   norm at most $n$, where the norm of a point is defined as the maximum of the absolute values of its coordinates. We define the \emph{protopower}  $\pp\real$ of $\real$ as the ultraproduct of the $\trunc\real n$. We extend this to any subset $X\sub \real^k$, by calling the ultraproduct of the truncations $\trunc Xn$ the    \emph{protopower} of $X$, and denote it $\pp X$. In other words, 
 $\pp X=\ul X\cap \pp\real^k$, where $\ul X$ is the ultrapower of $X$. In particular,  any protopower is bounded (in norm) by $\diag$. (To make this conform with the definitions in \cite[\S9]{SchUlBook}, one actually has to take the protoproduct of the structures $(\real,\frac 1n\norm\cdot)$, and the Archimedean hull of  our $\pp\real$ is then equal to this protoproduct.)

 By the \emph{trace} of a subset $\Xi\sub\ul\real^k$, denoted $\trace\Xi$, we mean the set of its real points, that is to say, $\trace \Xi=\Xi\cap\real^k$. If $\Xi$ is definable by a formula $\varphi$  in some expansion of $\ul\real$, we may use the slightly ambiguous notation $\varphi(\real)$ for its trace as well. The trace    of  a protopower $\pp X$ is equal to $X$, that is to say, $X=\trace{\pp X}$: indeed,   $\tuple a\in\real^k$ satisfies $\tuple a\in \pp X$, \iff\ $\tuple a\in\trunc Xn$ for almost all $n$  (by \los), \iff\  $\tuple a\in X$. 
 For given $n\in\nat$ and a  $k$-ary   function $f$, let us write $\trunc fn$ for the \emph{truncated} function defined by sending a point $\tuple a$ to $f(\tuple a)$ if   $\norm {\tuple a}\leq n$ and to zero otherwise (note that this is not the same as taking the truncation of the graph of $f$, since we allow values of arbitrary high norm).

 Let $\anlang$ be the language of ordered fields together with a function symbol for each everywhere convergent power series (also referred to as a \emph{globally analytic} function). Clearly, we may view $\real$ as an $\anlang$-structure, but this is not very useful, since $\zet$ is definable in it (as the zero set of $\sin(\pi x)$), and therefore neither \tame\ nor \omin. Instead, we approximate this $\anlang$-structure on $\real $ as follows.   
 Let $\anmodel_n$ be the $\anlang$-structure on $\real$
where each function symbol corresponding to a convergent power series $f$ is interpreted as its truncation $\trunc fn$. By \cite{DvdD}, each $\anmodel_n$ is o-minimal (where one usually denotes $\anmodel_1$ by $\real_{\text{an}}$), and hence their ultraproduct $\ul\anmodel$ is \omin. Moreover, $\ul{\anmodel}$ is \tame\ by Corollary~\ref{C:tame}. Although not part of the signature, power series with a smaller radius of convergence can also be encoded in this structure, at least in one variable: using a combination of linear  transformations $x\mapsto ax+b$, and the (inverse) trig functions $\tan x$ and $\arctan x$, any two open intervals (bounded or unbounded) are isomorphic via a globally analytic map. For instance, if $f$ is defined on the open interval $\oo{-1}1$, then $g(x):=f(\frac 2\pi\arctan x)$ is globally analytic, and hence is definable in $\anlang$.

\begin{definition}[\Taylor\ sets]\label{D:Taylor}
We call $X\sub\real^k$ a \emph{\Taylor} set, if there exists an $\anlang$-formula $\varphi(\tuple x,\tuple y)$ (without parameters), such that for each sufficiently large $n$, there exists a tuple of parameters   $\tuple b_n$ so that $\trunc Xn=\varphi(\anmodel_n,\tuple b_n)$.
\end{definition} 

Modifying $\varphi$ if necessary, we may even assume that this holds for all $n$, and that $\norm {\tuple x}\leq n$ is a conjunct in $\varphi$. If $\ul{\tuple b}$ is the ultraproduct of the $\tuple b_n$, then the protopower $\pp X$ is equal to $\varphi(\ul{\anmodel},\ul {\tuple b})$ by \los, and hence  $X=\trace{\pp X}$. 
Any set realized as a  \emph{protopower} of  a \Taylor\ set will be called an \emph{analytic protopower}, giving a one-one correspondence between \Taylor\ sets and analytic protopowers.
We refer to  the defining formula $\varphi(\tuple x,\ul{\tuple b})$ of $\pp X$ as the \emph{analytic} formula for $X$,  and we express this   by writing $X=\varphi(\real)$ (this does not mean that $X$ is definable, since the parameters might be non-standard; in the  terminology of \S\ref{s:locdef}, a \Taylor\ set is in general  only locally $\anlang$-definable). Not every definable subset  is an  analytic   protopower (equivalently, not every $\anlang(\ul\real)$-formula is analytic): let $\Theta$ be defined by $(\exists y)\, xy=1\en \sin(\pi y)=0$. Its trace $\trace\Theta$ is equal to the set of reciprocals of positive natural numbers and cannot be a \Taylor\ set by Lemma~\ref{L:discTaylor} below.  Any quantifier free $\anlang(\ul\real)$-formula is analytic, so that in particular, any globally real analytic variety is \Taylor. \Taylor\ sets are closed under (finite) Boolean combinations, but not under definable (analytic) images, nor under projections. In particular, the \Taylor\ sets do not form a first-order structure.
%
%

\begin{lemma}\label{L:discTaylor}
A real discrete subset   is \Taylor\ \iff\ it is closed. Moreover, a discrete \Taylor\ set intersects any bounded set in finitely many points.
\end{lemma}
\begin{proof}
If $X$ is discrete, then $\trunc Xn$ must be finite by o-minimality, and hence $X$ cannot have an accumulation point whence is closed. Conversely, if $X$ is discrete and closed, then it is the zero set of some analytic function  $f$ (taking sums of squares allows us to reduce  to a single equation), and hence $\trunc Xn$ is defined in $\anmodel_n$ by $\trunc fn(\tuple x)=0$,   and $\norm {\tuple x}\leq n$. 
\end{proof}

\begin{lemma}\label{L:Taylorul}
A subset $X\sub \real^k$ is \Taylor\ \iff\ its protopower $\pp X$ is $\ul{\anmodel}$-definable.
\end{lemma}
\begin{proof}
Recall that $\pp X$ is  the ultraproduct of the truncations $\trunc Xn$. One direction has already been observed.  Assume $\pp X$  is $\ul{\anmodel}$-definable, say $\pp X=\varphi(\ul{\anmodel},\tuple b)$, for some $\anlang$-formula $\varphi$ and some tuple of parameters $\tuple b$. Writing $\tuple b$ as the ultraproduct of tuples $\tuple b_n$, it follows from \los\ that   $\trunc Xn=\varphi(\anmodel_n,\tuple b_n)$ for almost all $n$. Enlarging the tuple of parameters if necessary, we may assume that $n$ is one of the entries of $\tuple b_n$.   Choosing for each $n$ some $m>n$ such that  $\trunc Xm=\varphi(\anmodel_m,\tuple b_m)$, we get  $\trunc Xn=\varphi(\anmodel_n,\tuple b_m)\en   \norm {\tuple x}\leq n$, showing that $X$ is \Taylor.
\end{proof}

We can rephrase this as a criterion for analytic protopowers:
\begin{corollary}\label{C:Taylorul}
A protopower $\pp X\sub \ul\real^k$ is   analytic   \iff\ it is $\ul{\anmodel}$-definable \iff\ its trace $X$ is \Taylor.
 \qed
\end{corollary}
%

In terms of formulae, we might paraphrase this as: an $\anlang(\ul\real)$-formula $\varphi$ is analytic \iff\ $\varphi(\ul\real)$ is the ultraproduct of the $\varphi(\trunc\real n)$. Thus, an open interval in $\ul\real$ is an (analytic) protopower \iff\ its endpoints are either real or equal to $\pm\diag$: indeed,   suppose $\oo\alpha\beta$ is a  protopower, and let $\cp\alpha$ and $\cp\beta$ be the respective standard parts of $\alpha$ and $\beta$. Hence $I:=\oo\alpha\beta\cap \real$ is a (not necessarily open) interval with endpoints $\cp\alpha$ and $\cp\beta$. If $\cp\alpha$ is finite, then $\trunc In$ is an interval with left endpoint $\cp\alpha$ for $n$ sufficiently large, and hence the same is true for the ultraproduct of these truncations. By Corollary~\ref{C:Taylorul}, this forces $\cp\alpha=\alpha$. In the other case, the left endpoint of $\trunc In$ is $-n$, and hence their ultraproduct has left endpoint $-\diag$, showing that $\alpha=-\diag$. The same argument applies to $\beta$, proving the claim.

\begin{example}\label{E:panal}
By Lemma~\ref{L:discTaylor}, every closed, discrete subset, whence in particular any subset of $\zet$, is \Taylor. To give a non-discrete example, consider the spiral $C\sub \real^2$ with parametric equations $x=\exp \tau\sin \tau$ and $y=\exp \tau\cos \tau$, for $\tau\in\real$. If $(x,y)\in\trunc Cn$, then $\exp\tau=\sqrt{x^2+y^2}\leq n\sqrt 2$ and hence $\tau\leq \log(n\sqrt 2)\leq n$. In particular, the negative values of $\tau$ can be larger in absolute value than $n$. Hence $C$ is not \Taylor. However, if $C^+$ is the `positive' part, given by the same equations but only for $\tau\geq 0$, then   $\trunc {C^+}n$ is defined in $\anmodel_n$ by  $x=\trunc \exp n (\tau)\trunc\sin n(\tau)$, $y=\trunc\exp n(\tau)\trunc\cos n(\tau)$, and  $\tau\leq \trunc\log n(n\sqrt 2)$, showing that $C^+$ is \Taylor\ (see Corollary~\ref{C:polarTaylor} below).
%
\end{example}

\begin{proposition}\label{P:closTaylor}
The closure, interior, frontier, and boundary of a \Taylor\ set is again \Taylor.
\end{proposition}
\begin{proof}
Since all concepts are obtained by either taking closures or Boolean combinations, it suffices to show that the closure $\bar X$ of a \Taylor\ set $X$ is again \Taylor. Let $\varphi(x,z)$ be an analytic formula for $X$, so that $\trunc Xn=\varphi(\anmodel_n,\tuple b_n)$, for some parameters $\tuple b_n$ and all $n$. If $\psi(x,z)$ is the formula $(\forall a>0)(\exists y)\norm{x-y}<a\en \varphi(y,z)$,  then $\psi(\anmodel_n,\tuple b_n)$ defines the closure of $\trunc Xn$. It is now easy to check that the latter is equal to $\trunc{\bar X}n$, showing that $\bar X$ is \Taylor.
\end{proof}

\begin{remark}\label{R:closTaylor}
From the proof it is also clear that if $\pp X$ is the  protopower of $X$, then the closure $\overline{\pp X}$ of $\pp X$ is the protopower of $\bar X$, and the analogous properties for the other topological operations. Inspecting the above proofs and examples, we can single out the following geometric feature of \Taylor\ sets.\footnote{The corresponding syntactic characterization of analytic formulae is   not yet clear to me.}
\end{remark}

\begin{proposition}\label{P:quantanalform}
Let $X\sub \real^{k+1}$ be a \Taylor\ set and let $Y\sub \real^k$ be its projection onto the first $k$ coordinates. If    there exists $K\in\nat$ such that  $\trunc Yn$ is contained in the projection of $\trunc X{Kn}$, for all sufficiently large $n$, then $Y$ is again \Taylor. 
%
\end{proposition}
\begin{proof}
Let $\varphi(\tuple x,y,\ul{\tuple c})$ be the analytic formula  defining $X$,  and choose tuples $\tuple c_n$ with  ultraproduct equal to $\ul{\tuple c}$, so that $\trunc Xn$ is defined in $\anmodel_n$ by $\varphi(\tuple x,y,\tuple c_n)$.    Let $\tilde \varphi(\tuple x,y,\ul{\tuple c})$  be the formula obtained from $\varphi$ by replacing every power series  $f(\tuple x,y)$ occurring in it by the power series$f(\tuple x, Ky)$, and put $\psi(\tuple x, \ul{\tuple c}):=(\exists y)\tilde \varphi(\tuple x,y,\ul{\tuple c})$.  I claim that $\psi$ is an analytic formula with $\psi(\real)=Y$. To this end, we have to show that $\trunc Yn=\psi(\anmodel_n,{\tuple c}_n)$, for almost all $n$. One inclusion is clear, so assume $\tuple a\in\trunc Yn$, for some $n$. Hence $\norm{\tuple a}\leq n$ and there exists $b\in\real$ such that $(\tuple a,b)\in X$. By assumption, we can choose $\norm b\leq Kn$. Let $b':=b/K$, so that $\norm{b'}\leq n$. Since then $\anmodel_n\models\tilde\varphi(\tuple a,b')$, as the point $(\tuple a,b')$ has norm at most $n$, whence agrees on any power series with its $n$-th truncation, we get $\anmodel_n\models\psi(\tuple a)$, as required.  
%
%
\end{proof}

Given a $C_1$-function $f\colon\real\to \real$   on an open interval $\oo ab$, we say that $f$ is \emph{increasing} at $b$ if 
$$
\lim_{x\to b^-}f'(x)>0
$$
with a similar definition for decreasing or at the left endpoint.

\begin{corollary}\label{C:polarTaylor}
Let $f$ be  a  power series converging on a half-open interval $\co ab$. If $f$ is  increasing at $b$, then the   curve $C\sub \real^2$  with  polar equation $R=f(\theta)$, for $a\leq\theta<b$, is \Taylor.
\end{corollary}
\begin{proof}
By the discussion at the beginning of this section, we may make an order-preserving, analytic change of variables so that $f$ becomes  convergent on $\real_{\geq0}$. 
By L'H\^ opital's rule, the limit of $f(x)/x$ for $x\to\infty$ exists and is positive. Hence, we may choose $K\in\nat$ large enough so that $1/K<f(x)/x$ for all $x\geq K$.   Let $X\sub\real^3$ be the semi-analytic set given by $x=f(z)\sin(z)$, $y=f(z)\cos(z)$, and $a\leq z<b$, so that $C$ is just the projection of $X$ onto the first two coordinates. By Proposition~\ref{P:quantanalform},  it suffices to show that $\trunc Cn$ is contained in the projection of $\trunc X{(K\sqrt 2)n}$, for all $n$. To this end, let $(a,b)\in\trunc Cn$, so that $a=f(\theta)\sin\theta$ and $b=f(\theta)\cos\theta$, for some $\theta\geq0$. In particular, $f(\theta)= \sqrt {a^2+b^2}\leq n\sqrt 2$. There is nothing to prove if $\theta\leq K$, so let  $\theta>K$ and hence $1/K<f(\theta)/\theta$.  The result now follows since  $\theta< Kf(\theta)\leq (K\sqrt 2)n$.
\end{proof}

Of course, a similar criterion exists if the domain is open at the left endpoint, where the function now has to be decreasing.  Any \Taylor\ set is  of the form $\varphi(\real)$ for some $\anlang(\ul\real)$-formula $\varphi$, that is to say,  is a trace  of an $\ul{\anmodel}$-definable subset. For each such trace $X:=\varphi(\real)$, we can define its \emph{dimension} $\op{dim}(X)$ to be the dimension of $\varphi(\ul\anmodel)$. In general, this notion is not   well behaved: the trace   of the discrete, zero-dimensional set given by the formula $(\exists y>0) \sin(\pi y)=0\en \sin(\pi xy)=0$ is equal to $\mathbb Q$, a non-discrete set. Fortunately,   \Taylor\ sets  behave  tamely, as witnessed, for instance, by the following planar trichotomy (compare with Theorem~\ref{T:triplanar}):

\begin{theorem}\label{T:planardimTaylor}
A non-empty \Taylor\ subset $X\sub \real^2$ is either
\begin{enumerate}
\item  zero-dimensional, discrete, and closed;
\item  one-dimensional, nowhere dense,  but at least one projection has non-empty interior;
\item two-dimensional with non-empty interior.
\end{enumerate}
\end{theorem}
\begin{proof}
Let $\pp X$ be the protopower of $X$ and $d$ its dimension. By Proposition~\ref{P:dimulomin}, almost all truncations $\trunc Xn$ have dimension $d$ . Hence, if $d=0$, then almost all (whence all) $\trunc Xn$ are finite and $X$ is closed and discrete. If $d=2$, then almost all (whence all) $\trunc Xn$ have non-empty interior, whence so does $X$. Finally, if $d=1$,   (almost) all $\trunc Xn$ are nowhere dense, and some projection has interior. Therefore, $X$ itself has the same properties.
\end{proof}

In view of Remark~\ref{R:closTaylor}, the dimension of the frontier $\fr X$ of a \Taylor\ set $X$ is strictly less than  its dimension   $\op{dim}(X)$. Hence, by the same argument as for Corollary~\ref{C:const}, we immediately get:

\begin{corollary}\label{C:constTaylor}
Any \Taylor\ set is constructible. \qed
\end{corollary}

Next, we study maps in this context. For $X\sub \real^k$ and $Y\sub \real^l$, let us call a map $f\colon X\to Y$ \emph\Taylor,    if its graph is a \Taylor\ set.

\begin{corollary}\label{C:defmapTaylor}
The domain and image of a \Taylor\ map are \Taylor, and so is any fiber. Likewise, if the graph of an $\ul{\anmodel}$-definable map $\gamma\colon\Xi\to \Theta$ is a protopower, then so are $\Xi$   and $\gamma(\Xi)$, as well as every fiber $\inverse\gamma {\tuple b}$ with $\tuple b\in Y$. Moreover, the trace of $\gamma$ induces a \Taylor\ map   $g\colon \trace\Xi\to\trace{\gamma(\Xi)}$, and any \Taylor\ map is obtained in this way.  
\end{corollary}
\begin{proof}
The first assertion   follows from the last assertion by  Corollary~\ref{C:Taylorul}. So assume  the graph $\Gamma(\gamma)$ is a  protopower.  Without loss of generality, we may assume $\Theta=\gamma(\Xi)$, that is to say, that $\gamma$ is surjective. 
Let $G:=\trace{\Gamma(\gamma)}$ be the trace, so that the ultraproduct of the $\trunc Gn$ is equal to $\Gamma(\gamma)$, and let $X:=\trace\Xi$ and $Y:=\trace \Theta$ be the respective traces of domain and image. It follows that $\trunc Xn$ is defined by the formula $(\exists \tuple y) (\tuple x,\tuple y)\in \trunc Gn$, and hence $X$ is \Taylor\ (alternatively, use Proposition~\ref{P:quantanalform}). Moreover, $\gamma$ restricted to $\trunc Xn$ takes values inside $\trunc Yn$, that is to say, induces a map $g_n\colon \trunc Xn\to \trunc Yn$. It follows that the ultraproduct of the $g_n$ is equal to $\gamma$. By \los, almost all $g_n$ are surjective. Therefore,   the respective ultraproducts of $\trunc Xn$ and $\trunc Yn$ are $\Xi$ and $\Theta$, proving that both sets are analytic protopowers  by  Corollary~\ref{C:Taylorul}. Moreover, the union $g$ of the $g_n$ is the restriction of $\gamma$ to $X$. Fix $\tuple b\in Y$, let $\Phi:=\inverse\gamma{\tuple b}$ its fiber and $F:=\trace\Phi$ the latter's trace. One checks that $\trunc Fn=\inverse{g_n}{\tuple b}$, and hence the ultraproduct of the $\trunc Fn$ is equal to $\Phi$, proving that $\Phi$ is an analytic protopower. 
\end{proof}

\begin{remark}[\Taylor\ cell decomposition]\label{R:Taylorcell}
In particular, a horizontal \Taylor\ $1$-cell in $\real^2$ must be the graph of a continuous, \Taylor\ map, and similarly,  a \Taylor\ $2$-cell in $\real^2$ is the region between two \Taylor\ graphs. Let $X$ be a \Taylor\ set with protopower $\pp X$. Since $\ul{\anmodel}$ is \tame, we can find a surjective, cellular map $\delta\colon \pp X\to \Delta$ with $\Delta$ a discrete, closed set. I conjecture that we may take $\delta$ to be a  protopower too. Assuming this, taking traces yields a \Taylor\ map $d\colon X\to \trace\Delta$, whose fibers are all \Taylor\ cells, and hence defined by means of continuous \Taylor\ maps. This yields a \Taylor\ cell decomposition of $X$ which is finite on each compact  subset by Lemma~\ref{L:discTaylor}.
\end{remark}

\begin{corollary}\label{C:DPPTaylor}
Any discrete \Taylor\ set $D$ satisfies   DPP in the sense that a \Taylor\ map  $D\to D$ is   injective \iff\ it is surjective. 
\end{corollary}
\begin{proof}
Let $g\colon D\to D$ be \Taylor\ and let $\pp g\colon \pp D\to \pp D$ be its protopower, that is to say, given as the ultraproduct of the restrictions $g_n:=\restrict g{\trunc Dn}$. By \los,  $\pp g$ is injective (surjective) \iff\ almost all $g_n$ are, whence  \iff\ $g$ is, and the result now follows easily from   the \omin\ DPP (Proposition~\ref{P:DPP}). 
\end{proof}

\begin{corollary}[Monotonicity for \Taylor\ maps]\label{C:discTaylor}
A \Taylor\ map $g\colon X\to Y$  is continuous outside a  set of dimension strictly less than the dimension of $X$. In particular,   one-variable \Taylor\ maps are monotone outside a discrete, closed (\Taylor) subset.
\end{corollary}
\begin{proof}
We may assume, for the purposes of this proof  that $g$ is  surjective, so that, in particular, both $X$ and $Y$ are \Taylor, by Corollary~\ref{C:defmapTaylor}. By the same result,  taking protopowers yields a definable map $\pp g\colon \pp X\to \pp Y$ whose restriction to $X$ is equal to $g$. By Theorem~\ref{T:disctu}, the set of discontinuities $\Delta$ of $\pp g$ has dimension strictly less than $\op{dim}(\pp X)=\op{dim}(X)$. Replacing $\Delta$ by its closure,\footnote{In fact, this is not needed since one can  show that $\Delta$ is already closed.} which does not change the dimension, we may assume $\Delta$ is closed. I claim that $g$ is continuous outside the trace $D:=\trace\Delta$. Indeed, if $a\in X\setminus D$, then by the non-standard criterion for continuity, we have to show that for every $\alpha$ infinitesimally close to $a$, their images under $\ul g$ remain infinitesimally close, where $\ul g$ is the ultrapower of $g$. However, since $\pp g$  is the ultrapower of the restrictions $\restrict g{\trunc Xn}$, both maps agree on bounded elements, and so we have to show that $\pp g(a)$ and $\pp g(\alpha)$ are infinitesimally close. This does hold indeed for $\alpha$ sufficiently close to $a$ since $a\notin \Delta$ and $\Delta$ is closed. 

In the one-variable case, we may choose $\Delta$  so that $\pp g$ is monotone on any interval with endpoints in $\Delta$, and clearly, $g$ is then monotone on $D$. It follows from Lemma~\ref{L:discTaylor} and Theorem~\ref{T:planardimTaylor} that $D$ is \Taylor. 
\end{proof}

\begin{remark}\label{R:genctuTaylor}
Using the discussion in Remark~\ref{R:genctu}, we can choose $\Delta$ in the above statement also to be \Taylor\ in higher dimensions. 
\end{remark}

 If $f\colon X\to Y$ is \Taylor\ and bijective, then its inverse is also \Taylor, and we will say that $X$ and $Y$ are  \emph{analytically isomorphic}. In the definition of a \gr, it was not necessary that the collection of subsets formed a first-order structure, only that they were preserved under Boolean combinations and products. Since this is true also of \Taylor\ sets, we can define the \emph{analytic} \gr\ $\grotan$ as the free Abelian group of analytic isomorphism classes of \Taylor\ sets modulo the scissor relations. 

\begin{proposition}\label{P:angr}
 Sending the class of a \Taylor\ set  $X$ to the class of its protopower $\pp X$ induces a natural ring \homo\ of \gr{s} $\grotan\to \grot{\ul\anmodel}$.
\end{proposition}
\begin{proof}
To show that the map $\class X\to \class{\pp X}$ is well-defined, suppose $f\colon X\to Y$ is an analytic isomorphism. The ultraproduct of the truncations $f_n\colon \trunc Xn\to \trunc Yn$ induces then a definable map $\pp f\colon  \pp X\to \pp Y$, and by \los, this is again a bijection.
\end{proof}
%

\providecommand{\bysame}{\leavevmode\hbox to3em{\hrulefill}\thinspace}
\providecommand{\MR}{\relax\ifhmode\unskip\space\fi MR }
\providecommand{\MRhref}[2]{%
  \href{http://www.ams.org/mathscinet-getitem?mr=#1}{#2}
}
\providecommand{\href}[2]{#2}


\end{document}